\documentclass{compositio}

\renewcommand{\contentsname}{Contents}

\usepackage{amsmath} 
\usepackage{amssymb}
\usepackage{amsthm}
\usepackage{mathrsfs} 
\usepackage{epsfig}
\usepackage{appendix}

\usepackage{caption} 

\usepackage{graphicx}
\usepackage{float}
\usepackage{amssymb}
\usepackage{mathrsfs}
\usepackage{epstopdf}
\usepackage[ruled, linesnumbered, vlined,commentsnumbered, algosection]{algorithm2e}
\SetAlgorithmName{Pseudocode}{Pseudocode}{Pseudocode}
\usepackage{algorithmic}
\usepackage{float}

\usepackage{bm}
\usepackage{color}
\usepackage{tikz-cd}
\usepackage{quiver}
\usepackage{capt-of}

\usepackage{matlab-prettifier}
\usepackage{tcolorbox}
\usepackage{listings}
\usepackage{xparse}
\tcbuselibrary{listings,skins,breakable,xparse}
\tcbset{%
 mylist/.style = {%
  colframe = gray,
  colback = white,
  coltitle = red!50!yellow!3!white,
  colbacktitle = white,
  listing only,
  attach boxed title to top center = {yshift = -\tcboxedtitleheight/2},
  enhanced,
  drop fuzzy shadow, 
  left = 6.5mm, 
  breakable, 
  enhanced jigsaw, 
  fonttitle=\small\bfseries\color{black}, 
  before skip=20pt plus 2pt, 
  after skip=20pt plus 2pt,
 },
 example/.style 2 args = {%
  mylist,
  title = {\thetcbcounter: #1},
  label = {#2},
 },
}
\newtcblisting[auto counter, number within = section, list inside = mcode]{matlab}[3][]{%
 listing options = {%
  style = MATLAB-editor,
  numbers = left,
  numberstyle = \footnotesize\color{darkgray}
 },
 overlay = {\begin{tcbclipinterior}\fill[blue!15!white] (frame.south west) rectangle ([xshift=5.3mm]frame.north west);\end{tcbclipinterior}},
 example = {#2}{#3}, #1,
}

\newtheorem{thm}{Theorem}[section] 
\newtheorem{cor}[thm]{Corollary}

\newtheorem{prop}[thm]{Proposition}
\newtheorem{pty}[thm]{Property}
\newtheorem{lem}[thm]{Lemma}

\theoremstyle{definition} 
\newtheorem{DEF}[thm]{Definition} 
\newtheorem{conj}[thm]{Conjecture}

\newtheorem{alg}[thm]{Algorithm}

\theoremstyle{remark} 
\newtheorem{rmk}[thm]{Remark}

\newtheorem{que}[thm]{Question}
\newtheorem{eg}[thm]{Example}

\usepackage{cite}
\usepackage[colorlinks,linkcolor=blue]{hyperref}

 \usepackage[usenames,dvipsnames]{pstricks}
 \usepackage{pstricks-add}
 \usepackage{epsfig}
 \usepackage{pst-grad} 
 \usepackage{pst-plot} 
 \usepackage[space]{grffile} 
 \usepackage{etoolbox} 
 \makeatletter 
 \patchcmd\Gread@eps{\@inputcheck#1 }{\@inputcheck"#1"\relax}{}{}
 \makeatother

\title[The approach of cluster symmetry to Diophantine equations]{The approach of cluster symmetry to\\ Diophantine equations}

\date{\today}

\begin{document}

\author{Leizhen Bao}
\email{mathbao@zju.edu.cn}
\address{1. School of Mathematical Sciences, 
Zhejiang University,  
Hangzhou, Zhejiang 310058, 
PR China\\
2. Department of Mathematical Sciences, 
Tsinghua University,  
Beijing 100084, 
PR China}
%
\author{Fang Li}
\email{fangli@zju.edu.cn}
\address{1. School of Mathematical Sciences, 
Zhejiang University,  
Hangzhou, Zhejiang 310058, 
PR China}

\classification{11D72 (primary), 13F60, 13A50 (secondary).}
\keywords{Diophantine equation, (generalized) cluster algebra, mutation, cluster symmetry, invariant ring.}

\begin{abstract}
This paper aims to employ a cluster-theoretic approach to provide a class of Diophantine equations whose solutions can be obtained by starting from initial solutions through mutations.

We establish a novel framework bridging cluster theory and Diophantine equations through the lens of cluster symmetry.  
On the one hand, we give the necessary and sufficient condition for Laurent polynomials to remain invariant under a given cluster symmetric map.
On the other hand, we construct a discriminant algorithm to determine whether a given Laurent polynomial has cluster symmetry and whether it can be realized in a generalized cluster algebra. As applications, we solve Markov-cluster equations, describe some invariant Laurent polynomial rings, and resolve the questions posed by Gyoda and Matsushita.
\end{abstract}

\maketitle

\vspace*{6pt}\tableofcontents

\section{Introduction}

In 1880, when Markov \cite{Markov} working on a Diophantine approximation problem\cite{Markov100},
he needed to consider the following Diophantine equation (\textbf{Markov equation})
    \begin{equation}\label{eq markov}
        M(x_1,x_2,x_3) := {x_1^2 + x_2^2 + x_3^2} - 3{x_1x_2x_3} = 0.  
    \end{equation}
To do this, he defined three transformations
    \begin{eqnarray*}
        m_1(x_1,x_2,x_3) &:=& (3x_2x_3 - x_1, x_2, x_3),\\
        m_2(x_1,x_2,x_3) &:=& (x_1, 3x_1x_3-x_2, x_3),\\
        m_3(x_1,x_2,x_3) &:=& (x_1, x_2, 3x_1x_2 - x_3).
    \end{eqnarray*}
He found that the orbits of the initial solution $(1,1,1)$ under the group $\langle m_1, m_2, m_3 \rangle$ are exactly the positive integer solutions of Equation \eqref{eq markov}. These transformations are important, but no one knew what they meant then. 

In 2002, Fomin and Zelenvinsky \cite{FZ1}, working on canonical bases of quantum groups,
abstracted out the so-called \textbf{cluster algebra}, which is also combinatorially called a \textbf{cluster pattern}. A cluster algebra can be defined by a skew-symmetric matrix $B$, a tuple of variables $\mathbf{x}$, and a set of transformations $\mu_1, \dots, \mu_n$ called \textbf{mutations}. The notion of mutations is the key to cluster theory. In fact, mutations can be realized on many mathematical subjects, such as flips of triangulations in Riemann surfaces \cite{FominShapiroThurston}, Bongartz completions of tilting modules\cite{Bongartz}, 
wall-crossing automorphisms of scattering diagrams\cite{ GrossHackingKeelKontsevich}, and so on. 
The theory of cluster algebras is widely associated with many related fields, such as quantum dilogarithms\cite{Keller2011}, Poisson geometry\cite{Gekhtman2010}, Donaldson-Thomas invariant theories\cite{DTandCA}, mirror symmetry theories\cite{Gross2015}, and other theories.

In 2012, Peng and Zhang\cite{PengZhang} revealed the connection between the Markov equation and a cluster algebra. Considering the skew-symmetric matrix
\begin{equation*} B := \left[\begin{array}{rrr}
  0 & 2 & -2 \\
  -2 & 0 & 2 \\
  2 & -2 & 0
\end{array} \right],
\end{equation*}
the corresponding mutations are
\begin{eqnarray*}
    \mu_1(x_1, x_2, x_3) &:=& \bigg(\frac{x_2^2+x_3^2}{x_1}, x_2, x_3\bigg),\\ 
    \mu_2(x_1, x_2, x_3) &:=& \bigg(x_1, \frac{x_1^2+x_3^2}{x_2}, x_3\bigg),\\
    \mu_3(x_1, x_2, x_3) &:=& \bigg(x_1, x_2,\frac{x_1^2+x_2^2}{x_3}\bigg).
\end{eqnarray*}
It is easy to check that if $(a,b,c)$ is a positive solution of Equation \eqref{eq markov}, then $m_i(a,b,c) = \mu_i(a,b,c)$ for any $i = 1,2,3$.
Then cluster algebra has had a connection with number theory.

However, the relationship between cluster algebra and number theory is not limited to this. In 2016, Lampe \cite{Lampe} connected a cluster algebra of rank $3$ to a Diophantine equation and connected a cluster algebra of rank $5$ to a Laurent polynomial. In 2024, Chen and Li\cite{ChenLi} considered cluster algebras of rank $2$. They constructed the corresponding Diophantine equations, found positive integer solutions to these equations, and classified the Diophantine equations for the cluster algebras of rank $2$.

There are still some other algebraic structures in cluster theory that can be associated with Diophantine equations.
In 2014, Chekhov and Shapiro\cite{ChekhovShapiro2014} generalized the cluster algebra into \textbf{generalized cluster algebra}; later in 2024, Gyoda and Matsushita\cite{Gyoda} considered some generalized cluster algebras of rank $3$, they constructed the corresponding Diophantine equations and solved the positive integer solutions of these equations. 
In 2009, Fock and Goncharov \cite{FockGoncharov} generalized cluster algebras into \textbf{cluster ensembles}; later in 2024, Kaufman \cite{Kaufman} considered cluster ensembles of affine ADE type, constructed the corresponding invariant Laurent polynomials for a composition of a permutation $\sigma_{(12)}$ and a mutation $\mu'_1$ and described the structure of the corresponding invariant rational function field.

In addition to this, some discrete dynamical systems can be related to Diophantine equations. For example, in 2008, Hone and Swart \cite{HoneSwart} considered the following two recurrence relations
\[
    u_{n+4}u_n = \alpha u_{n+3}u_{n+1} + \beta (u_{n+2})^2,\;\; 
    v_{n+5}v_n =\tilde{\alpha} v_{n+4}v_2 + \tilde{\beta}v_{n+2}v_{n+3}.
\]
On the one hand, these two recurrence relations generate Somos $4$ sequence $\{u_n\}$ and Somos $5$ sequence $\{v_n\}$, which are related to certain elliptic curves. On the other hand, according to the above recurrence relations, two mappings are defined as
\begin{eqnarray*}
    \psi_4(u_1, u_2, u_3, u_4) &:=& \bigg(u_2, u_3, u_4, \frac{\alpha u_4u_2 + \beta u_3^2}{u_1}\bigg),\\
    \psi_5(v_{1}, v_2, v_3, v_4, v_5) &:=& \bigg(v_2, v_3, v_4, v_5, \frac{\tilde{\alpha} v_5v_2 + \tilde{\beta}v_3v_4}{v_{1}}\bigg).
\end{eqnarray*}
They constructed two Laurent polynomials that are invariant under these maps
\begin{eqnarray*}
    F_4(x_1, x_2, x_3, x_4) &:=& \frac{x_1^2 x_4^2+\alpha(x_1 x_3^3+ x_2^3x_4)+\beta x_2^2 x_3^2}{x_1x_2x_3x_4}, \label{eqn somos4}\\
    F_5(x_1, x_2, x_3, x_4, x_5) &:=& \frac{x_1x_2^2x_5^2 + x_1^2x_4^2x_5 +\tilde{\alpha}(x_1x_3^2x_4^2  + x_2^2x_3^2x_5) + \tilde{\beta}x_2x_3^3x_4}{x_1x_2x_3x_4x_5}. \label{eqn somos5}
\end{eqnarray*}
That is, $F_4(\psi_4(\mathbf{x})) = F_4(\mathbf{x}), F_5(\psi_5(\mathbf{x})) = F_5(\mathbf{x})$. In addition, when the parameters $\alpha = \beta = \tilde{\alpha} = \tilde{\beta} = 1$, for this special case, the map $\psi_4$ is the composition of a permutation $\sigma_{(1234)}$ with the mutation $\mu^{(4)}_1$ in some cluster algebra of rank $4$, that is, $\psi_4 = \sigma_{(1234)}\mu^{(4)}_1$; the map $\psi_5$ is the composition of a permutation $\sigma_{(12345)}$ with the mutation $\mu^{(5)}_1$ in some cluster algebra of rank $5$, that is, $\psi_4 = \sigma_{(12345)}\mu^{(5)}_1$.

In summary, Diophantine equations related to cluster theory have frequently been found in recent years. Naturally, a question arises:

\textbf{Is there a systematic method to find Diophantine equations related to cluster theory?}

We note that the Diophantine equations related to cluster theory in the papers\cite{Markov, Lampe, Bao1, ChenLi, Gyoda, Kaufman, HoneSwart} are certain Laurent polynomials with initial vectors. For example, define the Laurent polynomial
    \begin{equation}\label{eq Laurent Markov}
        F_1(x_1,x_2,x_3) := \frac{x_1^2 + x_2^2 + x_3^2}{x_1x_2x_3},
    \end{equation}
the positive integer solutions of the Markov equation are the same as the positive integer solutions of the equation $F_1(x_1, x_2, x_3) = F_1(1,1,1)$.
Therefore, we focus on finding these Laurent polynomials related to cluster theory. 

We first classify Laurent polynomials. A Laurent polynomial $F(\mathbf{x})$ is of $\frac{\bm{\eta}}{\mathbf{d}}$ \textbf{type}, if $F(\mathbf{x}) = \frac{T(\mathbf{x})}{\mathbf{x}^{\mathbf{d}}}$, where $T(\mathbf{x}) \in \mathbb{Q}[\mathbf{x}]$, $\mathbf{d} \in \mathbb{Z}^n$, ${\eta_i}$ is the degree of $x_i$ in $T(\mathbf{x})$ and $x_i \nmid T(\mathbf{x})$ for all $i \in [1,n]$. For example, the Laurent polynomial $F_1(x,y,z)$ is of type $\frac{(2,2,2)}{(1,1,1)}$, $F_4(x_1, x_2, x_3, x_4)$ is of type $\frac{(2,3,3,2)}{(1,1,1,1)}$, and $F_5(x_1, x_2, x_3, x_4, x_5)$ is of type $\frac{(2,3,3,3,2)}{(1,1,1,1,1)}$.

These Laurent polynomials have the property that they are invariant under some special transformations related to cluster theory. For example, $F_1(\mu_i(\mathbf{x})) = F_1(\mathbf{x})$, $F_4(\psi_4(\mathbf{x})) = F_4(\mathbf{x})$, $F_5(\psi_5(\mathbf{x})) = F_5(\mathbf{x})$. From this, we introduce the \textbf{cluster symmetric map} $\psi_{\sigma, s, \omega_s}$ of the data $(\sigma, s, \omega_s)$ as 
\begin{equation*} 
    \psi_{\sigma, s, \omega_s}(\mathbf{x}) := \bigg(x_{\sigma(1)}, \cdots, x_{\sigma(t-1)}, \frac{\mathbf{x}^{r[-\mathbf{b}]_{+}}Z(\mathbf{x}^{\mathbf{b}})}{x_s}, x_{\sigma(t+1)}, \cdots, x_{\sigma(n)}\bigg),
\end{equation*}
where the meaning of the notations can be seen in  Definition \ref{seedlet}.

As Markov did, if a Laurent polynomial is invariant under a cluster symmetric map, then we can get a new solution from the initial one by applying the cluster symmetric map. Thus, we turn to the following question:

\textbf{How to find a Laurent polynomial that is invariant under a given cluster symmetric map?}

To this question, we give an affirmative answer in this paper: for a Laurent polynomial that is invariant under the action of a cluster symmetric map, we provide sufficient and necessary conditions to be satisfied by its coefficients.

\begin{thm}[(Theorem \ref{thm len1_T2f} and \ref{thm len1}, Remark \ref{rmk HLE})]
    Given a cluster symmetric map $\psi_{\sigma, s, \omega_s}$. Let  $F(\mathbf{x})$ be a Laurent polynomial of type $\frac{\bm{\eta}}{\mathbf{d}}$ in  $\mathbb{Q}[\mathbf{x}^{\pm}]$ and its expansion is 
        $$F(\mathbf{x}) = {\mathbf{x}^{-\mathbf{d}}}\sum_{\mathbf{j} \in \mathcal{N}} a_{\mathbf{j}}\mathbf{x}^{\mathbf{j}},$$
    where $\bm{\eta} \in \mathbb{Z}_{\geq 0}^n, \mathbf{d} \in \mathbb{Z}^n$ with $\mathbf{d} = \sigma{(\mathbf{d})}$ and $\eta_s = \eta_{\sigma^{-1}(s)} = 2d_s = 2d_{\sigma^{-1}(s)}$.
    \begin{equation} \label{eqn Fpsi=F}
        F(\psi_{\sigma, s, \omega_s}(\mathbf{x})) = F(\mathbf{x}),
    \end{equation}
    holds, if and only if, the coefficients $\{a_{\mathbf{j}} \in  \mathbb{Q} \mid \mathbf{j} \in \mathcal{N} \}$ of the Laurent polynomial $F(\mathbf{x})$ satisfy the system of homogeneous linear equations $HLE(\sigma, s, \omega_s, \bm{\eta}, \mathbf{d})$ defined in Remark \ref{rmk HLE}.
\end{thm}
\begin{rmk}
    When conditions $\mathbf{d} = \sigma{(\mathbf{d})}$ and $\eta_s = \eta_{\sigma^{-1}(s)} = 2d_s = 2d_{\sigma^{-1}(s)}$ are not satisfied, Relation \eqref{eqn Fpsi=F} also does not hold.
\end{rmk}

To solve the system of homogeneous linear equations $HLE(\sigma, s, \omega_s, \bm{\eta}, \mathbf{d})$ in the above theorem, we wrtie a MATLAB program attached to Appendix \ref{code ca2diop}, so that we can construct an invariant Laurent polynomial of the cluster symmetric maps efficiently and conveniently.

Then, we consider the opposite question: 

\textbf{How to find a cluster symmetric map such that a given Laurent polynomial is invariant under the map?} 

To do this, we collect all cluster symmetric maps of a given Laurent polynomial into a set. For a Laurent polynomial $F(\mathbf{x})$ of type $\frac{\bm{\eta}}{\mathbf{d}}$, the \textbf{cluster symmetric set} $\mathcal{S}(F)$ of $F(\mathbf{x})$ is defined as 
$\mathcal{S}(F) := \{ \psi_{\sigma,s,\omega_s} \mid F(\psi_{\sigma,s,\omega_s}(\mathbf{x})) = F(\mathbf{x}), \eta_s \ne 0\}$. Using an algorithm, we can determine this set.

\begin{thm}[(Algorithm \ref{algo diop2cs}, Proposition \ref{prop S(F)})] 
    Given a Laurent polynomial $F(\mathbf{x}) \in \mathbb{Q}[\mathbf{x}^{\pm}]$.
    The cluster symmetric set $\mathcal{S}(F)$ of $F(\mathbf{x})$ can be obtained by Algorithm \ref{algo diop2cs}.
\end{thm}

We provide a MATLAB program for this algorithm, attached to Appendix \ref{code diop2ca}.

From the above theorem, we can determine whether a Laurent polynomial corresponds to a seed or a generalized cluster algebra.

\begin{prop}[(Definition \ref{def seed and 1-csm}, \ref{def S1 of F}, \ref{def 1-cs seed}, Proposition \ref{prop find 1-csSeed})]
    Given a Laurent polynomial $F(\mathbf{x}) \in \mathbb{Q}[\mathbf{x}^{\pm}]$. Suppose that the cluster symmetric set $\mathcal{S}(F)$ is nonempty. If there exists a seed $\Omega$, such that for any $\psi_{\sigma,s,\omega_s} \in \mathcal{S}(F)$, the relations $\sigma\mu_s \in \mathcal{S}(\Omega)$ and $\omega_s = \pi_s(\Omega^{\pm})$ hold, then $\mathcal{S}(F) \subseteq \mathcal{S}(\Omega)$.
\end{prop}

As an application of our theoretical framework, we show that several other Diophantine equations share the same solution structure as the Markov equation \eqref{eq markov}.

\begin{thm}[(Theorem \ref{thm markov-cluster sol}, Definition \ref{def mc eqn})]
    For $i \in [1,10]$. Let $G_{3,i}$ be the group generated by the subset $\{ \mu_1, \mu_2, \mu_3 \}$ of the cluster symmetric set $\mathcal{S}(\Omega_{3,i})$. Then the set of positive integer solutions of the Markov-cluster equation $F_{3,i}(\mathbf{x}) = F_{3,i}(\mathbf{1})$ is exactly the orbit $G_{3,i}(\mathbf{1})$, that is,
    \begin{equation*}
        G_{3,i}(1,1,1) = \mathcal{V}_{\mathbb{Z}_{>0}}(F_{3,i}(x, y, z) - F_{3,i}(1,1,1)),
    \end{equation*}
    where the seeds $\Omega_{3,i}$ and the Laurent polynomials $ F_{3,i}$ are listed in Table \ref{tab all}.
\end{thm}

These seeds $\Omega_{3,1},\dots, \Omega_{3,10}$ share the same properties, that is, rank 3, the irreducibility of exchange matrices, and $\mu_1,\mu_2,\mu_3\in \mathcal{S} (\Omega_{3,i})$. 
Are there any other seeds satisfying these properties that can correspond to non-constant Laurent polynomials? The answer is negative.

\begin{thm}[(Corollary \ref{cor rank3 nonconstant}, Remark \ref{rmk irr B})]
    For any seed $\Omega := (B, \mathbf{x}, R, \mathbf{Z})$ of rank $n=3$ with an irreducible exchange matrix $B$. Suppose $\mu_1, \mu_2, \mu_3 \in \mathcal{S}(\Omega)$. 
    Then the relation $\mathbb{Q}[\mathbf{x}^{\pm}]^{\langle\mu_1, \mu_2, \mu_3\rangle} \ne \mathbb{Q}$ holds, if and only if, $\Omega = \sigma(\Omega_{3,i})$ for some $i \in [1,10]$, $\sigma \in \mathfrak{S}_3$, where $\Omega_{3,i}$ listed in Table \ref{tab all}.
\end{thm}

This paper is organized as follows. 

In Section \ref{sec 1-csm}, we define the cluster symmetric map (Definition \ref{seedlet}) and prove the main conclusion (Theorem \ref{thm len1}), which gives a method for constructing Laurent polynomials that are invariant under a given cluster symmetric map. 

In Section \ref{sec: ca}, we recall some definitions and some results about generalized cluster algebra, discuss the relationship between generalized cluster algebras and cluster symmetric maps (Proposition \ref{is CSM}), describe the invariant Laurent polynomial ring for some special cases (Proposition \ref{prop len0 B=0}, Theorem \ref{thm 4bounded}), and answer two questions posed by Gyoda and Matsushita in \cite{Gyoda} (Proposition \ref{poly gyoda}, Proposition \ref{rank3 nonconstant}). 

In Section \ref{sec diop}, we discuss solutions to Diophantine equations related to cluster theory. We prove that the Markov-cluster equations possess the same solution structure as the Markov equation, that is, its set of positive integer solutions coincides exactly with a group orbit (Theorem \ref{thm markov-cluster sol}).

In Section \ref{sec CSP}, we define cluster symmetric set of a Laurent polynomial (Definition \ref{def S1 of F}), and give an algorithm to find cluster symmetry of a Laurent polynomial (Algorithm \ref{algo diop2cs}). We then determine when a Laurent polynomial corresponds to a generalized cluster algebra (Proposition \ref{prop find 1-csSeed}). As a summary, we give Figure \ref{fig: CAN} which shows the relationships between the main concepts and theorems throughout this paper.

In the appendices, we show two MATLAB programs related to the main theorems of this paper; the program in Appendix \ref{code ca2diop} constructs invariant Laurent polynomials for a given cluster symmetric map, and the program in Appendix \ref{code diop2ca} finds all non-trivial cluster symmetric pairs of a given Laurent polynomial.\\

For convenience, we use the following notation. 

Fix a positive integer $n$. Let $S$ and $S'$ be sets of $n$-tuples, $\mathbf{v}$ be a $n$-tuple, $\sigma$ be a permutation in the symmetric group $\mathfrak{S}_n$. We denote three sets

$S + S' := \{ \mathbf{a} + \mathbf{b} \mid \mathbf{a} \in S, \mathbf{b}\in S'\}$, 

$S + \mathbf{v} := \{ \mathbf{a} + \mathbf{v} \mid \mathbf{a} \in S\}$, 

$\sigma(S) := \{\sigma(\mathbf{j}) \mid \mathbf{j} \in S \}.$

Denote $\pi_k(\mathbf{v})$ as the $k$-th component of the $n$-tuple $\mathbf{v}$. For any integer $i$, we denote the subset $\pi_k^{(i)}(S) := \{\mathbf{j} \in S ~|~ \pi_k(\mathbf{j}) = i\}$.

We denote that $\mathbf{x}:=(x_1, \dots, x_n)$  a tuple of $n$ indeterminates $x_1, \dots, x_n$, $\mathbf{x}^{\mathbf{v}}:= x_1^{v_1}\cdots x_n^{v_n}$, $\mathbb{Q}[\mathbf{x}] := \mathbb{Q}[x_1, \cdots, x_n]$ the polynomial ring and $\mathbb{Q}[\mathbf{x}^{\pm}] := \mathbb{Q}[x_1^{\pm}, \cdots, x_n^{\pm}]$  the Laurent polynomial ring. The invariant Laurent ring of a given group $G$ is
    $$\mathbb{Q}[\mathbf{x}^{\pm}]^G := \{F(\mathbf{x})\in \mathbb{Q}[\mathbf{x}^{\pm}] \mid F(g(\mathbf{x})) = F(\mathbf{x}), \text { for all } g \in G\}.$$

\section{Invariant Laurent polynomials of cluster symmetric maps} \label{sec 1-csm}

In this section, we will discuss, for a given cluster symmetric map, how to construct a Laurent polynomial that is invariant under the map. 

We first introduce the cluster symmetric map in Subsection \ref{setting}. Then we prove the main theorem, Theorem \ref{thm len1}, in Subsection \ref{CA2diop}. Last, we apply the main theorem for some examples in the Subsection \ref{app: ca2diop}.

\subsection{Cluster symmetric maps of datum}\label{setting}

We first define the cluster symmetric map of the data.

\begin{DEF} \label{seedlet}
Fix a positive integer $n$.

(i) For $s \in [1,n]$. A \textbf{seedlet at direction} $s$ is a triplet $\omega_s := (\mathbf{b}, r, Z)$, where
    \begin{itemize}
        \item[$\bullet$] $\mathbf{b} = (b_1, \cdots, b_n)$ is an $n$-tuple integer vector with $b_s = 0$;
        \item[$\bullet$] $r$ is an positive integer;
        \item[$\bullet$] $Z(u) = \sum_{i=0}^{r}z_iu^i \in \mathbb{Z}_{\geq 0}[u] $ is a polynomial satisfying 
            \begin{equation} \label{ne 0}
                z_0, z_r > 0
            \end{equation}
    \end{itemize}

(ii) The \textbf{exchange polynomial of the seedlet } $\omega_s$ is a polynomial $P_{\omega_s} \in \mathbb{Z}_{\geq 0}[\mathbf{x}]$ defined as
    \begin{equation}
        P_{\omega_s}(\mathbf{x})  := \mathbf{x}^{r[-\mathbf{b}]_{+}}Z(\mathbf{x}^{\mathbf{b}}) =  \sum_{i=0}^{r} z_{i}\mathbf{x}^{i[\mathbf{b}]_{+} + (r - i)[-\mathbf{b}]_{+}}, \label{Pk=sum ljy}
    \end{equation}
where $[\mathbf{b}]_+ := ([b_1]_+, \cdots, [b_n]_+)$ and $[b_i]_+ := \max\{b_i,0\}$.

(iii) For a permutation $\sigma \in \mathfrak{S}_n$, $s\in [1,n]$ and a seedlet $\omega_s$. We call $(\sigma, s, \omega_s)$ is a \textbf{data}. Let $t = \sigma^{-1}(s)$. The \textbf{cluster symmetric map of the data} $(\sigma, s, \omega_s)$ is defined as
    \begin{equation} \label{CSM}
        \psi_{\sigma, s, \omega_s}(\mathbf{x}) := \bigg(x_{\sigma(1)}, \cdots, x_{\sigma(t-1)}, \frac{P_{\omega_s}(\mathbf{x})}{x_s} , x_{\sigma(t+1)}, \cdots, x_{\sigma(n)}\bigg).
    \end{equation}
Briefly, the map is called \textbf{cluster symmetric map}.
\end{DEF}
\begin{rmk}
    For any $\sigma, \tau \in \mathfrak{S}_n$, we denote $\sigma(\mathbf{x}) := (x_{\sigma(1)}, \cdots, x_{\sigma(n)})$. Note that
    \begin{equation} \label{sigmatau}
        \sigma\tau(\mathbf{x}) = (x_{\sigma\tau(1)}, \cdots, x_{\sigma\tau(n)}) = \tau(x_{\sigma(1)}, \cdots, x_{\sigma(n)}) = \tau(\sigma(\mathbf{x})).
    \end{equation}
    Then Equation \eqref{CSM} can be written briefly as 
    \begin{equation}\label{CSM'}
    \psi_{\sigma, s, \omega_s}(\mathbf{x}) = \bigg(\sigma(\mathbf{x})\bigg)\bigg|_{x_s \gets \frac{P_{\omega_s}(\mathbf{x})}{x_s}}.
    \end{equation}
\end{rmk}

For readers who are familiar with generalized cluster algebra, the seedlet $\omega_s$ covers the information in the direction $s$ of a seed in the generalized cluster algebra, with the additional difference that the polynomial $Z(u)$ here does not require the reciprocity condition: $z_i = z_{r-i} $ for all $i \in \{0, \cdots, r \}$. The cluster symmetric map can be viewed as a composite map of a permutation and a mutation; the maps shown in \cite{Markov, Lampe, Bao1, ChenLi, Gyoda, Kaufman, HoneSwart} are all cluster symmetric maps. However, not all composite maps are cluster symmetric maps. In the next section, Proposition \ref{is CSM} will discuss it.

We begin with some examples when $\sigma =id$.

\begin{eg} \label{eg: 1csm, Gyoda} 
(i) Given a seedlet $\omega_2 := (\mathbf{b}', 1, Z)$, where $\mathbf{b}' = (1, 0, -2)$ and $Z(u) = 1 + u$. Then the cluster symmetric map of $(id, 2, \omega_2)$ is
    $$\psi_{id, 2, \omega_2}(\mathbf{x}) =  \bigg(x_1, \frac{x_1 + x_3^2}{x_2}, x_3\bigg).$$

(ii) Given a seedlet $\omega_3 := (\mathbf{b}'', 1, Z)$ , where $\mathbf{b}'' = (-1, 2, 0)$ and $Z(u) = 1 + u$. Then the cluster symmetric map of $(id, 3, \omega_3)$ is
    $$\psi_{id, 3, \omega_3}(\mathbf{x}) =  \bigg(x_1, x_2, \frac{x_1 + x_2^2}{x_3}\bigg).$$

(iii) Given a seedlet $\omega_1 := (\mathbf{b}, 4, Z)$, where $\mathbf{b} = (0,-1,1)$ and $Z(u) = k_0 + k_1u + k_2u^2 + k_3u^3 + k_4u^4$. Then the cluster symmetric map of $(id, 1, \omega_1)$ is
    $$\psi_{id, 1, \omega_1}(\mathbf{x}) =  \bigg( \frac{k_0x_2^4+k_{1}x_2^3x_3+k_{2}x_2^2x_3^2+k_{3}x_2x_3^3+k_4x_3^4}{x_1}, x_2, x_3\bigg).$$
\end{eg}

We then give some examples when $\sigma \ne id$.
\begin{eg} \label{eg CSM somos and marsh}
    (i) Given a seedlet $\omega_1 := (\mathbf{b}, 1, Z)$,where $\mathbf{b} = (0,1,1)$ and $Z(u) = 1+ u$. Then the cluster symmetric map of $(\sigma_{(12)}, 1, \omega_1)$ is
    $$\psi_{\sigma_{(12)}, 1, \omega_1}(\mathbf{x}) = \bigg(x_2, \frac{1 + x_2x_3}{x_1}, x_3\bigg).$$ 
    
    (ii) Given a seedlet $\omega_1 := (\mathbf{b}, 1, Z)$,where $\mathbf{b} = (0,1,1)$ and $Z(u) = 1+ u$. Then the cluster symmetric map of $(\sigma_{(123)}, 1, \omega_1)$ is
    $$\psi_{\sigma_{(123)}, 1, \omega_1}(\mathbf{x}) = \bigg(x_2, x_3, \frac{1 + x_2x_3}{x_1}\bigg).$$ 
    This map was studied by Fordy and Marsh in \cite{FordyMarsh} and is related to the primitive period 1 quiver.
    
    (iii) Given a seedlet $\omega_1 := (\mathbf{b}, 1, Z)$,where $\mathbf{b} = (0,1,-2,1)$ and $Z(u) = {\beta} + {\alpha}u$. Then the cluster symmetric map of $(\sigma_{(1234)}, 1, \omega_1)$ is 
    $$\psi_{\sigma_{(1234)}, 1, \omega_1}(\mathbf{x}) = \bigg(x_2, x_3, x_4, \frac{{\alpha}x_2x_4 + {\beta}x_3^2}{x_1}\bigg).$$
    This map was studied by Hone and Swart in \cite{HoneSwart} and is related to the Somos 4 sequence.
    
    (iv) Given a seedlet $\omega_1 := (\mathbf{b}, 1, Z)$, where $\mathbf{b} = (0,1,-1,-1,1)$ and $Z(u) = \tilde{\beta} + \tilde{\alpha}u$. Then the cluster symmetric map of $(\sigma_{(12345)}, 1, \omega_1)$ is 
    $$\psi_{\sigma_{(12345)}, 1, \omega_1}(\mathbf{x}) = \bigg(x_2, x_3, x_4, x_5,  \frac{\tilde{\alpha}x_2x_5 + \tilde{\beta}x_3x_4}{x_1}\bigg).$$
    This map was studied by Hone in \cite{HoneSigma} and is related to the Somos 5 sequence.
\end{eg}

We show some properties about cluster symmetric maps.

\begin{pty}\label{pty ex-poly} 
For $\sigma \in \mathfrak{S}_n$, $s\in [1,n]$ and a seedlet $\omega_s := (\mathbf{b}, r, Z)$. Let $\psi_{\sigma, s, \omega_s}(\mathbf{x})$ be a cluster symmetric map.

(i) Let $\omega_s' := (-\mathbf{b}, r, Z')$, where $Z'(u) = u^rZ(1/u)$. Then
    \begin{equation} \label{psi b = psi -b}
      \psi_{\sigma, s, \omega_s} = \psi_{\sigma, s, \omega_s'}.
    \end{equation}
 
(ii) For a permutation $\tau \in \mathfrak{S}_n$. Let $t:= \tau(s), \omega'_{t} := (\tau^{-1}(\mathbf{b}), r, Z)$. Then $\omega'_{t}$ is a seedlet at direction $t$ and we have
    \begin{eqnarray}
        P_{\omega'_t}(\tau^{-1}(\mathbf{x})) &=& P_{\omega_s}(\mathbf{x}), \label{Ptsigma = Ps ljy} \\
        \tau(\psi_{\sigma, s, \omega_s}(\mathbf{x})) &=& \psi_{\sigma\tau, s, \omega_s}(\mathbf{x}), \label{taupsi = psi}\\
        \psi_{\sigma, s, \omega_s}(\tau(\mathbf{x}))  &=& \psi_{\tau\sigma, t, \omega'_t}(\mathbf{x}). \label{psitau = psi}
    \end{eqnarray}
For the special case, when $\tau = \sigma^{-1}$, then $t = \sigma^{-1}(s)$, $\omega'_{t} := (\sigma(\mathbf{b}), r, Z)$ and we have
    \begin{equation} \label{psi inv}
        \psi_{\sigma, s, \omega_s}^{-1}(\mathbf{x}) = \psi_{\sigma^{-1}, t, \omega'_{t}}(\mathbf{x}).
    \end{equation}
\end{pty}

\begin{proof}
    (i) It is true, since 
    $$P_{\omega_s'}(\mathbf{x}) = \mathbf{x}^{r[\mathbf{b}]_{+}}Z'(\mathbf{x}^{-\mathbf{b}}) = \mathbf{x}^{r[\mathbf{b}]_{+}}\mathbf{x}^{-r\mathbf{b}}Z(\mathbf{x}^{\mathbf{b}}) = \mathbf{x}^{r[-\mathbf{b}]_{+}}Z(\mathbf{x}^{\mathbf{b}}) = P_{\omega_s}(\mathbf{x}).$$

    (ii) Clearly, $\pi_t(\tau^{-1}(\mathbf{b})) = b_s = 0$. So $\omega'_t$ is a seedlet at direction $t$. We have
    \begin{eqnarray*}
        P_{\omega'_t}(\tau^{-1}(\mathbf{x})) &=& \sum_{i=0}^{r} z_{i}(\tau^{-1}(\mathbf{x}))^{i[\tau^{-1}(\mathbf{b})]_{+} + (r - i)[-\tau^{-1}(\mathbf{b})]_{+}}\\
        &=& \sum_{i=0}^{r} z_{i}\mathbf{x}^{i[\mathbf{b}]_{+} + (r - i)[-\mathbf{b}]_{+}}\\
        &=& P_{\omega_s}(\mathbf{x}).
    \end{eqnarray*}
    By Equation \eqref{CSM'} and \eqref{sigmatau}, we have
    $$\tau(\psi_{\sigma, s, \omega_s}(\mathbf{x}))  = \bigg(\tau(\sigma(\mathbf{x}))\bigg)\bigg|_{x_s \gets \frac{P_{\omega_s}(\mathbf{x})}{x_s}} = \bigg(\sigma\tau(\mathbf{x})\bigg)\bigg|_{x_s \gets \frac{P_{\omega_s}(\mathbf{x})}{x_s}} = \psi_{\sigma\tau, s, \omega_s}(\mathbf{x}).$$

    Let $\mathbf{y} = \tau(\mathbf{x})$. Then 
    \begin{eqnarray*}
        \psi_{\sigma, s, \omega_s}(\tau(\mathbf{x})) &=& \psi_{\sigma, s, \omega_s}(\mathbf{y})\\
            &=& \bigg(\sigma(\mathbf{y})\bigg)\bigg|_{y_s \gets \frac{P_{\omega_s}(\mathbf{y})}{y_s}}   \qquad\qquad (\text{By } \eqref{CSM'})\\
            &=& \bigg(\sigma(\tau(\mathbf{x}))\bigg) \bigg|_{x_t \gets \frac{P_{\omega_s}(\tau(\mathbf{x}))}{x_t}} \\
            &=& \bigg(\tau\sigma(\mathbf{x}) \bigg)\bigg|_{x_t \gets \frac{P_{\omega'_t}(\mathbf{x})}{x_t}} \qquad\qquad (\text{By } \eqref{sigmatau}, \eqref{Ptsigma = Ps ljy})\\
            &=& \psi_{\tau\sigma, t, \omega'_t}(\mathbf{x}).
    \end{eqnarray*}
    
    When $\tau = \sigma^{-1}$. Denote that $\mathbf{x}' := \psi_{\sigma, s, \omega_s}(\mathbf{x})$ and $\mathbf{x}'' := \psi_{\sigma^{-1}, t, \omega'_{t}}(\mathbf{x}')$. For $i \ne s$, we know $x_i'' = x'_{\sigma^{-1}(i)} = x_i$. And we have
        $$x_s'' = \frac{P_{\omega'_t}(\mathbf{x}')}{x_t'} = \frac{P_{\omega_s}(\sigma^{-1}(\mathbf{x}'))}{x_t'} = \frac{P_{\omega_s}(\sigma^{-1}(\mathbf{x}'))}{P_{\omega_s}(\mathbf{x})/x_s} = x_s\frac{P_{\omega_s}(\mathbf{x})|_{x_s \gets x_s'}}{P_{\omega_s}(\mathbf{x})} = x_s.$$
    The last equality holds because $b_s = 0$ and the variable $x_s$ does not appear in $P_{\omega_s}(\mathbf{x})$. 
    
    So, we know $\psi_{\sigma^{-1}, t, \omega'_{t}}(\psi_{\sigma, s, \omega_s}(\mathbf{x})) = \mathbf{x}$.
\end{proof}

We aim to find a Laurent polynomial $F(\mathbf{x})$ which is invariant under a given cluster symmetric map $\psi_{\sigma,s,\omega_s}$. We define such Laurent polynomials.

\begin{DEF} \label{def 1-cspoly}
    Given a cluster symmetric map $\psi_{\sigma,s,\omega_s}$, if there exists a Laurent polynomial $F(\mathbf{x}) \in \mathbb{Q}[\mathbf{x}^{\pm}]$, such that $F(\mathbf{x})$ is invariant under the map $\psi_{\sigma,s,\omega_s}$, that is,  $F(\psi_{\sigma,s,\omega_s}(\mathbf{x})) = F(\mathbf{x})$, then we call $F(\mathbf{x})$ a \textbf{cluster symmetric polynomial\footnote{\label{footnote name} Although it would be more appropriate to call it ``cluster symmetric \textbf{Laurent} polynomial'', we think it would be better to drop the term ``Laurent''. Our considerations are as follows. First, an important property in cluster algebras is the positive Laurent phenomenon (Theorem \ref{thm PLP}), and we believe that the term ``cluster'' implies ``Laurent''. Second, the name ``cluster symmetric Laurent polynomial'' is too tedious.} about} $\psi_{\sigma,s,\omega_s}$, or briefly, a \textbf{cluster symmetric polynomial}. And for any constant $c \in \mathbb{Q}$, the Diophantine equation $F(\mathbf{x}) = c$ is a \textbf{cluster symmetric equation}.
\end{DEF}

For example, the Laurent polynomials $F_1$ in \eqref{eq Laurent Markov}, $F_{2,i}$ in Table \ref{tab rank 2 1-cspoly} and the Markov-cluster polynomial $F_{3,i}$ in Table \ref{tab all} are all cluster symmetric polynomials. The Markov equation \eqref{eq markov} is a cluster symmetric equation.

We first classify Laurent polynomials. To do it, we define two types of degree functions of a Laurent polynomial.

\begin{DEF} Fix $k \in [1,n]$. We define two functions $\deg^k, \deg_k$ as follows. For a Laurent polynomial $h(\mathbf{x}) \in \mathbb{Q}[\mathbf{x}^{\pm}]$, if $h(\mathbf{x}) = 0$, we define $\deg^k h(\mathbf{x}) := 0$ and $\deg_k h(\mathbf{x}) := 0$; if $h(\mathbf{x}) =  \sum_{a_{\mathbf{j}} \in \mathbb{Q}^*} a_{\mathbf{j}}\mathbf{x}^{\mathbf{j}} \ne 0$, we define 
\begin{eqnarray*}
    \deg^k h(\mathbf{x}) &:=& \max_{\mathbf{j} \text{ with } a_{\mathbf{j}} \ne 0}\{\text{the degree of } x_k \text{ in }\mathbf{x}^{\mathbf{j}}\},\\ 
    \deg_k h(\mathbf{x}) &:=& \min_{\mathbf{j} \text{ with } a_{\mathbf{j}} \ne 0}\{\text{the degree of } x_k \text{ in }\mathbf{x}^{\mathbf{j}}\}. 
\end{eqnarray*} 
\end{DEF}

For example, we consider the Laurent polynomial $h(\mathbf{x}) := \frac{x_1^2+x_2^3}{x_1} = x_1 + x_1^{-1}x_2^3$. Clearly, $\deg^1 h(\mathbf{x}) = 1$, $\deg^2 h(\mathbf{x}) = 3$, $\deg_1 h(\mathbf{x}) = -1$ and $\deg_2 h(\mathbf{x}) = 0$.

\begin{pty}
    For all $k \in [1,n], \sigma \in \mathfrak{S}_n$ and $h(\mathbf{x}) \in \mathbb{Q}[\mathbf{x}^{\pm}]$, we have
\begin{eqnarray}
    \deg^k h(\sigma(\mathbf{x})) &=& \deg^{\sigma^{-1}(k)} h(\mathbf{x}), \label{deg^k_sigma}\\
    \deg_k h(\sigma(\mathbf{x})) &=& \deg_{\sigma^{-1}(k)} h(\mathbf{x}).\label{deg_k_sigma}
\end{eqnarray}
\end{pty}

\begin{proof}
    If $h(\mathbf{x}) = 0$, it is obvious. 
    If $h(\mathbf{x}) \ne 0$, suppose $h(\mathbf{x}) =  \sum_{a_{\mathbf{j}} \in \mathbb{Q}^*} a_{\mathbf{j}}\mathbf{x}^{\mathbf{j}}$.
    Then 
    \begin{eqnarray*}
        \deg^kh(\sigma(\mathbf{x})) &=& \max_{\mathbf{j} \text{ with } a_{\mathbf{j}} \ne 0}\{\text{the degree of } x_k \text{ in }(\sigma(\mathbf{x}))^{\mathbf{j}}\}  \\
        &=& \max_{\mathbf{j} \text{ with } a_{\mathbf{j}} \ne 0}\{\text{the degree of } x_k \text{ in }\mathbf{x}^{\sigma^{-1}(\mathbf{j})}\}\\
        &=& \max_{\mathbf{j} \text{ with } a_{\mathbf{j}} \ne 0}\{\text{the degree of } x_{\sigma^{-1}(k)} \text{ in }\mathbf{x}^{\mathbf{j}}\}\\
        &=& \deg^{\sigma^{-1}(k)}h(\mathbf{x}).
    \end{eqnarray*}
    Similarly, Relation \eqref{deg_k_sigma} holds.
\end{proof}


\begin{DEF} \label{def lp}
Let $F(\mathbf{x}) \in \mathbb{Q}[\mathbf{x}^{\pm}]$. We call the Laurent polynomial $F(\mathbf{x})$ is of \textbf{type} $\frac{\bm{\eta}}{\mathbf{d}}$, if the unique expansion of $F(\mathbf{x})$ is
    $$F(\mathbf{x}) = \frac{T(\mathbf{x})}{\mathbf{x}^{\mathbf{d}}} = \frac{T(x_1, \cdots, x_n)}{x_1^{d_1} \cdots x_n^{d_n}},$$
where $T(\mathbf{x}) \in \mathbb{Q}[\mathbf{x}]$, $\mathbf{d}  \in  \mathbb{Z}^n$, $\bm{\eta} := (\deg^{1}T(\mathbf{x}), \cdots, \deg^{n}T(\mathbf{x})) \in \mathbb{Z}_{\geq 0}^n$ and
\begin{equation}\label{xnmidT}
    x_i \nmid T(\mathbf{x}), \quad \forall~i \in [1,n].
\end{equation}
\end{DEF}

\begin{rmk} Obviously, the relations \eqref{xnmidT} imply that
    \begin{equation}
        \deg_{k} T(\mathbf{x}) = 0, \quad \forall~ k \in [1,n]. \label{deg_kT=0}
    \end{equation}
\end{rmk}

\begin{eg} \label{eg eta d}
    (i) The Laurent polynomial 
    $F_1(\mathbf{x}) = \frac{x_1^2 + x_2^2 + x_3^2}{x_1x_2x_3}$ 
    is of type $\frac{(2,2,2)}{(1,1,1)}$. It is related to the Markov equation \eqref{eq markov}.

    (ii) The Laurent polynomial 
    $F_4(\mathbf{x}) = \frac{x_1^2x_4^2 + \alpha x_1x_3^3 + \alpha x_2^3x_4 + \beta x_2^2x_3^2}{x_1x_2x_3x_4}$
    is of type $\frac{(2,3,3,2)}{(1,1,1,1)}$. It was found by Hone and Swart in \cite{HoneSwart} and related to the Somos 4 sequence.

    (iii) The Laurent polynomial 
    $ F_{3,6}(\mathbf{x}) = \frac{x_1^2 + x_2^4 + x_3^4 + 2x_1x_2^2 + kx_2^2x_3^2 + 2x_1x_3^2}{x_1x_2^2x_3^2}$ 
    is of type $\frac{(2,4,4)}{(1,2,2)}$. It was found by Gyoda and Matsushita in \cite{Gyoda}.
\end{eg}

Once we find a cluster symmetric polynomial about a given cluster symmetric map, we can obtain another cluster symmetric polynomial about another corresponding cluster symmetric map.

\begin{prop}\label{prop new csp}
Let $\psi_{\sigma, s, \omega_s}$ be a cluster symmetric map, where $\omega_s := (\mathbf{b}, r, Z)$. Suppose $F(\mathbf{x})$ is a cluster symmetric polynomial about $\psi_{\sigma, s, \omega_s}$.

(i) Let $\omega_s' = (-\mathbf{b}, r, Z')$, where $Z'(u) = u^rZ(1/u)$. Then $F(\mathbf{x})$ is a cluster symmetric polynomial about $ \psi_{\sigma, s, \omega'_s}$.

(ii) For $\tau \in \mathfrak{S}_n$. Let $t:= \tau(s), \omega'_{t} := (\tau^{-1}(\mathbf{b}), r, Z)$ and $\widetilde{F}(\mathbf{x}) := F(\tau(\mathbf{x}))$. Then the Laurent polynomial $\widetilde{F}(\mathbf{x})$ is a cluster symmetric polynomial about  $\psi_{\tau\sigma\tau^{-1}, t, \omega'_t}$.

\end{prop}
\begin{proof}
    (i) It is obvious, since Equation \eqref{psi b = psi -b}.

    (ii) 
    Let $\mathbf{y} = \tau(\mathbf{x})$. Then
    \begin{eqnarray*}
        \widetilde{F}(\psi_{\tau\sigma\tau^{-1}, t, \omega'_t}(\mathbf{x})) &=& {F}(\tau(\psi_{\tau\sigma\tau^{-1}, t, \omega'_t}(\mathbf{x})))\\
            &=& {F}(\psi_{\tau\sigma, t, \omega'_t}(\mathbf{x}))  \quad\qquad\qquad (\text{By } \eqref{taupsi = psi}) \\
            &=& {F}(\psi_{\sigma, s, \omega_s}(\tau(\mathbf{x})))  \qquad\qquad (\text{By } \eqref{psitau = psi}) \\
            &=& {F}(\tau(\mathbf{x}))\\
            &=& \widetilde{F}(\mathbf{x}).
    \end{eqnarray*}
\end{proof}

\subsection{Construction of a cluster symmetric polynomial}\label{CA2diop}

In this subsection, we give a method to construct the cluster symmetric polynomial about the cluster symmetric map $\psi_{\sigma, s, \omega_s}$. We begin by establishing the notation for the expansion of the Laurent polynomial.

Let $F(\mathbf{x})$ be a Laurent polynomial of type $\frac{\bm{\eta}}{\mathbf{d}}$ in $\mathbb{Q}[\mathbf{x}^{\pm}]$. Suppose its expansion is 
    $$F(\mathbf{x}) = \frac{T(\mathbf{x})}{\mathbf{x}^{\mathbf{d}}} = {\mathbf{x}^{-\mathbf{d}}}\sum_{\mathbf{j} \in \mathcal{N}} a_{\mathbf{j}}\mathbf{x}^{\mathbf{j}},$$ 
where $\mathcal{N} := \{ \mathbf{j} \in \mathbb{Z}_{\geq 0}^n ~|~  0 \leq \pi_i(\mathbf{j}) \leq \pi_i(\bm{\eta}), \forall~ i \in [1,n]\}$  and $\pi_i(\mathbf{j})$ is meant to be the $i$-th component of the $n$-tuple $\mathbf{j}$. For $k \in [1,n]$ and $i \in \mathbb{Z}$, we define a subset of $\mathcal{N}$ as
    $\pi_k^{(i)}(\mathcal{N}) := \{\mathbf{j} \in \mathcal{N} ~|~ \pi_k(\mathbf{j}) = i\}$
and a polynomial in $\mathbb{Q}[\mathbf{x}]$ as
\begin{equation}
        f_{k,i}(\mathbf{x}) := \sum_{\mathbf{j} \in \pi_k^{(i)}(\mathcal{N})} a_{\mathbf{j}}\mathbf{x}^{\mathbf{j}-i\mathbf{e}_k},\label{f_sum}
\end{equation}
where $\mathbf{e}_k$' are standard basis. Then the polynomial $T(\mathbf{x})$ can be written as
\begin{equation}
    T(\mathbf{x}) = \sum_{i=0}^{\eta_k} f_{k, i}(\mathbf{x}) \mathbf{x}^{i\mathbf{e}_k}.\label{T=sumf}
\end{equation}
So the Laurent polynomial $F(\mathbf{x})$ can be written as 
\begin{equation}
    F(\mathbf{x}) = \mathbf{x}^{-\mathbf{d}} \sum_{i=0}^{\eta_k} f_{k, i}(\mathbf{x}) \mathbf{x}^{i\mathbf{e}_k}.\label{F=sumf}
\end{equation}

\begin{eg} \label{eg: F = sum f}
    We consider the Laurent polynomial $$F_4(\mathbf{x}) = \frac{x_1^2x_4^2 + \alpha x_1x_3^3 + \alpha x_2^3x_4 + \beta x_2^2x_3^2}{x_1x_2x_3x_4}$$
    in Example \ref{eg eta d} (ii). Then we have
    \begin{eqnarray*} 
        f_{1,0}(\mathbf{x}) = \alpha x_2^3x_4 + \beta x_2^2x_3^2, \quad & f_{1,1}(\mathbf{x}) = \alpha x_3^3, \quad & f_{1,2}(\mathbf{x}) = x_4^2,\\
        f_{4,0}(\mathbf{x}) = \alpha x_1x_3^3 + \beta x_2^2x_3^2, \quad & f_{4,1}(\mathbf{x}) = \alpha x_2^3, \quad & f_{4,2}(\mathbf{x}) = x_1^2.\\
    \end{eqnarray*}
\end{eg}

Regarding the polynomial $f_{k,i}(\mathbf{x})$, we have described some of their properties in the following two lemmas, which help to prove the main theorems of this subsection.

\begin{lem}\label{lem_f} 
(i) For any $k \in [1,n]$. We have
    \begin{eqnarray}
        f_{k, 0}(\mathbf{x}) &\ne&  0, \label{fk0}\\
        f_{k, \eta_k}(\mathbf{x}) &\ne&  0,\label{fknk}\\
        f_{k, j}(\mathbf{x}) &=& 0, \quad \text{for all } j \notin [0, \eta_k]. \label{f_kj_notin}
    \end{eqnarray}

(ii) For any $k \in [1,n], \sigma \in \mathfrak{S}_n, i \in \mathbb{Z}$. We have
    \begin{eqnarray}
        T(\sigma(\mathbf{x})) &=& \sum_{i=0}^{\eta_{\sigma^{-1}(k)}} f_{{\sigma^{-1}(k)}, i}(\sigma(\mathbf{x})) \mathbf{x}^{i\mathbf{e}_{k}},\label{Tsigma=sumf}\\
            f_{k,i}(\sigma(\mathbf{x})) &=& \sum_{\mathbf{t} \in \pi_k^{(i)}(\mathcal{N})} a_{\mathbf{t}}\mathbf{x}^{\sigma^{-1}(\mathbf{t})-i\mathbf{e}_{\sigma(k)}} = \sum_{\mathbf{j} \in \pi_{\sigma(k)}^{(i)}(\sigma^{-1}(\mathcal{N}))} a_{\sigma(\mathbf{j})}\mathbf{x}^{\mathbf{j} - i\mathbf{e}_{\sigma(k)}}.\label{f_sigma_sum}
    \end{eqnarray}
    
(iii) Given a cluster symmetric map $\psi_{\sigma, s, \omega_s}$ and an exchange polynomial $P_{\omega_s}$. Let $t := \sigma^{-1}(s)$. We have
    \begin{eqnarray}
        (\psi_{\sigma, s, \omega_s}(\mathbf{x}))^{\mathbf{d}} &=& \mathbf{x}^{\sigma^{-1}(\mathbf{d})}\left(\frac{P_{\omega_s}(\mathbf{x})}{x_s^2}\right)^{d_{t}} \label{sigma_mu_d} \\
        T(\psi_{\sigma, s, \omega_s}(\mathbf{x})) &=&   \sum_{i=0}^{\eta_{t}} f_{t,i}(\sigma(\mathbf{x})) \bigg(\frac{P_{\omega_s}(\mathbf{x})}{x_s}\bigg)^{i} \label{T_sigma_mu=sum}
    \end{eqnarray}
\end{lem}
\begin{proof} (i) Trivial.

    (ii) By Equation \eqref{T=sumf}, we have
    \begin{eqnarray*}
        T(\sigma(\mathbf{x})) &=& \sum_{i=0}^{\eta_{\sigma^{-1}(k)}} f_{{\sigma^{-1}(k)}, i}(\sigma(\mathbf{x})) (\sigma(\mathbf{x}))^{i\mathbf{e}_{\sigma^{-1}(k)}}\\
        &=& \sum_{i=0}^{\eta_{\sigma^{-1}(k)}} f_{{\sigma^{-1}(k)}, i}(\sigma(\mathbf{x}))\mathbf{x}^{\sigma^{-1}(i\mathbf{e}_{\sigma^{-1}(k)})}\\
        &=& \sum_{i=0}^{\eta_{\sigma^{-1}(k)}} f_{{\sigma^{-1}(k)}, i}(\sigma(\mathbf{x})) \mathbf{x}^{i\mathbf{e}_{k}}.
    \end{eqnarray*}
        
Since $\pi_k(\sigma(\mathbf{j})) = \pi_{\sigma(k)}(\mathbf{j})$, it is easy to check that 
    $$\mathbf{j} \in \pi_k^{(i)}(\mathcal{N})  ~\Leftrightarrow~ \sigma(\mathbf{j}) \in \pi_{\sigma^{-1}(k)}^{(i)}(\sigma(\mathcal{N})).$$
Then we have
    \begin{eqnarray*}
        f_{k,i}(\sigma(\mathbf{x})) &=& \sum_{\mathbf{t} \in \pi_k^{(i)}(\mathcal{N})} a_{\mathbf{t}}(\sigma(\mathbf{x}))^{\mathbf{t} - i\mathbf{e}_{k}}\\
        &=& \sum_{\mathbf{t} \in \pi_k^{(i)}(\mathcal{N})} a_{\mathbf{t}}\mathbf{x}^{\sigma^{-1}(\mathbf{t} - i\mathbf{e}_{k})}\\
        &=& \sum_{\mathbf{t} \in \pi_k^{(i)}(\mathcal{N})} a_{\mathbf{t}}\mathbf{x}^{\sigma^{-1}(\mathbf{t}) - i\mathbf{e}_{\sigma(k)}}\\
        &=& \sum_{\mathbf{j} \in \pi_{\sigma(k)}^{(i)}(\sigma^{-1}(\mathcal{N}))} a_{\sigma(\mathbf{j})}\mathbf{x}^{\mathbf{j} - i\mathbf{e}_{\sigma(k)}}.
    \end{eqnarray*}

(iii)  By Equation \eqref{CSM'}, we have
\begin{eqnarray*} 
    (\psi_{\sigma, s, \omega_s}(\mathbf{x}))^{\mathbf{d}} &=& (\sigma(\mathbf{x}))^{\mathbf{d}}|_{x_s \gets {P_{\omega_s}(\mathbf{x})}/x_s}\\
    &=& \mathbf{x}^{\sigma^{-1}(\mathbf{d})}|_{x_s \gets {P_{\omega_s}(\mathbf{x})}/x_s}\\
    &=& \mathbf{x}^{\sigma^{-1}(\mathbf{d})}\left(\frac{P_{\omega_s}(\mathbf{x})}{x_s^2}\right)^{d_{\sigma^{-1}(s)}}
\end{eqnarray*}
and
\begin{eqnarray*}
    T(\psi_{\sigma, s, \omega_s}(\mathbf{x})) &=& \bigg(T(\sigma(\mathbf{x}))\bigg)\bigg|_{x_s \gets \frac{P_{\omega_s}(\mathbf{x})}{x_s}} \\
        &=&  \sum_{i=0}^{\eta_{t}} \bigg( f_{t,i}(\sigma(\mathbf{x})) \mathbf{x}^{i\mathbf{e}_s} \bigg) \bigg|_{x_s \gets \frac{P_{\omega_s}(\mathbf{x})}{x_s}}  \qquad\qquad\quad (\text{By } \eqref{Tsigma=sumf})\\
        &=&  \sum_{i=0}^{\eta_{t}} \bigg( f_{t,i}(\sigma(\mathbf{x})) \bigg) \bigg|_{x_s \gets \frac{P_{\omega_s}(\mathbf{x})}{x_s}} \bigg(\frac{P_{\omega_s}(\mathbf{x})}{x_s}\bigg)^{i} \\
        &=&  \sum_{i=0}^{\eta_{t}} f_{t,i}(\sigma(\mathbf{x})) \bigg(\frac{P_{\omega_s}(\mathbf{x})}{x_s}\bigg)^{i}. \qquad\quad\quad\quad\quad\quad\quad  (\text{By } \eqref{lem_degf})\\
\end{eqnarray*}
\end{proof}

\begin{lem} For all $k, j \in [1,n]$, the following relations hold,

(i)
    \begin{equation}
            0 \leq \deg_k f_{j,i}(\mathbf{x}) \leq \deg^k f_{j,i}(\mathbf{x}) \leq \eta_k, ~\forall~i \in [0, \eta_j]. \label{deg_f}
    \end{equation}

(ii) 
    \begin{eqnarray}
        \deg^k f_{k,i}(\mathbf{x}) &=& 0, ~\text{for all }~i \in [0, \eta_k]. \label{lem_degf} \\
        \text{If } k \ne j, \text{then } \deg^k f_{j, i_k}(\mathbf{x}) &=& \eta_k, ~\text{for some } i_k \in [0, \eta_k]. \label{deg^k=nk}
    \end{eqnarray}

(iii)
    \begin{equation}
        \deg_k f_{j, i_k}(\mathbf{x}) = 0, ~\quad \text{for some } i_k \in [0, \eta_k]. \label{deg_kf=0}
    \end{equation}

(iv) For a seedlet $\omega_s := (\mathbf{b}, r, Z)$. We have
    \begin{eqnarray}
        \deg^j P_{\omega_s}(\mathbf{x}) &=& r|b_{j}|,\label{deg_p_br}\\
        \deg_j P_{\omega_s}(\mathbf{x}) &=& 0,\label{deg_p_0}
    \end{eqnarray}

\end{lem}

\begin{proof}
    (i) By Equation \eqref{deg_kT=0}, we have 
    $$0 = \deg_k T(\mathbf{x}) \leq \deg_k f_{j,i}(\mathbf{x}) \leq \deg^k f_{j,i}(\mathbf{x}) \leq \deg^k T(\mathbf{x}) = \eta_k.$$ 

    (ii) When $T(\mathbf{x}) = 0$, it is true. Suppose $T(\mathbf{x}) \ne 0$. By the definition of $f_{k,i}(\mathbf{x})$ in \eqref{f_sum}, it is obvious that Equation \eqref{lem_degf} holds.
    If $k \ne j$, by Equation \eqref{T=sumf}, we have
    
    \centerline{$\eta_k = \deg^k T(\mathbf{x}) =  \deg^k \sum_{i=0}^{\eta_j} f_{j, i}(\mathbf{x}) = \mathop{\max_{0\leq i \leq n_j}}_{f_{j,i} \ne 0} \{ \deg_k f_{j, i}(\mathbf{x})\}.$}
    
    Hence there exists $i_k \in [0, \eta_k]$, such that, $\deg^k f_{j, i_k}(\mathbf{x}) = 0$.
    
    (iii) When $T(\mathbf{x}) = 0$, it is true. Suppose $T(\mathbf{x}) \ne 0$. If $k = j$. By (i) and (ii), we have $\deg_k f_{k, i}(\mathbf{x}) = 0$ for all $i \in [0, \eta_k]$. If $k \ne j$, by Equation \eqref{deg_kT=0} and \eqref{T=sumf}, we have
    $$0 = \deg_k T(\mathbf{x}) = \deg_k \sum_{i=0}^{\eta_j} f_{j, i}(\mathbf{x}) \mathbf{x}^{i\mathbf{e}_j} = \deg_k \sum_{i=0}^{\eta_j} f_{j, i}(\mathbf{x}) = \mathop{\min_{0\leq i \leq \eta_j}}_{f_{j,i} \ne 0} \{ \deg_k f_{j, i}(\mathbf{x})\}.$$
    Hence there exists $i_k \in [0, \eta_k]$, such that, $\deg_k f_{j, i_k}(\mathbf{x}) = 0$.

(iv) Since the equations \eqref{ne 0}, we have
    \begin{eqnarray*}
        \deg^j P_{\omega_s}(\mathbf{x}) &=& \deg^j \left( \sum_{i=0}^{r} z_{i}\mathbf{x}^{i[\mathbf{b}]_{+} + (r_k - i)[-\mathbf{b}]_{+}} \right)\\
            &=& \max_{0 \leq i \leq r, z_{i} \ne 0}\{ \pi_j(i[\mathbf{b}]_{+} + (r - i)[-\mathbf{b}]_{+}) \}\\
            &=& \max_{0 \leq i \leq r, z_{i} \ne 0}\{ i[b_j]_+ + (r - i)[-b_j]_+ \}\\
            &=& \max\{ r[b_j]_+, r[-b_j]_+ \}\\
            &=& r|b_{j}|
    \end{eqnarray*}
and
    \begin{eqnarray*}
        \deg_j P_{\omega_s}(\mathbf{x}) &=& \deg_j \left( \sum_{i=0}^{r} z_{i}\mathbf{x}^{i[\mathbf{b}]_{+} + (r_k - i)[-\mathbf{b}]_{+}} \right)\\
            &=& \min_{0 \leq i \leq r, z_{i} \ne 0}\{ \pi_j(i[\mathbf{b}]_{+} + (r - i)[-\mathbf{b}]_{+}) \}\\
            &=& \min_{0 \leq i \leq r, z_{i} \ne 0}\{ i[b_j]_+ + (r - i)[-b_j]_+ \}\\
            &=& \min\{ r[b_j]_+, r[-b_j]_+ \}\\
            &=& 0.
    \end{eqnarray*}
\end{proof}

By giving equivalence conditions under which the relation $F(\psi_{\sigma, s, \omega_s}(\mathbf{x})) = F(\mathbf{x})$ holds, we describe the cluster symmetric polynomials.

\begin{thm} \label{thm len1_T2f} 
Given a cluster symmetric map $\psi_{\sigma, s, \omega_s}$.
Let $F(\mathbf{x})$ be a Laurent polynomial of type $\frac{\bm{\eta}}{\mathbf{d}}$ in $\mathbb{Q}[\mathbf{x}^{\pm}]$, and suppose that its expansion is
    $$F(\mathbf{x}) = \mathbf{x}^{-\mathbf{d}} \sum_{i=0}^{\eta_k} f_{k, i}(\mathbf{x}) \mathbf{x}^{i\mathbf{e}_k}$$
as shown in Equation \eqref{F=sumf}, where $k \in [1,n]$. Then the relation 
    \begin{equation}
        F(\psi_{\sigma, s, \omega_s}(\mathbf{x})) = F(\mathbf{x})\label{len1}
    \end{equation}
holds, if and only if, the following relations      
    \begin{eqnarray}
        f_{\sigma^{-1}(s),i}(\sigma(\mathbf{x})) &=& f_{s,\eta_s-i}(\mathbf{x})P_{\omega_s}^{d_s - i}(\mathbf{x}),~ \forall~i \in [0, \eta_s], \label{len1_f=fp@s}\\
        \mathbf{d} &=& \sigma{(\mathbf{d})}, \label{len1_d=sigmad}\\
        \eta_s = \eta_{\sigma^{-1}(s)} &=& 2d_s = 2d_{\sigma^{-1}(s)} \label{len1_eta_s=2d_s}
    \end{eqnarray}
hold.
\end{thm}

\begin{proof} For convenience, we denote $t := \sigma^{-1}(s)$, $\psi := \psi_{\sigma, s, \omega_s}$, $P := P_{\omega_s}$ and
    \begin{eqnarray*}
        \delta(\mathbf{d}, \sigma, s) &:=& \mathbf{d} - \sigma^{-1}(\mathbf{d}) + (d_{\sigma^{-1}(s)} -d_s)\mathbf{e}_s.
    \end{eqnarray*}

\textbf{STEP 1.} We claim that the equation $F(\psi(\mathbf{x})) = F(\mathbf{x})$ holds, if and only if, the following equations
    \begin{eqnarray}
        \mathbf{x}^{\delta(\mathbf{d}, \sigma, s)} f_{t,i}(\sigma(\mathbf{x}))  = f_{s,d_s+d_t-i}(\mathbf{x})P^{d_t - i}(\mathbf{x}),&~ \forall~i \in [0, \eta_s], \label{len1_xf=fp@s}\\
        \eta_s = \eta_t = d_s + d_t,& \label{len1_sigma_n}
    \end{eqnarray}
hold.

By the equations \eqref{T=sumf}, \eqref{sigma_mu_d} and \eqref{T_sigma_mu=sum}, we have    
    \begin{eqnarray*}
        && (\psi(\mathbf{x}))^{\mathbf{d}}\mathbf{x}^{\mathbf{d}}\bigg[F(\psi(\mathbf{x})) - F(\mathbf{x})\bigg]\\
        &=& \mathbf{x}^{\mathbf{d}}T(\psi(\mathbf{x})) - (\psi(\mathbf{x}))^{\mathbf{d}}T(\mathbf{x})\\
        &=& \mathbf{x}^{\mathbf{d}}\sum_{i=0}^{\eta_t}f_{t,i}(\sigma(\mathbf{x})) \bigg(\frac{P(\mathbf{x})}{x_s}\bigg)^{i} - \left(\frac{P(\mathbf{x})}{x_s^2}\right)^{d_t}\mathbf{x}^{\sigma^{-1}(\mathbf{d})}\sum_{i=0}^{\eta_s}f_{s,i}(\mathbf{x})x_{s}^{i}\\
        &=& \mathbf{x}^{\mathbf{d} - d_s\mathbf{e}_s} \sum_{i=0}^{\eta_t}f_{t,i}(\sigma(\mathbf{x})) P^i(\mathbf{x}) x_s^{d_s - i} - P^{d_t}(\mathbf{x}) \mathbf{x}^{\sigma^{-1}(\mathbf{d}) - d_t\mathbf{e}_s} \sum_{i=0}^{\eta_s}f_{s,i}(\mathbf{x})x_{s}^{i - d_t}\\
        &=& \mathbf{x}^{\mathbf{d} - d_s\mathbf{e}_s} \sum_{i=0}^{\eta_t}f_{t,i}(\sigma(\mathbf{x})) P^i(\mathbf{x}) x_s^{d_s - i}  - P^{d_t}(\mathbf{x}) \mathbf{x}^{\sigma^{-1}(\mathbf{d}) - d_t\mathbf{e}_s} \sum_{i=d_s+d_t-\eta_s}^{d_s+d_t} f_{s,d_s+d_t-i}(\mathbf{x})x_{s}^{d_s - i}\\
        &=& \sum_{i=0}^{\eta_t} \left(\mathbf{x}^{\mathbf{d} - d_s\mathbf{e}_s}f_{t,i}(\sigma(\mathbf{x})) P^i(\mathbf{x})\right) x_s^{d_s - i} - \sum_{i=d_s+d_t-\eta_s}^{d_s+d_t} \left(\mathbf{x}^{\sigma^{-1}(\mathbf{d}) - d_t\mathbf{e}_s}f_{s,d_s+d_t-i}(\mathbf{x})P^{d_t}(\mathbf{x}) \right)x_{s}^{d_s - i}.
    \end{eqnarray*}

$(\Leftarrow):$ By the above equation, it is clear that the equations $d_s + d_t = \eta_s = \eta_t$ and 
    $$\mathbf{x}^{\delta(\mathbf{d}, \sigma, s)} f_{t,i}(\sigma(\mathbf{x})) = f_{s,d_s+d_t-i}(\mathbf{x})P^{d_t - i}(\mathbf{x}),~ \forall i \in [0, \eta_s]$$
imply the equation $F(\psi(\mathbf{x})) = F(\mathbf{x})$.

$(\Rightarrow):$ Suppose Equation \eqref{len1} holds, then by the above equation we have
    \footnotesize{\begin{equation*}
        \sum_{i=0}^{\eta_t} \left(\mathbf{x}^{\mathbf{d} - d_s\mathbf{e}_s}f_{t,i}(\sigma(\mathbf{x})) P^i(\mathbf{x})\right) x_s^{d_s - i} = \sum_{i=d_s+d_t-\eta_s}^{d_s+d_t} \left(\mathbf{x}^{\sigma^{-1}(\mathbf{d}) - d_t\mathbf{e}_s}f_{s,d_s+d_t-i}(\mathbf{x})P^{d_t}(\mathbf{x}) \right)x_{s}^{d_s - i}.
    \end{equation*}}

    Since  $x_s, P(\mathbf{x}) \ne 0$ and the relations 
    $$f_{t,0}, f_{t, \eta_t}, f_{s, 0}, f_{s, \eta_s} \ne 0,$$
    which from the equations \eqref{fk0} and \eqref{fknk}, the above equation implies that 
    $$ 0 \leq d_s+d_t-\eta_s,\quad \eta_t \leq d_s+d_t,\quad d_s+d_t \leq \eta_t,\quad d_s+d_t-\eta_s \leq 0,$$
    or
    $$\eta_s = \eta_t = d_s + d_t.$$ 
    That is, Equation \eqref{len1_sigma_n} holds. Then we have 
    $$\mathbf{x}^{\mathbf{d} - d_s\mathbf{e}_s}f_{t,i}(\sigma(\mathbf{x})) P^i(\mathbf{x}) = \mathbf{x}^{\sigma^{-1}(\mathbf{d}) - d_t\mathbf{e}_s}f_{s,d_s+d_t-i}(\mathbf{x})P^{d_t}(\mathbf{x}), ~\forall~i \in [0, \eta_s].$$
    So, Equation \eqref{len1_xf=fp@s} holds.

\textbf{STEP 2.} We prove that the equations
    \begin{align*}
        &\mathbf{x}^{\delta(\mathbf{d}, \sigma, s)} f_{t,i}(\sigma(\mathbf{x})) = f_{s,d_s+d_t-i}(\mathbf{x})P^{d_t - i}(\mathbf{x}),~ \forall~i \in [0, \eta_s], & \eqref{len1_xf=fp@s}\\
        &~~~~\qquad\qquad\qquad \eta_s = \eta_t = d_s + d_t, & \eqref{len1_sigma_n}
    \end{align*}
hold, if and only if, the equations
\begin{align*}
        &\qquad\qquad f_{t,i}(\sigma(\mathbf{x})) = f_{s,\eta_s-i}(\mathbf{x})P^{d_s - i}(\mathbf{x}),~ \forall~i \in [0, \eta_s] & \qquad \eqref{len1_f=fp@s}\\
        &\qquad\qquad\quad\qquad\qquad~ \mathbf{d} = \sigma{(\mathbf{d})}, & \qquad\eqref{len1_d=sigmad}\\
        &\qquad\qquad \quad \eta_s = \eta_t = 2d_s = 2d_t. & \qquad \eqref{len1_eta_s=2d_s}
\end{align*}
hold.

($\Leftarrow$): Since $\mathbf{d} = \sigma{(\mathbf{d})}$ implies $\delta(\mathbf{d}, \sigma, s) = \mathbf{d} - \sigma^{-1}(\mathbf{d}) + (d_t -d_s)\mathbf{e}_s = \mathbf{0}$, Equation \eqref{len1_xf=fp@s} holds.

($\Rightarrow$): For any $k \in [1,n]$, we have
    \begin{eqnarray*}
        \pi_k(\delta(\mathbf{d}, \sigma, s)) &=& \deg_k \mathbf{x}^{\delta(\mathbf{d}, \sigma, s)} f_{t,i}(\sigma(\mathbf{x})) -  \deg_k f_{t,i}(\sigma(\mathbf{x})) \\
        &=& \deg_k \mathbf{x}^{\delta(\mathbf{d}, \sigma, s)} f_{t,i}(\sigma(\mathbf{x})) -  \deg_{\sigma^{-1}(k)} f_{t,i}(\mathbf{x}) ~~~~~\qquad (\text{By }\eqref{deg_k_sigma})\\ 
        &=& \deg_k f_{s,\eta_t-i}(\mathbf{x})P^{d_t - i}(\mathbf{x}) - \deg_{\sigma^{-1}(k)} f_{t,i}(\mathbf{x}) \quad\qquad (\text{By }\eqref{len1_xf=fp@s})\\
        &=& \deg_k f_{s,\eta_s-i}(\mathbf{x}) - \deg_{\sigma^{-1}(k)} f_{t,i}(\mathbf{x}). \quad\qquad\quad\quad\qquad (\text{By } \eqref{deg_p_0})
    \end{eqnarray*}
Since Relation \eqref{deg_kf=0}, there exist $i' \in [0, \eta_s]$ and $i'' \in [0, \eta_t]$, such that
        $$\deg_k f_{s, i'}(\mathbf{x}) = 0, \quad \deg_{\sigma^{-1}(k)} f_{t, i''}(\mathbf{x}) = 0.$$
    Then we have
        \begin{eqnarray*}
            \pi_k(\delta(\mathbf{d}, \sigma, s)) &=& -  \deg_{\sigma^{-1}(k)} f_{t,\eta_s - i'}(\mathbf{x}),\\
            \pi_k(\delta(\mathbf{d}, \sigma, s)) &=& \deg_k f_{s,\eta_s-i''}(\mathbf{x}).
        \end{eqnarray*}
    So by the equations \eqref{deg_f}, we know 
        $$\pi_k(\delta(\mathbf{d}, \sigma, s)) \in [-\eta_{\sigma^{-1}(k)}, 0] \cap [0, \eta_k].$$
    
    Hence for all $k \in [1,n]$, we have $\pi_k(\delta(\mathbf{d}, \sigma, s)) = 0$, that is, $\delta(\mathbf{d}, \sigma, s) = \mathbf{0}$. Then Relation \eqref{len1_f=fp@s} holds and 
        \begin{eqnarray*}
            0 = \pi_k(\delta(\mathbf{d}, \sigma, s)) = \pi_k (\mathbf{d} - \sigma^{-1}(\mathbf{d}) + (d_t -d_s)\mathbf{e}_s) = \begin{cases}
                d_k - d_{\sigma^{-1}(k)}, &\text{if } k \ne s,\\
                0, &\text{if } k = s.
            \end{cases}
        \end{eqnarray*}
    Then we have
        $$d_t = d_{\sigma^{-1}(t)} = \cdots = d_{\sigma^{-(\operatorname{ord}(\sigma)-2)}(t)} = d_{\sigma^{-(\operatorname{ord}(\sigma)-1)}(t)} = d_s.$$
    So equations \eqref{len1_d=sigmad} and \eqref{len1_eta_s=2d_s} hold.    
\end{proof}

\begin{eg} \label{eg f = fp}
    It is easy to check that the Laurent polynomial 
    \begin{eqnarray*}
        F_4(\mathbf{x}) = \frac{x_1^2x_4^2 + \alpha x_1x_3^3 + \alpha x_2^3x_4 + \beta x_2^2x_3^2}{x_1x_2x_3x_4}
    \end{eqnarray*}
    shown in Example \ref{eg eta d} (ii) is invariant under the cluster symmetric map $\psi_{\sigma_{(1234)}, 1, \omega_1}$ defined in Example \ref{eg CSM somos and marsh}. By Example \ref{eg: F = sum f}, we have
    \begin{eqnarray*} 
        &f_{4,0}(\sigma_{(1234)}(\mathbf{x})) = \alpha x_2x_4^3 + \beta x_3^2x_4^2 = x_4^2(\alpha x_2x_4 + \beta x_3^2) = f_{1,2}(\mathbf{x})P_{\omega_s}(\mathbf{x}),\\
        &f_{4,1}(\sigma_{(1234)}(\mathbf{x})) = \alpha x_3^3 = f_{1,1}(\mathbf{x}),\\
        &f_{4,2}(\sigma_{(1234)}(\mathbf{x}))P_{\omega_s}(\mathbf{x}) = x_2^2(\alpha x_2x_4 + \beta x_3^2) = f_{1,0}(\mathbf{x}).
    \end{eqnarray*}
    So Relation \eqref{len1_f=fp@s} holds. 
\end{eg}

The above theorem urges us to describe the following relations,
$$f_{t,i}(\sigma(\mathbf{x})) = f_{s,\eta_s-i}(\mathbf{x})P_{\omega_s}^{d_s - i}(\mathbf{x}), \qquad \forall~i \in [0, \eta_s].$$
To do it, we introduce a lemma.

\begin{lem} \label{lem def c,b}
Given a seedlet $\omega_s := (\mathbf{b}, r, Z)$ and an exchange polynomial $P_{\omega_s}(\mathbf{x})$. For $k,l,i \in \mathbb{Z}_{\geq 0}$, we denote a coefficient
    $$c_{k, l} := \begin{cases}
        \sum_{l_1, \cdots, l_k \in [0, r], l_1 + \cdots + l_k = l} z_{l_1}\cdots z_{l_k}, &\text{ if } k > 0,\\
        1, &\text{ if } k = 0,
    \end{cases}$$
and a $n$-tuple $\mathbf{b}_{s,k,l}^{(i)} := l[\mathbf{b}]_{+} + (kr - l)[-\mathbf{b}]_{+} - i\mathbf{e}_s.$ Then,
    \begin{eqnarray}
        (P_{\omega_s}(\mathbf{x}))^k &=& \sum_{i=0}^{kr} c_{k,l}\mathbf{x}^{\mathbf{b}_{s,k,l}^{(0)}}, \label{P_k^j}\\
        \pi_s(\mathbf{b}_{s,k,l}^{(i)}) &=& -i.  \label{pikB_ki}
    \end{eqnarray}
\end{lem}

\begin{proof}
Since 
    $$(Z(u))^k = (\sum_{l=0}^{r}z_{l}u^l)^k = \sum_{l=0}^{kr}\bigg(\mathop{\sum_{l_1, \cdots, l_k \in [0, r]}}_{l_1 + \cdots + l_k = l} z_{l_1}\cdots z_{l_k} \bigg) u^l = \sum_{l=0}^{kr} c_{k,l} u^l,$$
we have
\begin{eqnarray*}
    (P_{\omega_s}(\mathbf{x}))^k &=&  \bigg(\mathbf{x}^{r[-\mathbf{b}]_{+}}Z(\mathbf{x}^{\mathbf{b}})\bigg)^k \qquad\qquad (\text{By } \eqref{Pk=sum ljy})\\
        &=& \mathbf{x}^{kr[-\mathbf{b}]_{+}}\sum_{l=0}^{kr} c_{k, l} \mathbf{x}^{l\mathbf{b}} \\
        &=& \sum_{l=0}^{kr} c_{k, l}\mathbf{x}^{l[\mathbf{b}]_{+} + (kr - l)[-\mathbf{b}]_{+}} \\
        &=& \sum_{l=0}^{kr} c_{k, l} \mathbf{x}^{\mathbf{b}_{s,k,l}^{(0)}}.
\end{eqnarray*}
By Definition \ref{seedlet}, we know $b_s = 0$, so $\pi_s(\mathbf{b}_{s,k,l}^{(i)}) = l[b_{s}]_{+} + (kr - l)[-b_{s}]_{+} - i = -i$.
\end{proof}

Theorem \ref{thm len1_T2f} formally describes the cluster symmetric polynomials, while the following theorem is used to construct them concretely. 

\begin{thm} \label{thm len1}
    Given a seedlet $\omega_s := (\mathbf{b}, r, Z)$ and a cluster symmetric map $\psi_{\sigma, s, \omega_s}$. 
    For any $\bm{\eta} \in \mathbb{Z}_{\geq 0}^n, \mathbf{d} \in \mathbb{Z}^n$ with $\mathbf{d} = \sigma{(\mathbf{d})}$ and $\eta_s = \eta_{t} = 2d_s = 2d_{t}$. Let $F(\mathbf{x})$ be a $\frac{\bm{\eta}}{\mathbf{d}}$ type Laurent polynomial in $\mathbb{Q}[\mathbf{x}^{\pm}]$ and its expansion is $F(\mathbf{x}) = {\mathbf{x}^{-\mathbf{d}}}\sum_{\mathbf{j} \in \mathcal{N}} a_{\mathbf{j}}\mathbf{x}^{\mathbf{j}}$, where $\mathcal{N} := \{ \mathbf{j} \in \mathbb{Z}_{\geq 0}^n ~|~  0 \leq \pi_i(\mathbf{j}) \leq \pi_i(\bm{\eta}), \forall~ i \in [1,n]\}$ and $a_{\mathbf{j}} \in \mathbb{Q}$ for all $\mathbf{j} \in \mathcal{N}$.
    Then the relation 
      $$F(\psi_{\sigma, s, \omega_s}(\mathbf{x})) = F(\mathbf{x})$$
    holds, if and only if, for any $k \in [0, d_s]$, the Laurent polynomial's coefficients $\{a_{\mathbf{j}} \in  \mathbb{Q} \mid \mathbf{j} \in \mathcal{N} \}$ satisfy the system of homogeneous linear equations $HLE(\sigma, s, \omega_s, \bm{\eta}, \mathbf{d},  k)$:
    \small{\begin{eqnarray*}
    \left\{ \begin{aligned}
        & 0 = a_{\sigma(\mathbf{j})} - \mathop{\sum_{0 \leq l \leq kr}}_{ \mathbf{j} - \mathbf{b}_{s,k,l}^{(2k)} \in \mathcal{N}} a_{\mathbf{j} - \mathbf{b}_{s,k,l}^{(2k)}}c_{k, l}, & \text{ if } \mathbf{j} \in \pi_s^{(d_s-k)}\bigg(\sigma^{-1}(\mathcal{N})\cap \bigcup_{0 \leq l \leq kr} (\mathcal{N} + \mathbf{b}_{s,k,l}^{(2k)})\bigg),\\
        &0 = \mathop{\sum_{0 \leq l \leq kr}}_{ \mathbf{j} - \mathbf{b}_{s,k,l}^{(2k)} \in \mathcal{N}} a_{\mathbf{j} - \mathbf{b}_{s,k,l}^{(2k)}}c_{k, l}, &\text{ if } \mathbf{j} \in \pi_s^{(d_s-k)}\bigg(\bigcup_{0 \leq l \leq kr} (\mathcal{N} + \mathbf{b}_{s,k,l}^{(2k)}) \setminus \sigma^{-1}(\mathcal{N})\bigg),\\
        &0 = a_{\sigma(\mathbf{j})}, &\text{ if } \mathbf{j} \in \pi_s^{(d_s-k)}\bigg(\sigma^{-1}(\mathcal{N}) \setminus \bigcup_{0 \leq l \leq kr} (\mathcal{N} + \mathbf{b}_{s,k,l}^{(2k)})\bigg),\\
    \end{aligned}
    \right.
    \end{eqnarray*}}
    and the system of homogeneous linear equations $HLE(\sigma^{-1}, t, \omega_t, \bm{\eta}, \mathbf{d}, k)$:
    \small{\begin{eqnarray*}
    \left\{ \begin{aligned}
        & 0 = a_{\sigma^{-1}(\mathbf{j})} - \mathop{\sum_{0 \leq l \leq kr}}_{ \mathbf{j} - \mathbf{v}_{t,k,l}^{(2k)} \in \mathcal{N}} a_{\mathbf{j} - \mathbf{v}_{t,k,l}^{(2k)}}c_{k, l}, & \text{ if } \mathbf{j} \in \pi_t^{(d_t-k)}\bigg(\sigma(\mathcal{N})\cap \bigcup_{0 \leq l \leq kr} (\mathcal{N} + \mathbf{v}_{t,k,l}^{(2k)})\bigg),\\
        &0 = \mathop{\sum_{0 \leq l \leq kr}}_{ \mathbf{j} - \mathbf{v}_{t,k,l}^{(2k)} \in \mathcal{N}} a_{\mathbf{j} - \mathbf{v}_{t,k,l}^{(2k)}}c_{k, l}, &\text{ if } \mathbf{j} \in \pi_t^{(d_t-k)}\bigg(\bigcup_{0 \leq l \leq kr} (\mathcal{N} + \mathbf{v}_{t,k,l}^{(2k)}) \setminus \sigma(\mathcal{N})\bigg),\\
        &0 = a_{\sigma^{-1}(\mathbf{j})}, &\text{ if } \mathbf{j} \in \pi_t^{(d_t-k)}\bigg(\sigma(\mathcal{N}) \setminus \bigcup_{0 \leq l \leq kr} (\mathcal{N} + \mathbf{v}_{t,k,l}^{(2k)})\bigg),\\
    \end{aligned}
    \right.
    \end{eqnarray*}}
    where $t := \sigma^{-1}(s)$, $\mathbf{v} := \sigma(\mathbf{b})$, $\omega_t := (\mathbf{v}, r, Z)$, $\pi_s^{(k)}(\mathcal{N}) := \{ \mathbf{j} \in \mathcal{N} \mid \pi_s(\mathbf{j}) = k \}$, $\mathbf{b}_{s,k,l}^{(i)} := l[\mathbf{b}]_{+} + (kr- l)[-\mathbf{b}]_{+} - i\mathbf{e}_s$ and 
    $$c_{k, l} := \begin{cases}
        \sum_{l_1, \cdots, l_k \in [0, r], l_1 + \cdots + l_k = l} z_{l_1}\cdots z_{l_k}, &\text{ if } k > 0,\\
        1, &\text{ if } k = 0.
    \end{cases}$$
\end{thm}
\begin{proof}
    (i) As in Equation \eqref{f_sum}, we denote polynomials
    $$f_{s,i}(\mathbf{x}) := \sum_{\mathbf{j} \in \pi_s^{(i)}(\mathcal{N})} a_{\mathbf{j}}\mathbf{x}^{\mathbf{j}-i\mathbf{e}_s}, \qquad f_{\sigma^{-1}(s),i}(\mathbf{x}) := \sum_{\mathbf{j} \in \pi_{\sigma^{-1}(s)}^{(i)}(\mathcal{N})} a_{\mathbf{j}}\mathbf{x}^{\mathbf{j}-i\mathbf{e}_{\sigma^{-1}(s)}},$$
where $\mathcal{N} := \{ \mathbf{j} \in \mathbb{Z}_{\geq 0}^n ~|~  0 \leq \pi_i(\mathbf{j}) \leq \pi_i(\bm{\eta}), \forall~ i \in [1,n]\}$. We claim that for  any $k \in [0, d_s]$, the equation
        \begin{equation} \label{len1_ftoa1}
    E(\sigma, s, \omega_s, \bm{\eta}, \mathbf{d},  k):~    f_{\sigma^{-1}(s), d_s-k}(\sigma(\mathbf{x})) = f_{s,d_s+k}(\mathbf{x})P_{\omega_s}^{k}(\mathbf{x})  
        \end{equation}
     holds, if and only if, the coefficients $a_{\mathbf{j}}$'s  satisfy the system of homogeneous linear equations $HLE(\sigma, s, \omega_s, \bm{\eta}, \mathbf{d},  k)$.

  By Equations \eqref{f_sigma_sum} and \eqref{P_k^j}, we have
    \small{\begin{eqnarray*}
        && f_{\sigma^{-1}(s), d_s-k}(\sigma(\mathbf{x})) -f_{s, d_s+k}(\mathbf{x}) P_{\omega_s}^{k}(\mathbf{x})\\
        &=&  \sum_{\mathbf{j} \in \pi_{s}^{(d_s-k)}(\sigma^{-1}(\mathcal{N}))} a_{\sigma(\mathbf{j})}\mathbf{x}^{\mathbf{j} - (d_s-k)\mathbf{e}_s}
        - \bigg(\sum_{\mathbf{t} \in \pi_s^{(d_s+k)}(\mathcal{N})} a_{\mathbf{t}}\mathbf{x}^{\mathbf{t} - (d_s+k)\mathbf{e}_{s}}\bigg)\bigg( \sum_{l=0}^{kr} c_{s,k, l}\mathbf{x}^{\mathbf{b}_{s,k,l}^{(0)}} \bigg)\\
        &=& \sum_{\mathbf{j} \in \pi_{s}^{(d_s-k)}(\sigma^{-1}(\mathcal{N}))} a_{\sigma(\mathbf{j})}\mathbf{x}^{\mathbf{j} - (d_s-k)\mathbf{e}_s } - \sum_{l=0}^{kr} \sum_{\mathbf{t} \in \pi_s^{(d_s+k)}(\mathcal{N})} a_{\mathbf{t}}c_{k,l}\mathbf{x}^{\mathbf{t} + \mathbf{b}_{s,k,l}^{(d_s+k)}}\\
        &=& \mathbf{x}^{ -(d_s-k)\mathbf{e}_s} \bigg(\sum_{\mathbf{j} \in \pi_{s}^{(d_s-k)}(\sigma^{-1}(\mathcal{N}))} a_{\sigma(\mathbf{j})}\mathbf{x}^{\mathbf{j} } - \sum_{l=0}^{kr} \sum_{\mathbf{t} \in \pi_s^{(d_s+k)}(\mathcal{N})} a_{\mathbf{t}}c_{k,l}\mathbf{x}^{\mathbf{t} + \mathbf{b}_{s,k,l}^{(2k)}}\bigg)
    \end{eqnarray*}}

Then the equation $f_{\sigma^{-1}(s),d_s-k}(\sigma(\mathbf{x})) = f_{s, d_s+k}(\mathbf{x}) P_{\omega_s}^{k}(\mathbf{x})$ holds, if and only if, relation
    \begin{equation} \label{sum=sumsum_t'}
        \sum_{\mathbf{j} \in \pi_{s}^{(d_s-k)}(\sigma^{-1}(\mathcal{N}))} a_{\sigma(\mathbf{j})}\mathbf{x}^{\mathbf{j} } = \sum_{l=0}^{kr} \sum_{\mathbf{t} \in \pi_s^{(d_s+k)}(\mathcal{N})} a_{\mathbf{t}}c_{k,l}\mathbf{x}^{\mathbf{t} + \mathbf{b}_{s,k,l}^{(2k)}}
    \end{equation}
holds. Denote $S_L$ the set of all the exponent vectors of the terms on the left-hand side of the above equation, and $S_R$ the set of all the exponent vectors of the terms on the right-hand side of the above equation. Clearly, $S_L = \pi_{s}^{(d_s-k)}(\sigma^{-1}(\mathcal{N}))$ and $S_R = \pi_s^{(d_s+k)}(\mathcal{N}) + \cup_{0\leq l \leq kr} \mathbf{b}_{s,k,l}^{(2k)}$. By Equation \eqref{pikB_ki}, for all $l \in [0, kr]$, we know $\pi_s(\mathbf{j}) = d_s-k$, if and only if, $\pi_s(\mathbf{j} - \mathbf{b}_{s,k,l}^{(2k)}) = d_s+k$. Then we have 
\begin{eqnarray*}
    S_L \cap S_R &=& \pi_{s}^{(d_s-k)}\bigg(\sigma^{-1}(\mathcal{N})\bigg)\cap \pi_s^{(d_s-k)}\bigg(\mathcal{N} + \bigcup_{0\leq l \leq kr} \mathbf{b}_{s,k,l}^{(2k)}\bigg)\\
     &=& \pi_{s}^{(d_s-k)}\bigg(\sigma^{-1}(\mathcal{N})\cap \bigcup_{0\leq l \leq kr} (\mathcal{N} + \mathbf{b}_{s,k,l}^{(2k)})\bigg),\\
    S_R \setminus S_L &=& \pi_s^{(d_s-k)}\bigg(\bigcup_{0 \leq l \leq kr} (\mathcal{N} + \mathbf{b}_{s,k,l}^{(2k)}) \setminus \sigma^{-1}(\mathcal{N})\bigg), \\
    S_L \setminus S_R &=& \pi_s^{(d_s-k)}\bigg(\sigma^{-1}(\mathcal{N}) \setminus \bigcup_{0 \leq l \leq kr} (\mathcal{N} + \mathbf{b}_{s,k,l}^{(2k)})\bigg).
\end{eqnarray*}
Hence, by comparing the coefficients of Equation \eqref{sum=sumsum_t'}, it is easy to check that relation \eqref{sum=sumsum_t'} holds, if and only if, the coefficients $a_{\mathbf{j}}$'s satisfy the system of homogeneous linear equations $HLE(\sigma, s, \omega_s, \bm{\eta}, \mathbf{d},  k)$.

    (ii) Under the conditions $\mathbf{d} = \sigma{(\mathbf{d})}$ and $\eta_s = \eta_{t} = 2d_s = 2d_{t}$, by Theorem \ref{thm len1_T2f}, we know that the relation $F(\psi_{\sigma, s, \omega_s}(\mathbf{x})) = F(\mathbf{x})$ holds,  if and only if, for all $i \in [0, \eta_s]$ the following relations
    \begin{equation} \label{ffP}
        f_{t, i}(\sigma(\mathbf{x})) = f_{s,\eta_s-i}(\mathbf{x})P_{\omega_s}^{d_s - i}(\mathbf{x})
    \end{equation}
    hold. Clearly, for all $i \in [0, \eta_s]$, relation \eqref{ffP} hold, if and only if, for $k \in [0, d_s]$, relations 
    \begin{eqnarray}
        f_{t,d_s-k}(\sigma(\mathbf{x})) &=& f_{s,d_s+k}(\mathbf{x})P_{\omega_s}^{k}(\mathbf{x}), \label{ffP1} \\
        f_{s,d_s-k}(\mathbf{x}) &=& f_{t,d_s+k}(\sigma(\mathbf{x}))P_{\omega_s}^{k}(\mathbf{x}).\label{ffP2}
    \end{eqnarray}
    hold.
    
    In Equation \eqref{Ptsigma = Ps ljy} we know that the relation $P_{\omega_t}(\mathbf{x}) = P_{\omega_s}(\sigma^{-1}(\mathbf{x}))$ holds. So we know that relation \eqref{ffP2} holds, if and only if, the relation
        \begin{equation} \label{ffP2'}
            f_{s, d_t-k}(\sigma^{-1}(\mathbf{x})) = f_{t,d_t+k}(\mathbf{x})P_{\omega_t}^{k}(\mathbf{x})
        \end{equation}
    holds. 
    
    In (i), relation \eqref{ffP1} is briefly written as $E(\sigma, s, \omega_s, \bm{\eta}, \mathbf{d},  k)$, then Relation \eqref{ffP2'} can be written as $E(\sigma^{-1}, t, \omega_t, \bm{\eta}, \mathbf{d}, k)$. Hence, for all $i \in [0, \eta_s]$, relation \eqref{ffP} holds, if and only if, for all $k \in [0, d_s]$, the relations $E(\sigma, s, \omega_s, \bm{\eta}, \mathbf{d},  k)$ and $E(\sigma^{-1}, t, \omega_t, \bm{\eta}, \mathbf{d}, k)$ hold.
    
    By (i), we know that the relation $E(\sigma, s, \omega_s, \bm{\eta}, \mathbf{d},  k)$ holds, if and only if, the coefficients $\{a_{\mathbf{j}} \in  \mathbb{Q} \mid \mathbf{j} \in \mathcal{N} \}$ satisfy the system of homogeneous linear equations $HLE(\sigma, s, \omega_s, \bm{\eta}, \mathbf{d},  k)$; the equation $E(\sigma^{-1}, t, \omega_t, \bm{\eta}, \mathbf{d}, k)$ holds, if and only if, the coefficients $\{a_{\mathbf{j}} \in  \mathbb{Q} \mid \mathbf{j} \in \mathcal{N} \}$ satisfy the system of homogeneous linear equations $HLE(\sigma^{-1}, t, \omega_t, \bm{\eta}, \mathbf{d}, k)$. It is thus proved.
\end{proof}

\begin{rmk}\label{rmk HLE}
    (i) We denote that $HLE(\sigma, s, \omega_s, \bm{\eta}, \mathbf{d})$ be the system of homogeneous linear equations containing the system of homogeneous linear equations $HLE(\sigma, s, \omega_s, \bm{\eta}, \mathbf{d},  k)$ and the system of homogeneous linear equations $HLE(\sigma^{-1}, \sigma^{-1}(s), \omega_{\sigma^{-1}(s)}, \bm{\eta}, \mathbf{d},  k)$ for all $k \in [0, d_s]$. That is, 
    $$ HLE(\sigma, s, \omega_s, \bm{\eta}, \mathbf{d}): \left\{ \begin{array}{l}
         HLE(\sigma, s, \omega_s, \bm{\eta}, \mathbf{d},  0),\\
         \qquad\qquad\cdots\cdots\\
         HLE(\sigma, s, \omega_s, \bm{\eta}, \mathbf{d},  d_s), \\
         HLE(\sigma^{-1}, \sigma^{-1}(s), \omega_{\sigma^{-1}(s)}, \bm{\eta}, \mathbf{d},  0),\\
         \qquad\qquad\cdots\cdots\\
         HLE(\sigma^{-1}, \sigma^{-1}(s), \omega_{\sigma^{-1}(s)}, \bm{\eta}, \mathbf{d}, d_s).
    \end{array} \right.$$
    
    (ii) The problem of finding an invariant Laurent polynomial of a given cluster symmetric map $\psi_{\sigma,s,\omega_s}$ is converted to the problem of solving a system of homogeneous linear equations $HLE(\sigma, s, \omega_s, \bm{\eta}, \mathbf{d})$. However, solving this system of equations is tedious. Therefore, we wrtie a MATLAB program attached to Appendix \ref{code ca2diop} so that we can find a cluster symmetric polynomial efficiently and conveniently.  
    
    (iii) Since $F(\mathbf{x})$ is of type $\frac{\bm{\eta}}{\mathbf{d}}$, the coefficients $\{a_{\mathbf{j}} \in  \mathbb{Q} \mid \mathbf{j} \in \mathcal{N} \}$ must satisfy the condition:
    \begin{equation} \label{con eta ne 0}
        \text{for all } i \in \{i \in [1,n] \mid \eta_i \ne 0\}, \text{ there exists } \mathbf{j} \in \pi_i^{(\eta_i)}(\mathcal{N}), \text{ such that } a_{\mathbf{j}} \ne 0.
    \end{equation}
    Hence, once the system of homogeneous linear equations $HLE(\sigma, s, \omega_s, \bm{\eta}, \mathbf{d})$ has been solved, the above conditions must be checked.

    (iv) If Condition \eqref{con eta ne 0} is not checked, then a fundamental solution of $HLE(\sigma, s, \omega_s, \bm{\eta}, \mathbf{d})$ corresponds to an invariant Laurent polynomial of type $\frac{\bm{\eta}'}{\mathbf{d}'}$, where $\eta'_i \leq \eta_i$, $d'_i \leq d_i$ for all $i$. See Example \ref{eg: Gyoda}(i). 
\end{rmk}

\subsection{Examples and practice-level discussion}\label{app: ca2diop}

In this subsection, we apply Theorem \ref{thm len1} to compute several examples and introduce some practice-level propositions.

Some cluster symmetric polynomials are trivial. For example, we consider the cluster symmetric map $\psi_{\sigma_{(12),1,\omega_1}}$ defined in Example \ref{eg CSM somos and marsh}(i). Clearly, the polynomial $F(x_1, x_2, x_3) := x_3$ is a cluster symmetric polynomial about $\psi_{\sigma_{(12),1,\omega_1}}$. However, for $F(\mathbf{x})$, the map serves only as the permutation $\sigma_{(12)}$, there is no substitution of variables here. So, we classify such cluster symmetric polynomials as follows.

\begin{DEF} 
\label{def non-trivial 1-cspoly}
    Given a cluster symmetric polynomial $F(\mathbf{x})$ about $\psi_{\sigma,s,\omega_s}$. Suppose $F(\mathbf{x})$ is of type $\frac{\bm{\eta}}{\mathbf{d}}$. If $\eta_s = 0$, we call $F(\mathbf{x})$ is \textbf{trivial}. If $\eta_s \ne 0$, we call $F(\mathbf{x})$ is \textbf{non-trivial}. 
\end{DEF}

\begin{prop}
    Let $F(\mathbf{x})$ be a trivial cluster symmetric polynomial about cluster symmetric map $\psi_{\sigma, s, \omega_s}$. Then $F(\mathbf{x})$ is invariant under the permutation $\sigma$, that is, $F(\mathbf{x}) \in \mathbb{Q}[\mathbf{x}^{\pm}]^{\langle \sigma \rangle}$.
\end{prop}
\begin{proof}
    Suppose $F(\mathbf{x})$ is of type $\frac{\bm{\eta}}{\mathbf{d}}$. Then $\eta_s = 0$. By Theorem \ref{thm len1_T2f}, we know that $\eta_{\sigma^{-1}(s)} = 0$, $\sigma(\mathbf{d}) = \mathbf{d}$, and $f_{\sigma^{-1}(s),0}(\sigma(\mathbf{x})) = f_{s,0}(\mathbf{x})$. Since the expansions of $F(\mathbf{x})$ are $F(\mathbf{x}) = \mathbf{x}^{-\mathbf{d}} f_{s, 0}(\mathbf{x})$ and $F(\mathbf{x}) = \mathbf{x}^{-\mathbf{d}} f_{\sigma^{-1}(s), 0}(\mathbf{x})$, we know that 
        $F(\sigma(\mathbf{x})) = (\sigma(\mathbf{x}))^{-\mathbf{d}} f_{\sigma^{-1}(s), 0}(\sigma(\mathbf{x}))  =  \mathbf{x}^{-\mathbf{d}} f_{s, 0}(\mathbf{x}) = F(\mathbf{x})$.
    Hence $F(\mathbf{x}) \in \mathbb{Q}[\mathbf{x}^{\pm}]^{\langle \sigma \rangle}$.
\end{proof}

By the above proposition, we only need to consider non-trivial cluster symmetric polynomials. In the end of this subsection, we will provide concrete steps for finding non-trivial cluster symmetric polynomials. To do it, we first consider how to choose the tuple $\mathbf{d}$.

\begin{prop}\label{ddd}
Given a cluster symmetric map $\psi_{\sigma, s, \omega_s}$. For any $i \in [1,n]$, we denote an $n$-tuple $\mathbf{e}_{\sigma,i} := \sum_{j \in \langle \sigma \rangle(i)}\mathbf{e}_j$. If $F(\mathbf{x}) \in \mathbb{Q}[\mathbf{x}^{\pm}]^{\langle \psi_{\sigma, s, \omega_s} \rangle}$ and $i \notin \langle \sigma \rangle(s)$, then $\mathbf{x}^{d\mathbf{e}_{\sigma,i}}F(\mathbf{x}) \in \mathbb{Q}[\mathbf{x}^{\pm}]^{\langle \psi_{\sigma, s, \omega_s} \rangle}$ for any $d \in \mathbb{Z}$. Specifically, suppose $F(\mathbf{x})$ is of type $\frac{\bm{\eta}}{\mathbf{d}}$, then $\mathbf{x}^{\mathbf{d} - d_s\mathbf{e}_{\sigma, s}}F(\mathbf{x}) \in \mathbb{Q}[\mathbf{x}^{\pm}]^{\langle \psi_{\sigma, s, \omega_s} \rangle}$.
\end{prop}
\begin{proof}
    Suppose $F(\mathbf{x}) := \frac{T(\mathbf{x})}{\mathbf{x}^{\mathbf{d}}}$. For $i \notin \langle \sigma \rangle(s)$, let $F'(\mathbf{x}) := \mathbf{x}^{d\mathbf{e}_{\sigma,i}}F(\mathbf{x})$ and $t := \sigma^{-1}(s)$. By Theorem \ref{thm len1_T2f}, we have $\sigma(\mathbf{d}) = \mathbf{d}$. Then it is easy to check that $\sigma(\mathbf{d}-d\mathbf{e}_{\sigma,i}) = \mathbf{d}-d\mathbf{e}_{\sigma,i}$.
    By Equation \eqref{sigma_mu_d}, we have
    $$\frac{T(\mathbf{x})}{\mathbf{x}^{\mathbf{d}}} = F(\mathbf{x}) = F(\psi_{\sigma, s, \omega_s}(\mathbf{x})) = \frac{T(\psi_{\sigma, s, \omega_s}(\mathbf{x}))}{(\psi_{\sigma, s, \omega_s}(\mathbf{x}))^{\mathbf{d}}} = \frac{T(\psi_{\sigma, s, \omega_s}(\mathbf{x}))}{\left(\frac{P_{\omega_s}(\mathbf{x})}{x_s^2}\right)^{d_{t}}\mathbf{x}^{\sigma^{-1}(\mathbf{d})}}.$$
    Hence $T(\psi_{\sigma, s, \omega_s}(\mathbf{x}))\left(\frac{P_{\omega_s}(\mathbf{x})}{x_s^2}\right)^{-d_{t}} = T(\mathbf{x})$. Since $(\psi_{\sigma, s, \omega_s}(\mathbf{x}))^{d\mathbf{e}_{\sigma,i}} = \mathbf{x}^{d\mathbf{e}_{\sigma,i}}$, we have
    \begin{eqnarray*}
        F'(\psi_{\sigma, s, \omega_s}(\mathbf{x})) &=& \frac{T(\psi_{\sigma, s, \omega_s}(\mathbf{x}))}{(\psi_{\sigma, s, \omega_s}(\mathbf{x}))^{\mathbf{d}- d\mathbf{e}_{\sigma,i}}} = \frac{T(\psi_{\sigma, s, \omega_s}(\mathbf{x}))}{\left(\frac{P_{\omega_s}(\mathbf{x})}{x_s^2}\right)^{d_{t}}\mathbf{x}^{\sigma^{-1}(\mathbf{d})-d\mathbf{e}_{\sigma,i}}}\\
            &=& \frac{T(\mathbf{x})}{\mathbf{x}^{\sigma^{-1}(\mathbf{d})-d\mathbf{e}_{\sigma,i}}} = F'(\mathbf{x}).
    \end{eqnarray*}
    That is, $\mathbf{x}^{d\mathbf{e}_{\sigma,i}}F(\mathbf{x}) = F'(\mathbf{x}) \in \mathbb{Q}[\mathbf{x}^{\pm}]^{\langle \psi_{\sigma, s, \omega_s} \rangle}$. Suppose $[1,n] = \langle \sigma \rangle(s)\sqcup \langle \sigma \rangle(i_1)\sqcup \cdots \sqcup \langle \sigma \rangle(i_m)$, where $m \in [1,n-1]$, $i_1, \cdots i_m \in [1,n]\setminus\{s\}$. So
    $$ \mathbf{d} - d_s\mathbf{e}_{\sigma, s} = \sum_{i \in [1,n] \setminus \langle \sigma \rangle(s)} d_i\mathbf{e}_i = \sum_{j=1}^m d_{i_j}\mathbf{e}_{\sigma, i_j}.$$
    Hence, we know $\mathbf{x}^{\mathbf{d} - d_s\mathbf{e}_{\sigma, s}}F(\mathbf{x}) = \mathbf{x}^{\sum_{j=1}^m d_{i_j}\mathbf{e}_{\sigma, i_j}}F(\mathbf{x}) \in \mathbb{Q}[\mathbf{x}^{\pm}]^{\langle \psi_{\sigma, s, \omega_s} \rangle}$.
\end{proof}

In the above proposition, the Laurent polynomial $\mathbf{x}^{\mathbf{d} - d_s\mathbf{e}_{\sigma, s}}F(\mathbf{x})$ is of type $\frac{\bm{\eta}}{d_s\mathbf{e}_{\sigma, s}}$, and by Equation \eqref{len1_eta_s=2d_s}, we know that $d_s = \eta_{s}/2 \geq 0$. Hence, we only need to consider the non-negative $n$-tuple $\mathbf{d} := d_s\mathbf{e}_{\sigma, s}$ with $d_s > 0$.

In the following, we consider examples of the permutation being the identity, that is, $\sigma = id$.

\begin{eg} \label{eg: Gyoda} 
(i) Consider the cluster symmetric map
    $$\psi_{id, 2, \omega_2}(\mathbf{x}) =  \bigg(x_1, \frac{x_1 + x_3^2}{x_2}, x_3\bigg),$$
where the seedlet $\omega_2 := (\mathbf{b}',1, Z')$ defined from Example \ref{eg: 1csm, Gyoda}(i). Let $\bm{\eta} = (1,2,2)$ and $\mathbf{d} = (0,1,0)$.
Then applying Theorem \ref{thm len1} or the corresponding MATLAB program in Appendix \ref{code ca2diop}, we find the solutions of the system of homogeneous linear equations $HLE(id, 2, \omega_2, \bm{\eta}, \mathbf{d})$ as follows 
$$a_{\mathbf{j}} = \begin{cases}
    t_1, & \text{  if } \mathbf{j} \in \{(0,0,2), (0,2,0),(1,0,0)\},\\
    t_2, & \text{  if } \mathbf{j} = (0,1,0).\\
    t_3, & \text{  if } \mathbf{j} = (0,1,1).\\
    t_4, & \text{  if } \mathbf{j} = (0,1,2).\\
    t_5, & \text{  if } \mathbf{j} = (1,1,0).\\
    t_6, & \text{  if } \mathbf{j} = (1,1,1).\\
    t_7, & \text{  if } \mathbf{j} = (1,1,2).\\
    0,   & \text{  otherwise},
\end{cases}$$
where $t_1, \cdots, t_7 \in \mathbb{Q}$ and $t_1 \ne 0$. We denote that 
    $$F_2(\mathbf{x}) := \frac{x_1 + x_2^2 + x_3^2}{x_2} ~\text{ and }~ H_2(\mathbf{x}) := t_2 + t_3x_3 + t_4x_3^2 + t_5x_1 + t_6x_1x_3 + t_7x_1x_3^2.$$
Hence the $\frac{\bm{\eta}}{\mathbf{d}}$ type cluster symmetric polynomial about $\psi_{id, 2, \omega_2}$ is $t_1F_2(\mathbf{x}) + H_2(\mathbf{x}).$

(ii) Consider the cluster symmetric map
    $$\psi_{id, 3, \omega_3}(\mathbf{x}) =  \bigg(x_1, x_2, \frac{x_1 + x_2^2}{x_3}\bigg),$$
where the seedlet $\omega_3 := (\mathbf{b}'',1, Z'')$ defined from Example \ref{eg: 1csm, Gyoda}(ii). Let $\bm{\eta} = (1,2,2)$ and $\mathbf{d} = (0,0,1)$. We denote that
    $$F_3(\mathbf{x}) := F_2(\sigma_{(23)}(\mathbf{x})) = \frac{x_1 + x_2^2 + x_3^2}{x_3} ~\text{ and }~ H_3(\mathbf{x}) := H_2(\sigma_{(23)}(\mathbf{x})).$$
Since $\sigma_{(23)}(\mathbf{b}'') = -\mathbf{b}'$, by Proposition \ref{prop new csp}, we know that the cluster symmetric polynomial of type $\frac{\bm{\eta}}{\mathbf{d}}$ about $\psi_{id, 2, \omega_3}$ is $t_1F_3(\mathbf{x}) + H_3(\mathbf{x})$.

(iii) Consider the cluster symmetric map
    $$\psi_{id, 1, \omega_1}(\mathbf{x}) =  \bigg( \frac{k_0x_2^4+k_{1}x_2^3x_3+k_{2}x_2^2x_3^2+k_{3}x_2x_3^3+k_4x_3^4}{x_1}, x_2, x_3\bigg),$$
where the seedlet $\omega_1 := (\mathbf{b}, 1, Z)$ defined from Example \ref{eg: 1csm, Gyoda}(iii). Let $\bm{\eta} = (2,4,4)$ and $\mathbf{d} = (1,0,0)$. 
Based on the results of running the MATLAB program in Appendix \ref{code ca2diop}, we denote a Laurent polynomial
    $$F_1(\mathbf{x}) := \frac{x_1^2 + k_0x_2^4+k_{1}x_2^3x_3+k_{2}x_2^2x_3^2+k_{3}x_2x_3^3+k_4x_3^4}{x_1}.$$
Then the $\frac{\bm{\eta}}{\mathbf{d}}$ type cluster symmetric polynomial about $\psi_{id, 1, \omega_1}$ is
$$aF_1(\mathbf{x}) + H_1(\mathbf{x}),$$
where $a \in \mathbb{Q}_{\ne 0}$ and $H_1(\mathbf{x})$ is a polynomial in $\mathbb{Q}[\mathbf{x}]$ with $\deg^1H_1(\mathbf{x}) = 0, \deg^2 H_1(\mathbf{x}) \leq 4$ and $\deg^3H_1(\mathbf{x}) \leq 4$.
\end{eg}

Note that the cluster symmetric polynomials in the above examples can be written as 
$$t_1 \frac{P_{\omega_s}(\mathbf{x}) + x_s^2}{x_s} + H(\mathbf{x}),$$
where $\deg^s H(\mathbf{x}) = 0$. This is due to the following proposition.

\begin{prop} \label{F in ring: id}
    Given a cluster symmetric map $\psi_{id, s, \omega_s}$. Let $F(\mathbf{x}) := \frac{T(\mathbf{x})}{\mathbf{x}^{\mathbf{d}}}$ be a $\frac{\bm{\eta}}{\mathbf{d}}$ type Laurent polynomial in $\mathbb{Q}[\mathbf{x}^{\pm}]^{\langle \psi_{id, s, \omega_s} \rangle}$. Then
    \begin{equation} \label{eq T/d in Q}
        \frac{T(\mathbf{x})}{\mathbf{x}^{d_s\mathbf{e}_s}} \in \mathbb{Q}\bigg[\frac{P_{\omega_s}(\mathbf{x}) + x_s^2}{x_s}, \mathbf{x}\setminus x_s\bigg], 
    \end{equation}
    and the invariant Laurent polynomial ring $\mathbb{Q}[\mathbf{x}^{\pm}]^{\langle \psi_{id, s, \omega_s} \rangle}$ is the polynomial ring in  $\frac{P_{\omega_s}(\mathbf{x}) + x_s^2}{x_s}$ and $x_i^{\pm}$ for all $i \ne s$, that is,
    \begin{equation}\label{eq Q^psi = Q[P,x]}
        \mathbb{Q}[\mathbf{x}^{\pm}]^{\langle \psi_{id, s, \omega_s} \rangle} = \mathbb{Q}\bigg[\frac{P_{\omega_s}(\mathbf{x}) + x_s^2}{x_s}, \mathbf{x}^{\pm}\setminus x_{s}^{\pm}\bigg].
    \end{equation}
\end{prop}
\begin{proof} 
    Suppose $T(\mathbf{x}) = \sum_{i=0}^{\eta_s} f_{s,i}(\mathbf{x})x_s^{i-d_s}$, where $f_{s,i}(\mathbf{x})$ is defined in Equation \eqref{f_sum}. Then by Equation \eqref{len1_f=fp@s}, we have
    \begin{eqnarray*}
    \frac{T(\mathbf{x})}{\mathbf{x}^{d_s\mathbf{e}_s}} &=& \sum_{i=0}^{\eta_s} f_{s,i}(\mathbf{x})x_s^{i-d_s}\\
    &=& f_{s,d_s}(\mathbf{x}) + \sum_{i=0}^{d_s-1} f_{s,i}(\mathbf{x})x_s^{i-d_s} + \sum_{i=d_s+1}^{\eta_s} f_{s,i}(\mathbf{x})x_s^{i-d_s}\\
    &=& f_{s,d_s}(\mathbf{x}) + \sum_{i=0}^{d_s-1} f_{s,\eta_s-i}(\mathbf{x}) \bigg(\frac{P_{\omega_s}(\mathbf{x})}{x_s}\bigg)^{d_s-i} + \sum_{i=d_s+1}^{\eta_s} f_{s,i}(\mathbf{x})x_s^{i-d_s}\\
    &=& f_{s,d_s}(\mathbf{x}) + \sum_{i=0}^{d_s-1} f_{s,\eta_s-i}(\mathbf{x}) \bigg(\bigg(\frac{P_{\omega_s}(\mathbf{x})}{x_s}\bigg)^{d_s-i} + x_s^{d_s-i}\bigg)\bigg).
    \end{eqnarray*}

Let $H_i(u,v) := u^{d_s - i} + v^{d_s - i}$. 
Since $H_i(u,v)$ is a symmetric polynomial, 
by the fundamental theorem on symmetric polynomials (Theorem \ref{FTSP}), we know that there exists $\widetilde{H}_i(u,v) \in \mathbb{Q}[u,v]$, such that
$H_i(u,v) = \widetilde{H}_i(S_{2,1}(u,v), S_{2,2}(u,v))$.
Then
\begin{eqnarray*}
    \bigg(\frac{P_{\omega_s}(\mathbf{x})}{x_s}\bigg)^{d_s-i} + x_s^{d_s-i} &=& \widetilde{H}_i\bigg(S_{2,1}(\frac{P_{\omega_s}(\mathbf{x})}{x_s}, x_s), S_{2,2}(\frac{P_{\omega_s}(\mathbf{x})}{x_s},x_s)\bigg)\\
        &=& \widetilde{H}_i\bigg(\frac{P_{\omega_s}(\mathbf{x}) + x_s^2}{x_s}, P_{\omega_s}(\mathbf{x})\bigg)
\end{eqnarray*}
and 
\begin{equation*}
    \frac{T(\mathbf{x})}{\mathbf{x}^{d_s\mathbf{e}_s}} =  f_{s,d_s}(\mathbf{x}) + \sum_{i=0}^{d_s} f_{s,\eta_s-i}(\mathbf{x})\widetilde{H}_i\bigg(\frac{P_{\omega_s}(\mathbf{x}) + x_s^2}{x_s}, P_{\omega_s}(\mathbf{x})\bigg).
\end{equation*}
By Equations \eqref{lem_degf} and \eqref{deg_p_br}, we know that
$$\deg^s f_{s,d_s}(\mathbf{x}) = \deg^s f_{s,\eta_s-i}(\mathbf{x}) = 0 \text{ and } \deg^s P_{\omega_s}(\mathbf{x}) = r|b_s| = 0.$$
Hence $\frac{T(\mathbf{x})}{\mathbf{x}^{d_s\mathbf{e}_s}} \in \mathbb{Q}\bigg[x_1, \cdots, x_{s-1}, \frac{P_{\omega_s}(\mathbf{x}) + x_s^2}{x_s}, x_{s+1}, \cdots, x_n \bigg]$ and $$F(\mathbf{x}) = \frac{1}{\mathbf{x}^{\mathbf{d}-d_s\mathbf{e}_s}}\frac{T(\mathbf{x})}{\mathbf{x}^{d_s\mathbf{e}_s}} \in \mathbb{Q}\bigg[x_1^{\pm}, \cdots, x_{s-1}^{\pm}, \frac{P_{\omega_s}(\mathbf{x}) + x_s^2}{x_s}, x_{s+1}^{\pm}, \cdots, x_n^{\pm}\bigg].$$
So we have $\mathbb{Q}[\mathbf{x}^{\pm}]^{\langle \psi_{id, s, \omega_s} \rangle} \subset \mathbb{Q}\bigg[\frac{P_{\omega_s}(\mathbf{x}) + x_s^2}{x_s}, \mathbf{x}^{\pm}\setminus x_{s}^{\pm}\bigg]$. Clearly, $\frac{P_{\omega_s}(\mathbf{x}) + x_s^2}{x_s}, x_i^{\pm}\in \mathbb{Q}[\mathbf{x}^{\pm}]^{\langle \psi_{id, s, \omega_s} \rangle}$, where $i \in [1,n]$ and $i \ne s$. Hence, Equation \eqref{eq Q^psi = Q[P,x]} holds.
\end{proof}

Next, we show some examples of different $\bm{\eta}$.

\begin{eg} \label{somos5}
    Consider the cluster symmetric map
    $$\psi_{\sigma_{(12345)}, 1, \omega_1}(\mathbf{x}) = \bigg(x_2, x_3, x_4, x_5, \frac{\tilde{\alpha}x_2x_5 + \tilde{\beta}x_3x_4}{x_1}\bigg),$$
    where the seedlet $\omega_1 := (\mathbf{b},1, Z)$ defined from Example \ref{eg CSM somos and marsh}(iv). This map is related to the Somos 5 sequence\cite{HoneSigma}. Let $\mathbf{d} = (1,1,1,1,1)$.
    
    (i) When $\bm{\eta} = (2,2,2,2,2)$. After computing, there are no $\frac{\bm{\eta}}{\mathbf{d}}$ type cluster symmetric polynomial about the map $\psi_{\sigma_{(12345)}, 1, \omega_1}$.    
    
    (ii) When $\bm{\eta} = (2,2,3,2,2)$. After computing, we denote that
        $$F_1(\mathbf{x}) := \frac{  x_1x_2^2x_5^2 + x_1^2x_4^2x_5 +\tilde{\alpha}(x_1x_3^2x_4^2  + x_2^2x_3^2x_5) + \tilde{\beta}x_2x_3^3x_4}{x_1x_2x_3x_4x_5}.$$
    Then the $\frac{\bm{\eta}}{\mathbf{d}}$ type cluster symmetric polynomial about $\psi_{\sigma_{(12345)}, 1, \omega_1}$ is 
        $t_1F_1(\mathbf{x}) + t_0$, 
    where $t_1, t_0 \in \mathbb{Q}$ and $t_1 \ne 0$.
    
    (iii) When $\bm{\eta} = (2,3,3,3,2)$. After computing, we denote that
        $$F_2(\mathbf{x}) := \frac{x_1^2x_3x_5^2 + \tilde{\alpha}(x_1x_2x_4^3 + x_1x_3^3x_5 + x_2^3x_4x_5) + \tilde{\beta}x_2^2x_3x_4^2}{x_1x_2x_3x_4x_5}.$$
    Then the $\frac{\bm{\eta}}{\mathbf{d}}$ type cluster symmetric polynomial about $\psi_{\sigma_{(12345)}, 1, \omega_1}$ is 
        $q_2F_2(\mathbf{x}) + q_1F_1(\mathbf{x}) + q_0$
    where $q_2, q_1, q_0 \in \mathbb{Q}$ and $q_2 \ne 0$. 
\end{eg}

The above example shows that for a fixed $\mathbf{d}$, different $\bm{\eta}$ will give different results. However, $\bm{\eta}$ is an arbitrary non-negative $n$-tuple except that it satisfies the relation $\eta_s = \eta_{\sigma^{-1}(s)} = 2d_s = 2d_{\sigma^{-1}(s)}$. How can we further restrict the range of $\bm{\eta}$? We have the following proposition.

\begin{prop}\label{d and n}
    Given a seedlet $\omega_s := (\mathbf{b}, r, Z)$ and a cluster symmetric map $\psi_{\sigma, s, \omega_s}$. Let $F(\mathbf{x})$ be a Laurent polynomial of type $\frac{\bm{\eta}}{\mathbf{d}}$ in $\mathbb{Q}[\mathbf{x}^{\pm}]$. Suppose that the equation 
    $F(\psi_{\sigma, s, \omega_s}(\mathbf{x})) = F(\mathbf{x})$ holds. Then $\bm{\eta}$ and $\mathbf{d}$ satisfy $\eta_s = \eta_{\sigma^{-1}(s)} = 2d_s = 2d_{\sigma^{-1}(s)}, \sigma(\mathbf{d}) = \mathbf{d}$ and
    \begin{equation}
        2\min\{\eta_k, \eta_{\sigma^{-1}(k)}\} \geq \eta_sr|b_{k}| \geq 2|\eta_k - \eta_{\sigma^{-1}(k)}|,\label{eta_k>=drb}
    \end{equation}
    for all $k \in [1,n]$.
\end{prop}
\begin{proof} Let $t := \sigma^{-1}(s)$. By Theorem \ref{thm len1_T2f}, we have $\eta_s = \eta_{\sigma^{-1}(s)} = 2d_s = 2d_{\sigma^{-1}(s)}$ and $\sigma(\mathbf{d}) = \mathbf{d}$. If $k = s$, since $b_{s} = 0$, we have $\eta_s, \eta_{t} \geq 0 = |\eta_s - \eta_t|$.
If $k \ne s$. Suppose that the expansions of $F(\mathbf{x})$ are
    $$F(\mathbf{x}) = \mathbf{x}^{-\mathbf{d}} \sum_{i=0}^{\eta_s} f_{s, i}(\mathbf{x}) \mathbf{x}^{i\mathbf{e}_s} = \mathbf{x}^{-\mathbf{d}} \sum_{i=0}^{\eta_t} f_{t, i}(\mathbf{x}) \mathbf{x}^{i\mathbf{e}_t}$$
as shown in Equation \eqref{F=sumf}.
    
Since Theorem \ref{thm len1_T2f}, the following relations
    $$f_{t,i}(\sigma(\mathbf{x})) = f_{s,\eta_s-i}(\mathbf{x})P_{\omega_s}^{d_s - i}(\mathbf{x}),~\forall~i \in [0, \eta_s]$$ 
hold. We apply the function $\deg^k$ to the above equation, then by equations \eqref{deg^k_sigma} and \eqref{deg_p_br}, we have 
    \begin{equation}
        \deg^{\sigma^{-1}(k)} f_{t,i}(\mathbf{x}) = \deg^k f_{s,\eta_s-i}(\mathbf{x}) + (d_s - i)r|b_{k}|. \label{degt=degs+drb}
    \end{equation}
    
    Observing the above equation, on the one hand, by the equations \eqref{deg_f} and \eqref{len1_eta_s=2d_s}, for $i \in [0, \eta_s]$ we have
    \begin{eqnarray*}
        \deg^k f_{s,i}(\mathbf{x}) &\in& [0, \eta_k]\cap[ (d_s - i)r|b_{k}|, \eta_{\sigma^{-1}(k)} + (d_s - i)r|b_{k}|],\\
        \deg^{\sigma^{-1}(k)} f_{t,i}(\mathbf{x}) &\in& [0, \eta_{\sigma^{-1}(k)}]\cap[(d_s - i)r|b_{k}|, \eta_{k} + (d_s - i)r|b_{k}|].
    \end{eqnarray*}
    Since the sets on the right side of the above relations are not empty, we have
    $$\eta_k \geq r|d_s - i||b_{k}|  \quad \text{and}\quad \eta_{\sigma^{-1}(k)} \geq r|d_s - i||b_{k}|$$
    for all  $i \in [0,\eta_s]$. Then 
    $ \eta_k \geq d_sr|b_{k}| $ and $\eta_{\sigma^{-1}(k)} \geq d_sr|b_{k}|.$ So we have $2\min\{\eta_k, \eta_{\sigma^{-1}(k)}\} \geq \eta_sr|b_{k}|$. 
    
    On the other hand, by Relation \eqref{deg^k=nk}, there exist $i_k,j_k \in [0, \eta_s]$, such that
        $$\deg^k f_{s, \eta_s -i_k}(\mathbf{x}) = \eta_k \quad \text{and}\quad \deg^{\sigma^{-1}(k)} f_{t, j_k}(\mathbf{x}) = \eta_{\sigma^{-1}(k)}.$$
    Then by Equation \eqref{degt=degs+drb}, we have
    
    \centerline{$\eta_k  - d_sr|b_{k}| \leq \eta_k + (d_s - i_k)r|b_{k}| = \deg^{\sigma^{-1}(k)} f_{t,i_k}(\mathbf{x}) \leq \eta_{\sigma^{-1}(k)},$}
    
    \centerline{$\eta_{\sigma^{-1}(k)} = \deg^k f_{s,\eta_s-i}(\mathbf{x}) + (d_s - j_k)r|b_{k}|  \leq \eta_k + (d_s - j_k)r|b_{k}| \leq \eta_k + d_sr|b_{k}|.$}
    
    Therefore, $d_sr|b_{k}| \geq |\eta_k - \eta_{\sigma^{-1}(k)}|$ and $\eta_sr|b_{k}| \geq 2|\eta_k - \eta_{\sigma^{-1}(k)}|$.
\end{proof}

If there exists a non-trivial cluster symmetric polynomial that is about two cluster symmetric maps, then the conditions that the two maps need to satisfy are immediately known by the above proposition.
\begin{cor} \label{cor 4 eta}
Given two cluster symmetric maps $\psi_{\sigma, s, \omega_s}$, $\psi_{\tau, s', \omega_{s'}}$, where $\omega_s := (\mathbf{b}, r, Z)$ and $\omega_{s'} := (\mathbf{b}', r', Z')$. Let $F(\mathbf{x})$ be a Laurent polynomial of type $\frac{\bm{\eta}}{\mathbf{d}}$ in $\mathbb{Q}[\mathbf{x}^{\pm}]^{\langle \psi_{\sigma, s, \omega_s}, \psi_{\tau, s', \omega_{s'}} \rangle}$. If $\eta_s \ne 0$, then 
    \begin{equation} \label{eqn 4>rrbb}
        4 \geq rr'\max\{|b_{s'}|,|b_{\sigma(s')}|\}\max\{|b'_s|,|b'_{\tau(s)}|,|b'_{\sigma^{-1}(s)}|,|b'_{\tau(\sigma^{-1}(s))}|\}.
    \end{equation}
\end{cor}
\begin{proof}
    By Proposition \ref{d and n}, for all $k \in [1,n]$, we have
    $$2\eta_k \geq \eta_sr\max\{|b_k|,|b_{\sigma(k)}|\} \quad \text{ and } \quad  2\eta_k \geq \eta_{s'}r'\max\{|b'_k|,|b'_{\tau(k)}|\}.$$
    Then we have
    $$4\eta_k \geq 2\eta_{s'}r'\max\{|b'_k|,|b'_{\tau(k)}|\} \geq \eta_srr'\max\{|b_{s'}|,|b_{\sigma(s')}|\}\max\{|b'_k|,|b'_{\tau(k)}|\}.$$
    Taking $k = s$ and $k = \sigma^{-1}(s)$, we have
    \begin{eqnarray*}
        4\eta_s &\geq& \eta_srr'\max\{|b_{s'}|,|b_{\sigma(s')}|\}\max\{|b'_s|,|b'_{\tau(s)}|\},\\
        4\eta_{\sigma^{-1}(s)} &\geq&\eta_srr'\max\{|b_{s'}|,|b_{\sigma(s')}|\}\max\{|b'_{\sigma^{-1}(s)}|,|b'_{\tau(\sigma^{-1}(s))}|\}.
    \end{eqnarray*}
    Since $\eta_s = \eta_{\sigma^{-1}(s)} \ne 0$, we know Relation \eqref{eqn 4>rrbb} holds.
\end{proof}

Finally, we show some examples of different $\mathbf{d}$ and $\bm{\eta}$.

\begin{eg} 
    Consider the cluster symmetric map
    $$\psi_{\sigma_{(123)}, 1, \omega_1}(\mathbf{x}) = \bigg(x_2, x_3, \frac{1 + x_2x_3}{x_1}\bigg),$$ 
    where the seedlet $\omega_1 := (\mathbf{b},1, Z)$ defined from Example \ref{eg CSM somos and marsh}(ii). This map was studied by Fordy and Marsh in \cite{FordyMarsh} and is related to the primitive period 1 quiver. After computing, 
    we denote that
        \begin{eqnarray*}
            F_1(\mathbf{x}) &:=& \frac{x_{1}^2x_{3}+x_{1}x_{2}^2+x_{1}x_{3}^2+x_{2}^2x_{3}+x_{2}}{x_1x_2x_3},\\
            F_2(\mathbf{x}) &:=& \frac{x_{2}x_{1}^2+x_{1}+x_{2}x_{3}^2+x_{3}}{x_1x_2x_3}.
        \end{eqnarray*}
    
    (i) For $\mathbf{d} = (1,1,1)$, using Theorem \ref{thm len1} and Proposition \ref{d and n}, we can check that only when $\bm{\eta} = (2,2,2)$ there exists a cluster symmetric polynomial. The $\frac{\bm{\eta}}{\mathbf{d}}$ type cluster symmetric polynomial about the cluster symmetric map $\psi_{\sigma_{(123)}, 1, \omega_1}$ is
        $$t_1F_1(\mathbf{x}) + t_2F_2(\mathbf{x}) + t_3$$
    where $t_1, t_2, t_3 \in \mathbb{Q}$ and $t_1 \ne 0$.

    (ii) For $\mathbf{d} = (2,2,2)$, we can check that only when $\bm{\eta} = (4,4,4)$ there exists a cluster symmetric polynomial. The $\frac{\bm{\eta}}{\mathbf{d}}$ type cluster symmetric polynomial about $\psi_{\sigma_{(123)}, 1, \omega_1}$ is
        $$t_1F_1(\mathbf{x})^2 + t_2F_1(\mathbf{x})F_2(\mathbf{x}) + t_3F_3(\mathbf{x})^2 + t_4F_1(\mathbf{x}) + t_5F_2(\mathbf{x}) + t_6$$
    where $t_1, \cdots, t_6 \in \mathbb{Q}$ and $t_1 \ne 0$.

    (iii)  For $\mathbf{d} = (3,3,3)$, we can check that only when $\bm{\eta} = (6,6,6)$ there exists a cluster symmetric polynomial. The $\frac{\bm{\eta}}{\mathbf{d}}$ type cluster symmetric polynomial about $\psi_{\sigma_{(123)}, 1, \omega_1}$ is
    $H(F_1(\mathbf{x}), F_2(\mathbf{x}))$,
    where $H(u,v)$ is a polynomial with $\deg^1 H(u,v) = 3$ and $\deg^2 H(u,v) \leq 3$.
    
\end{eg}

\begin{rmk} \label{rmk step of constuct}
To summarize this section, find a non-trivial cluster symmetric polynomial about a given cluster symmetric map $\psi_{\sigma,s,\omega_s}$ in the following steps:
\begin{enumerate}
    \item Choose a $n$-tuple $\mathbf{d}$.
    
    By Proposition \ref{ddd}, we only need to consider the non-negative $n$-tuple $\mathbf{d} := d\mathbf{e}_{\sigma,i} = d\sum_{j \in \langle \sigma \rangle(s)}\mathbf{e}_j$, where $d \geq 0$. By Definition \ref{def non-trivial 1-cspoly} and Equation \eqref{len1_eta_s=2d_s}, the number $d$ should be a positive integer.
    \item Choose a $n$-tuple $\bm{\eta}$. 
    
    By Proposition \ref{d and n}, the tuple $\bm{\eta}$ should satisfy two conditions $\eta_s = \eta_{\sigma^{-1}(s)} = 2d$ and $\min\{\eta_k, \eta_{\sigma^{-1}(k)}\} \geq dr|b_{k}| \geq |\eta_k - \eta_{\sigma^{-1}(k)}|$ for all $k \in [1,n]$. 
        
    \item Solve the system $HLE(\sigma, s, \omega_s, \bm{\eta}, \mathbf{d})$. 
    
    Applying the MATLAB program in Appendix \ref{code ca2diop}, we obtain the solutions of the homogeneous linear equation system $HLE(\sigma, s, \omega_s, \bm{\eta}, \mathbf{d})$. A fundamental solution of $HLE(\sigma, s, \omega_s, \bm{\eta}, \mathbf{d})$ corresponds to an invariant Laurent polynomial of type $\frac{\bm{\eta}'}{\mathbf{d}'}$, where $\eta'_i \leq \eta_i$, $d'_i \leq d_i$ for all $i \in [1,n]$.
\end{enumerate}
\end{rmk}

\section{Cluster symmetric maps and generalized cluster algebras}\label{sec: ca}

In this section, we first set up the notion of a cluster symmetric map of a seed, similarly to that in the case of data. Here, a given cluster symmetric map of a seed is abstracted from the composition of a permutation and a mutation coming from a generalized cluster algebra. However,  we will see that not all such compositions are cluster symmetric maps. We will discuss when this is true. In the end, we will answer two questions posed by Gyoda and Matsushita in \cite{Gyoda}.

\subsection{Generalized cluster algebra}
In this subsection, we recall some definitions and theorems of the generalized cluster algebra \cite{ChekhovShapiro2014, Nakanishi2015}. We fix a positive integer $n$. Let $x_1, \dots, x_n$ be indeterminates and $\mathcal{F} := \mathbb{Q}(x_1, \dots, x_n)$, we call $\mathcal{F}$ \textbf{ambient field}. We first define the seed.

\begin{DEF}\label{def seed} 
A \textbf{seed} in $\mathcal{F}$ is a quadruplet $\Omega := (B, \mathbf{x}, R, \mathbf{Z})$, where
\begin{itemize}
    \item $B = \left(b_{i j}\right)$ is an $n \times n$ integer skew-symmetrizable matrix, called an \textbf{exchange matrix};
    \item $\mathbf{x}=\left(x_1, \ldots, x_n\right)$ is an $n$-tuple such that $\left\{x_1, \ldots, x_n\right\}$ is a free generating set of $\mathcal{F}$. We call $\mathbf{x}$ the \textbf{cluster} and $x_1, \ldots, x_n$ the \textbf{cluster variables} of $\Omega$;
    \item $R = \operatorname{diag}\left(r_1, \cdots, r_n\right)$ is a diagonal integer matrix with $r_i>0$, called a \textbf{mutation degree matrix};
    \item $\mathbf{Z} = (Z_1, \cdots, Z_n)$ is an $n$-tuple of polynomials, where for $k \in [1,n]$, 
        $$Z_k(u) := \sum_{i=0}^{r_k} z_{k,i}u^i =  z_{k,0} + z_{k,1}u+ \cdots + z_{k,r_k}u^{r_k} \in \mathbb{ Z}_{\geq 0}[u]$$
    satisfying the reciprocity condition 
        \begin{equation}\label{recip}
            z_{k, t} = z_{k, r_k-t} \text{ for } t \in \{1, \cdots, r_k-1\}
        \end{equation}
    and $z_{k, 0} = z_{k, r_k} = 1$.
    We call $\mathbf{Z}$ the \textbf{mutation polynomial tuple} and $Z_1, \cdots, Z_n$ the \textbf{mutation polynomials} of the seed $\Omega$.
\end{itemize}
\end{DEF}

\begin{rmk}
An integer matrix $B_{n \times n}$ is \textbf{skew-symmetrizable} if there is a positive integer diagonal matrix $S$ such that $SB$ is skew-symmetric. This $S$ is said to be a \textbf{skew-symmetrizer} of $B$. A positive integer diagonal matrix $S$ is said to be a \textbf{skew-symmetrizer} of the seed $\Omega = (B, \mathbf{x}, R, \mathbf{Z})$, if $S$ is a skew-symmetrizer of $BR$, that is, $SBR = -(SBR)^T.$
\end{rmk}

A seed can induce $n$ seedlets. The following property is trivial.
\begin{pty} \label{seed to seedlet}
    Given a seed $\Omega := (B, \mathbf{x}, R, \mathbf{Z})$. For $s \in [1,n]$, we denote a map as
        $$\pi_s(\Omega) := (B_s, r_s, Z_s),$$
    where $B_s := (b_{1s}, \cdots, b_{ns})$ be the transpose of the $s$-th column of the matrix $B$. Then the triplet $\pi_s(\Omega)$ is a seedlet at direction $s$.
\end{pty}


There are two types of transformation of seed, mutation and permutation. We first define the mutation.
\begin{DEF} \label{def mu}
    Let $\Omega := (B, \mathbf{x}, R, \mathbf{Z})$ be a seed. The \textbf{mutation of the  seed $\Omega$ at direction} $s \in [1,n]$ is defined to be the new seed 
        $\mu_s(B, \mathbf{x}, R, \mathbf{Z}) := (\mu_s(B), \mu_s(\mathbf{x}), \mu_s(R), \mu_s(\mathbf{Z})) := (B', \mathbf{x}', R, \mathbf{Z})$
    given by    
\begin{eqnarray*}
    b_{i j}^{\prime} &=& \begin{cases}-b_{i j}, & \text{if } i=s \text { or } j=s, \\ 
    b_{i j} + r_s\left(\left[b_{i s}\right]_{+} b_{s j}+b_{i s}\left[-b_{s j}\right]_{+}\right), & \text {otherwise.}\end{cases}\\
    x_j^{\prime} &=& \begin{cases} x_s^{-1}P_{\Omega,s}(\mathbf{x}), & ~~\qquad\qquad\qquad\qquad \text{if } j=s,\\ x_j, & ~\qquad\qquad\qquad\qquad \text{ otherwise, }\end{cases}
\end{eqnarray*}
where $P_{\Omega,s}(\mathbf{x}) \in \mathbb{Z}_{\geq 0}[\mathbf{x}]$ is the \textbf{exchange polynomial of $\Omega$ at direction} $s$ defined by 
    $$P_{\Omega,s}(\mathbf{x}) :=  \mathbf{x}^{r_s[-B_s]_{+}}Z_{s}(\mathbf{x}^{B_s}) = \sum_{i=0}^{r_s} z_{s,i}\mathbf{x}^{i[B_s]_{+} + (r_s - i)[-B_s]_{+}}$$
and the \textbf{exchange polynomial tuple of the seed} $\Omega$  is defined as
    $$\mathbf{P}(\Omega) := (P_{\Omega,1}, \cdots, P_{\Omega, n}).$$
\end{DEF}

\begin{rmk}\label{rmk POmega}
    (i) Given a seed $\Omega := (B, \mathbf{x}, R, \mathbf{Z})$, we denote a seed $\Omega^{-} := (-B, \mathbf{x}, R, \mathbf{Z})$. Then by Condition \eqref{recip}, for all $s \in [1,n]$ we have
    \begin{eqnarray*}
        P_{\Omega^{-},s}(\mathbf{x}) &=& \sum_{i=0}^{r_s} z_{s,i}\mathbf{x}^{i[-B_s]_{+} + (r_s - i)[B_s]_{+}}\\
            &=& \sum_{j=0}^{r_s} z_{s,r_s-j}\mathbf{x}^{(r_s - j)[-B_s]_{+} + j[B_s]_{+}}\\
            &=& \sum_{j=0}^{r_s} z_{s,j} \mathbf{x}^{j[B_s]_{+}+(r_s - j)[-B_s]_{+} }\\
            &=& P_{\Omega,s}(\mathbf{x}).
    \end{eqnarray*}
Hence, the following two exchange polynomial tuples coincide
    \begin{equation} \label{P(-B)=P(B)}
        \mathbf{P}(\Omega) = \mathbf{P}(\Omega^{-}).
    \end{equation}

    (ii) The exchange polynomial of $\Omega$ at direction $s$ is the same as the exchange polynomial of the seedlet $\pi_s(\Omega)$, that is, $P_{\Omega,s}(\mathbf{x}) = P_{\pi_s(\Omega)}(\mathbf{x})$.

    (iii) When $R = I_n$, then $\mathbf{Z} = (1+u, \cdots, 1+u)$ and for $s \in [1,n]$ we have
        $$ P_s(B, \mathbf{x}, I_n, \mathbf{Z}) =  \mathbf{x}^{[B_k]_{+}} + \mathbf{x}^{[-B_k]_{+}} = \prod_{i=1}^n x_i^{\left[-b_{i k}\right]_{+}} +  \prod_{i=1}^n x_i^{[b_{i k}]_{+}}.$$
    Hence, the seed $(B, \mathbf{x}, I_n, \mathbf{Z})$ is the classic seed, the mutation is the classic mutation defined by Fomin and Zelevinsky in \cite{FZ1}.
\end{rmk}

The second transformation of the seed is permutation.
\begin{DEF}
    Let $\Omega := (B, \mathbf{x}, R, \mathbf{Z})$ be a seed. For any permutation $\sigma \in \mathfrak{S}_n$. The \textbf{permutation} $\sigma$ of the seed $\Omega$ is defined to be the new seed 
    $$\sigma(B, \mathbf{x}, R, \mathbf{Z}) := (\sigma(B), \sigma(\mathbf{x}), \sigma(R), \sigma(\mathbf{Z})) := (B', \mathbf{x}', R', \mathbf{Z}'),$$
    where $ b_{ij}' = b_{\sigma(i)\sigma(j)}, x_i' = x_{\sigma(i)}, r_i' = r_{\sigma(i)}, Z_i' = Z_{\sigma(i)}.$
\end{DEF}

We have the following property to facilitate the computation of the composite of a permutation and a mutation.

\begin{pty} 
    Given a seed $\Omega := (B, \mathbf{x}, R, \mathbf{Z})$, for any $\sigma \in \mathfrak{S}_n, k \in [1,n]$, we have
    \begin{eqnarray} \label{musigma=sigmamu}
        \sigma\mu_k(\Omega) &=& \mu_{\sigma^{-1}(k)}\sigma(\Omega). 
    \end{eqnarray}
\end{pty}

Now we give the definition of generalized cluster algebras. 

\begin{DEF} \label{def geneCA}
    For any two seeds $\Omega := ({B}, {\mathbf{x}}, R, \mathbf{Z})$, $\Omega' := ({B}', {\mathbf{x}}', R', \mathbf{Z}')$, if there exists a finite-length sequence of mutations $\mu_{s_1}, \cdots, \mu_{s_m}$, such that $\mu_{s_m}\cdots \mu_{s_1}(\Omega) = \Omega'$, then we call the two seeds $\Omega$ and $\Omega'$ is \textbf{mutation equivalent}, denoted as $\Omega \sim \Omega'$. Let
\begin{equation*}
	\mathcal{X}(\Omega) = \mathcal{X}({B}, {\mathbf{x}}, R, \mathbf{Z}) := \bigcup_{({B}, {\mathbf{x}}, R, \mathbf{Z}) \sim ({B}', {\mathbf{x}}', R', \mathbf{Z}')} \{x_1', \dots, x_n'\}
\end{equation*}
be the set of cluster variables for all seeds that are mutation equivalent to $\Omega$. The $\mathbb{Q}$-subalgebra generated by $\mathcal{X}(\Omega)$ of the
ambient field $\mathcal{F}$ is the \textbf{generalized cluster algebra}, we denote it as $\mathcal{A}(\Omega)$.
\end{DEF}

One of the main results of generalized cluster algebras is the positive Laurent phenomenon. That is, after arbitrarily mutating an initial cluster, the resulting new cluster variables can always be expressed as a Laurent polynomial of the initial cluster variables, and the coefficients of the Laurent polynomials are positive. We restate this in our notation.

\begin{thm}[(Positive Laurent phenomenon {\cite[Theorem 5.8]{moulang2024}})]
\label{thm PLP}
    Given a seed $\Omega := (B, \mathbf{x}, R, \mathbf{Z})$. Let $\mathbf{x}' := \mu_{s_m} \cdots \mu_{s_1}(\mathbf{x})$, where $s_1, \cdots, s_m \in [1,n], m \in \mathbb{Z}_{\geq 0}$. Then $x_i' \in \mathbb{Z}_{\geq 0}[\mathbf{x}^{\pm}]$ for all $i \in [1,n]$. 
\end{thm}

\subsection{Cluster symmetric maps of a seed}

In general, the exchange polynomial tuple $\mathbf{P}(\Omega)$ may not be preserved under permutations or mutations. 
For example, we consider the seed $\Omega := (B, \mathbf{x}, I_3, \mathbf{Z})$, where 
    \begin{center}
    $B := \left[
     	\begin{array}{rrr}
     	 0 & 1 & 1 \\
     	-1 & 0 & 0 \\
     	-1 & 0 & 0 \\
     	\end{array}
    \right]$.
    \end{center}
Then $\mathbf{P}(\Omega) = (P_{\Omega,1}, P_{\Omega,2}, P_{\Omega,3})$, where
    $$P_{\Omega,1}(\mathbf{x}) = \mathbf{x}^{(0,1,1)} + 1, \quad P_{\Omega,2}(\mathbf{x}) = \mathbf{x}^{(1,0,0)} + 1, \quad P_{\Omega,3}(\mathbf{x}) = \mathbf{x}^{(1,0,0)} + 1.$$
    
Let $\mu_2(\Omega) = (\bar{B}, \bar{\mathbf{x}}, I_3, \mathbf{Z})$, where
\begin{center}
    $\bar{B} := \left[
     	\begin{array}{rrr}
     	 0 & -1 & 1 \\
     	 1 & 0 & 0 \\
     	-1 & 0 & 0 \\
     	\end{array}
    \right]$.
    \end{center}
Then $\mathbf{P}(\mu_2(\Omega)) = (P_{\mu_2(\Omega),1}, P_{\mu_2(\Omega),2}, P_{\mu_2(\Omega),3})$, where
    $$P_{\mu_2(\Omega),1}(\bar{\mathbf{x}}) = \bar{\mathbf{x}}^{(0,1,0)} + \bar{\mathbf{x}}^{(0,0,1)}, P_{\mu_2(\Omega),2}(\bar{\mathbf{x}}) = \bar{\mathbf{x}}^{(1,0,0)} + 1, P_{\mu_2(\Omega),3}(\bar{\mathbf{x}}) = \bar{\mathbf{x}}^{(1,0,0)} + 1.$$
Since $P_{\Omega,1} \ne P_{\mu_2(\Omega),1}$, we have $\mathbf{P}(\Omega) \ne \mathbf{P}(\mu_2(\Omega))$. 

However, under some special actions, the exchange polynomial tuple will be preserved. For example, let $\mu_1(\Omega) = (\tilde{B}, \tilde{\mathbf{x}}, I_3, \mathbf{Z})$, where
    \begin{center}
    $\tilde{B} := \left[
     	\begin{array}{rrr}
     	 0 & -1 & -1 \\
     	 1 & 0 & 0 \\
     	 1 & 0 & 0 \\
     	\end{array}
    \right].$
    \end{center}
Then $\mathbf{P}(\mu_1(\Omega)) = (P_{\mu_1(\Omega),1}, P_{\mu_1(\Omega),2}, P_{\mu_1(\Omega),3})$, where
    $$P_{\mu_1(\Omega),1}(\tilde{\mathbf{x}}) = \tilde{\mathbf{x}}^{(0,1,1)} + 1, \quad P_{\mu_1(\Omega),2}(\tilde{\mathbf{x}}) = \tilde{\mathbf{x}}^{(1,0,0)} + 1, \quad P_{\mu_1(\Omega),3}(\tilde{\mathbf{x}}) = \tilde{\mathbf{x}}^{(1,0,0)} + 1.$$
Hence we have $\mathbf{P}(\Omega) = \mathbf{P}(\mu_1(\Omega))$. And we can also check $\mathbf{P}(\Omega) = \mathbf{P}(\sigma_{(23)}(\Omega))$.

From the above observations, we define a group that can preserve the exchange polynomial tuple under the permutations or mutations.

\begin{prop} \label{P=Pg}
Given a seed $\Omega := (B, \mathbf{x}, R, \mathbf{Z})$. Denote $\overline{\mathcal{G}}(\Omega)$ be the set
    $$\{ g := \sigma\mu_{s_m}\cdots\mu_{s_1} \mid g(B, \mathbf{x}, R, \mathbf{Z}) = (\pm B, \mathbf{x}', R, \mathbf{Z}), \sigma \in \mathfrak{S}_n, m \geq 0, s_i \in [1,n] \}.$$
Then we have

    (i) $\overline{\mathcal{G}}(B, \mathbf{x}, R, \mathbf{Z}) = \overline{\mathcal{G}}(-B, \mathbf{x}, R, \mathbf{Z}).$ 
    
    (ii) $\overline{\mathcal{G}}(\Omega)$ is a group. We call $\overline{\mathcal{G}}(\Omega)$ the \textbf{complete cluster symmetric group of the seed} $\Omega$.

    (iii) The action of the complete cluster symmetric group $\overline{\mathcal{G}}(\Omega)$ preserves the exchange polynomial tuple of $\Omega$, that is, 
    for any $g \in \overline{\mathcal{G}}(\Omega)$, we have
        $\mathbf{P}(\Omega) = \mathbf{P}(g(\Omega))$.
\end{prop}
\begin{proof} (i) Let $h \in \overline{\mathcal{G}}(\Omega)$ with $h = \sigma \in \mathfrak{S}_n$ or $h = \mu_k $ for some $k \in [1,n]$. Denote $(B', \mathbf{x}', R', \mathbf{Z}') := h(B, \mathbf{x}, R, \mathbf{Z})$. It is easy to check that $h(-B, \mathbf{x}, R, \mathbf{Z}) = (-B', \mathbf{x}', R', \mathbf{Z}')$. Then for $g \in \overline{\mathcal{G}}(B, \mathbf{x}, R, \mathbf{Z})$, we have 
$$g(-B, \mathbf{x}, R, \mathbf{Z}) = (\mp B, \mathbf{x}, R, \mathbf{Z}),$$ 
that is, $g \in \overline{\mathcal{G}}(-B, \mathbf{x}, R, \mathbf{Z})$.

(ii) Let $id$ be the identity of $\mathfrak{S}_n$, then we have $id(\Omega) = \Omega$. So $id \in \overline{\mathcal{G}}(\Omega)$.

Let $g := \sigma\mu_{s_m}\cdots\mu_{s_1}, g_1 := \tau\mu_{t_p}\cdots\mu_{t_1} \in \overline{\mathcal{G}}(\Omega)$, then by Equation \eqref{musigma=sigmamu}, we have
    $$gg_1 = (\sigma\mu_{s_m}\cdots\mu_{s_1})(\tau\mu_{t_p}\cdots\mu_{t_1}) = \sigma\tau\mu_{\tau(s_m)}\cdots\mu_{\tau(s_1)}\mu_{t_p}\cdots\mu_{t_1}$$
and
    $$gg_1(B, \mathbf{x}, R, \mathbf{Z}) = \begin{cases}
        g(B, \mathbf{x}', R, \mathbf{Z})\\
        g(-B, \mathbf{x}', R, \mathbf{Z})
    \end{cases} = \begin{cases}
        (\pm B, \mathbf{x}', R, \mathbf{Z})\\
        (\mp B, \mathbf{x}', R, \mathbf{Z})
    \end{cases},$$
where the last equality is by (i). So $gg_1 \in \overline{\mathcal{G}}(\Omega)$.
    
Let $g' := \sigma^{-1}\mu_{\sigma^{-1}(s_1)}\cdots\mu_{\sigma^{-1}(s_m)}$. Then by Equation \eqref{musigma=sigmamu}, we have
    \begin{eqnarray*}
        gg' &=& (\sigma\mu_{s_m}\cdots\mu_{s_1})(\sigma^{-1}\mu_{\sigma^{-1}(s_1)}\cdots\mu_{\sigma^{-1}(s_m)})\\
            &=& \sigma\sigma^{-1}\mu_{\sigma^{-1}(s_m)}\cdots\mu_{\sigma^{-1}(s_1)}\mu_{\sigma^{-1}(s_1)}\cdots\mu_{\sigma^{-1}(s_m)} = id
    \end{eqnarray*}
and 
    $$g'(B, \mathbf{x}, R, \mathbf{Z}) = g'(g(\pm B, \mathbf{x}', R, \mathbf{Z})) = id(\pm B, \mathbf{x}', R, \mathbf{Z}) = (\pm B, \mathbf{x}', R, \mathbf{Z}).$$
So $g$ has an inverse $g' \in \overline{\mathcal{G}}(\Omega)$. Hence $\overline{\mathcal{G}}(\Omega)$ is a group.

(iii) For any $g \in \overline{\mathcal{G}}(\Omega)$. Let $g(\Omega) = (B', \bar{\mathbf{x}}, R, \mathbf{Z})$ where $B' = \pm B$. Fix $s \in [1,n]$. If $B' = B$. Then $P_{\Omega,s}(\mathbf{x}) =  \mathbf{x}^{r_s[-B_s]_{+}}Z_{s}(\mathbf{x}^{B_s})$ and $P_{\Omega,s}(\mathbf{x}') =  (\mathbf{x}')^{r_s[-B_s]_{+}}Z_{s}((\mathbf{x}')^{B_s})$. So $P_{\Omega,s} = P_{g(\Omega),s}$. If $B' = -B$. By Remark \ref{rmk POmega}(i), we know $P_{\Omega,s} = P_{g(\Omega),s}$.    
\end{proof}

\begin{rmk}
    (i) Although any action $g$ of the group $\overline{\mathcal{G}}(\Omega)$ is a transformation between seeds, according to the above proposition, the action $g$ can be regarded as a transformation between clusters, that is, $g(\mathbf{x}) = \mathbf{x}'$. 

    (ii) The complete cluster symmetric group is a subgroup of the mutation group defined by King and Pressland\cite {King2017}. The mutation-periodic group of an exchange matrix defined by Liu and Li\cite{LiuLi} is a subgroup of the complete cluster symmetric group.
\end{rmk}

In practice, the complete cluster symmetric group is not easy to describe, but we can easily calculate some subset of it. 

\begin{DEF} \label{def csm}
Given a seed $\Omega$. The \textbf{cluster symmetric set of the seed $\Omega$} is defined as
    $$\mathcal{S}(\Omega) := \{ \sigma\mu_{s} \mid \sigma\mu_{s}(B, \mathbf{x}, R, \mathbf{Z}) = (\pm B, \mathbf{x}', R, \mathbf{Z}), \sigma \in \mathfrak{S}_n, s \in [1,n] \}.$$
The \textbf{cluster symmetric group of the seed $\Omega$} be the group ${\mathcal{G}}(\Omega)$ generated by the set $\mathcal{S}(\Omega)$. The element in the set $\mathcal{S}(\Omega)$ is called the \textbf{cluster symmetric map of the seed} $\Omega$.
\end{DEF}

Following this definition, clearly, the cluster symmetric group ${\mathcal{G}}(\Omega)$ is a subgroup of the complete cluster symmetric group $\overline{\mathcal{G}}(\Omega)$. 

What is the relationship between the cluster symmetric maps of the seed $\Omega$ and the cluster symmetric map of the data defined in Definition \ref{seedlet}(iii)? The following proposition answers: a cluster symmetric map of a seed $\Omega$ is a cluster symmetric map of a data. We begin with a lemma.

\begin{lem} \label{lem Bt Bs}
    Given a seed $\Omega := (B, \mathbf{x}, R, \mathbf{Z})$. Suppose $\sigma\mu_s \in \mathcal{S}(\Omega)$. Let $t := \sigma^{-1}(s)$. Then
    
    (i) $\sigma^{-1}\mu_t \in \mathcal{S}(\Omega)$.

    (ii) $B_t = \pm \sigma(B_s)$, where $B_k := (b_{1k}, \dots, b_{nk})$.

    (iii) $\pi_t(\Omega) = (\sigma(B_s), r_s, Z_s)$ or $\pi_t(\Omega) = (-\sigma(B_s), r_s, Z_s)$.
\end{lem}
\begin{proof}
    Since $\sigma\mu_s(B, \mathbf{x}, R, \mathbf{Z}) = (\pm B, \mathbf{x}', R, \mathbf{Z})$, we know that $\sigma\mu_s(B) = \pm B$, $r_s = r_t$ and $Z_s = Z_t$. Since $\sigma^{-1}\mu_{t}(B) = \sigma^{-1}\mu_{t}(\pm \sigma\mu_s(B)) = \pm \sigma^{-1}\sigma\mu_{s}\mu_s(B)= \pm B,$
    we have $\sigma^{-1}\mu_t \in \mathcal{S}(\Omega)$ and $\mu_t(B) = \pm \sigma(B)$. Considering the transpose of the $t$-th column of matrices of two sides of the equation $\mu_t(B) = \pm (b_{\sigma(i)\sigma(j)})$, we have $-B_t = \pm (b_{\sigma(1)s}, \cdots, b_{\sigma(n)s}) = \pm \sigma(B_s)$.
    So $\pi_t(\Omega) = (B_t, r_t, Z_t) = (\pm \sigma(B_s), r_s, Z_s)$.
\end{proof}

\begin{prop}\label{is CSM}
Given a seed $\Omega := (B, \mathbf{x}, R, \mathbf{Z})$. 
Any cluster symmetric map of the seed $\Omega$, when treated as a transformation of the cluster, is a cluster symmetric map of a data defined in Definition \ref{seedlet}. That is, suppose $\sigma\mu_s \in \mathcal{S}(\Omega)$, then 
    $$\psi^m_{\sigma, s, \pi_s(\Omega)}(\mathbf{x}) = (\sigma\mu_s)^m(\mathbf{x}), \quad \text{ for all } m \in \mathbb{Z},$$
where the map $\pi_s$ is defined in Property \ref{seed to seedlet}. So $\psi_{\sigma, s, \pi_s(\Omega)} \in \mathcal{S}(\Omega)$.
\end{prop}
\begin{proof} 
(i) For $m>0$, we prove it by induction on $m$. When $m = 1$, by Remark \ref{rmk POmega}(iii), we have
    $$\psi_{\sigma, s, \pi_s(\Omega)}(\mathbf{x}) = \bigg(\sigma(\mathbf{x})\bigg)\bigg|_{x_s \gets \frac{P_{\pi_s(\Omega)}(\mathbf{x})}{x_s}} = \bigg(\sigma(\mathbf{x})\bigg)\bigg|_{x_s \gets \frac{P_{\Omega,s}(\mathbf{x})}{x_s}} = \sigma\mu_s(\mathbf{x}).$$
Assume it is true for $m = k-1$. Let $\mathbf{y} := \psi^{k-1}_{\sigma, s, \pi_s(\Omega)}(\mathbf{x})$. By Proposition \ref{P=Pg}(iii), we have
    $$\sigma\mu_s(\mathbf{y}) = \bigg(\sigma(\mathbf{y})\bigg)\bigg|_{y_s \gets \frac{P_{(\sigma\mu_s)^{k-1}(\Omega),s}(\mathbf{y})}{y_s}} = \bigg(\sigma(\mathbf{y})\bigg)\bigg|_{y_s \gets \frac{P_{\Omega,s}(\mathbf{y})}{y_s}} = \psi_{\sigma, s, \pi_s(\Omega)}(\mathbf{y}).$$

(ii) For $m < 0$, we prove it by induction in $m$. When $m = -1$, let $t := \sigma^{-1}(s), \omega'_t := (\sigma(B_s), r_s, Z_s)$. Then by Equation \eqref{psi inv}, we have
    $$\psi^{-1}_{\sigma, s, \pi_s(\Omega)}(\mathbf{x}) = \psi_{\sigma^{-1}, t, \omega'_t}(\mathbf{x}).$$
By Lemma \ref{lem Bt Bs}, we have $\omega'_t = (B_t, r_t, Z_t) = \pi_t(\Omega)$ or $\omega'_t = (-B_t, r_t, Z_t) = \pi_t(\Omega^{-})$, where $\Omega^{-} := (-B, \mathbf{x}, R, \mathbf{Z})$. By Equation \eqref{P(-B)=P(B)}, we have $P_{\omega'_t} = P_{\pi_t(\Omega)} = P_{\pi_t(\Omega^{-})}$. Hence
$$\psi_{\sigma^{-1}, t, \omega'_t}(\mathbf{x}) = \bigg(\sigma^{-1}(\mathbf{x})\bigg)\bigg|_{x_t \gets \frac{P_{\pi_t(\Omega)}(\mathbf{x})}{x_t} } = \bigg(\sigma^{-1}(\mathbf{x})\bigg)\bigg|_{x_t \gets \frac{P_{\Omega,t}(\mathbf{x})}{x_t}} = \sigma^{-1}\mu_t(\mathbf{x}).$$
Assume that it is true for $m = k+1$. Let $\mathbf{y} := \psi^{k+1}_{\sigma, s, \pi_s(\Omega)}(\mathbf{x})$. By Proposition \ref{P=Pg}(iii), we have
    \begin{eqnarray*}
        \sigma^{-1}\mu_t(\mathbf{y}) &=& \bigg(\sigma^{-1}(\mathbf{y})\bigg)\bigg|_{y_t \gets \frac{P_{(\sigma^{-1}\mu_t)^{k+1}(\Omega),t}(\mathbf{y})}{y_t}}\\
        &=& \bigg(\sigma^{-1}(\mathbf{y})\bigg)\bigg|_{y_t \gets \frac{P_{\Omega,t}(\mathbf{y})}{y_t}} = \psi_{\sigma^{-1}, t, \pi_t(\Omega)}(\mathbf{y}).
    \end{eqnarray*}
\end{proof}

When the cluster symmetric set $\mathcal{S}(\Omega)$ is not empty, by the above proposition, we know that the seed $\Omega$ can correspond to cluster symmetric maps. Conversely, when can a cluster symmetric map correspond to a seed? We give the following definition and property.

\begin{DEF} \label{def seed and 1-csm}
    Given a cluster symmetric map $\psi_{\sigma,s,\omega_s}$. 
    If there exists a seed $\Omega := (B, \mathbf{x}, R, \mathbf{Z})$, such that $\sigma\mu_s \in \mathcal{S}(\Omega)$ and $\omega_s = \pi_s(\Omega^{\pm})$ where $\Omega^{\pm} := (\pm B, \mathbf{x}, R, \mathbf{Z})$, then we call $\psi_{\sigma,s,\omega_s}$ \textbf{corresponds to} the seed $\Omega$ and the seed $\Omega$ \textbf{corresponds to} the map $\psi_{\sigma,s,\omega_s}$. In this situation, by Proposition \ref{is CSM}, we know that $\psi^m_{\sigma,s,\omega_s}(\mathbf{x}) = (\sigma\mu_s)^m(\mathbf{x})$ for $m \in \mathbb{Z}$, so $\psi_{\sigma,s,\omega_s} \in \mathcal{S}(\Omega)$.
\end{DEF}

\begin{rmk}
    In general, $\pi_s(\Omega^+) \ne \pi_s(\Omega^-)$, but we have $\psi_{\sigma,s,\pi_s(\Omega^+)} = \psi_{\sigma,s, \pi_s(\Omega^-)}$ since Proposition \ref{pty ex-poly}(i) and Condition \eqref{recip}. So we require $\omega_s = \pi_s(\Omega^{\pm})$ instead of $\omega_s = \pi_s(\Omega)$ in the definition.
\end{rmk}

\begin{pty}
    Given a seed $\Omega := (B, \mathbf{x}, R, \mathbf{Z})$ with a nonempty cluster symmetric set $\mathcal{S}(\Omega)$. If $\sigma\mu_s \in \mathcal{S}(\Omega)$, then the cluster symmetric map $\psi_{\sigma,s,\pi_s(\Omega)}$ corresponds to the seed $\Omega$ and its inverse cluster symmetric map $\psi^{-1}_{\sigma,s,\pi_s(\Omega)}$ corresponds to the seed $\Omega$.
\end{pty}
\begin{proof}
    Obviously, $\psi_{\sigma,s,\pi_s(\Omega)}$ corresponds to the seed $\Omega$. 
    By Equation \eqref{psi inv}, we have $\psi^{-1}_{\sigma,s,\pi_s(\Omega)} = \psi_{\sigma^{-1},t,\omega_t}$, where $t := \sigma^{-1}(s)$ and $\omega_t := (\sigma(B_s), r_s, Z_s)$.
    We claim that $\psi_{\sigma^{-1},t,\omega_t}$ corresponds to the seed $\Omega$. Since by Lemma \ref{lem Bt Bs}, we have $\sigma^{-1}\mu_t \in \mathcal{S}(\Omega)$ and $\pi_t(\Omega) = (\pm \sigma(B_s), r_s, Z_s)$ which implies $\omega_t = \pi_t(\Omega^{\pm})$.
\end{proof}

We show some examples.

\begin{eg} \label{eg 1-csm not to seed}
    Denote a seedlet $\omega_1 = ((0,1,-2,1),1,a+bu)$ with $(a, b) \ne (1, 1)$. Since the polynomial $a+bu$ is not a mutation polynomial defined in Definition \ref{def seed}, there does not exist a seed $\Omega$, such that $\omega_1 = \pi_1(\Omega^{\pm})$. So for any $\sigma \in \mathfrak{S}_4$, the cluster symmetric map $\psi_{\sigma,1,\omega_1}$ does not correspond to any seeds.
\end{eg}

\begin{eg} \label{eg 1-csm to seed} 
    Denote a seedlet $\omega_1 = ((0,1,-2,1),1,1+u)$ and a seed $\Omega := (B, \mathbf{x}, R, \mathbf{Z})$ where
    
    \begin{equation*}
        B =  \begin{bmatrix}
        		0 & -1 & 2 & -1\\
			    1 &  0 & b_{23} & b_{24}\\
			   -2 & b_{32} & 0 & b_{34}\\
                1 & b_{42} & b_{43} &  0
			\end{bmatrix}, R = \begin{bmatrix}
        		1 &  &  & \\
			    & r_2 & &  \\
			    &  & r_3 & \\
                &  &  &  r_4  
			\end{bmatrix}, \begin{array}{l}
        Z_1(u) = 1 + u, \\
        Z_2(u) = \sum_{i=0}^{r_2}z_{2,i}u^i, \\
        Z_3(u) = \sum_{i=0}^{r_3}z_{3,i}u^i, \\
        Z_4(u) = \sum_{i=0}^{r_4}z_{4,i}u^i. \end{array}
    \end{equation*}
    It is clearly that $\pi_1(\Omega) = \omega_1$ and
    $$\mu_1(B) = \begin{bmatrix}
        		0 & 1 & -2 & 1\\
			    -1 &  0 & b_{23}+2 & b_{24}\\
			    2 & b_{32}-2 & 0 & b_{34}-2\\
                -1 & b_{42} & b_{43}+2 &  0
			\end{bmatrix}.$$
    
    (i) We consider the cluster symmetric map $\psi_{\sigma_{(24)},1,\omega_1}$. Assume $\sigma_{(24)}\mu_1 \in \mathcal{S}(\Omega)$, we know that
    \begin{equation*}
        B = \begin{bmatrix}
        		0 & -1 & 2 & -1\\
			    1 &  0 & -c & -d\\
			   -2 &  c & 0 & 2-c\\
                1 &  d & c-2 &  0
			\end{bmatrix}, R = \begin{bmatrix}
        		1 &  &  & \\
			    & r_2 & &  \\
			    &  & r_3 & \\
                &  &  &  r_2  
			\end{bmatrix}, \begin{array}{l}
        Z_1(u) = 1 + u, \\
        Z_2(u) = \sum_{i=0}^{r_2}z_{2,i}u^i, \\
        Z_3(u) = \sum_{i=0}^{r_3}z_{3,i}u^i, \\
        Z_4(u) = Z_2(u), \end{array}
    \end{equation*}
    where $c,d \in \mathbb{Z}$. Hence $\psi_{\sigma_{(24)},1,\omega_1}$ corresponds to the seed $\Omega$.

    (ii) It is easy to check that $\sigma_{(12)}\mu_1(B, \mathbf{x}, R, \mathbf{Z}) \ne (\pm B, \mathbf{x}', R, \mathbf{Z})$, then the cluster symmetric map $\psi_{\sigma_{(12)},1,\omega_1}$ does not correspond to any seeds.

    (iii) We consider the cluster symmetric map $\psi_{\sigma_{id},1,\omega_1}$. Assume $\mu_1 \in \mathcal{S}(\Omega)$, we have
    $$B = \begin{bmatrix}
        		0 & -1 & 2 & -1\\
			    1 &  0 & -1 & 0\\
			    -2 & 1 & 0 & 1\\
                1 & 0 & -1 &  0
			\end{bmatrix}, R = \begin{bmatrix}
        		1 &  &  & \\
			    & r_2 & &  \\
			    &  & r_3 & \\
                 &  &  &  r_4  
			\end{bmatrix}, \begin{array}{l}
        Z_1(u) = 1 + u, \\
        Z_2(u) = \sum_{i=0}^{r_2}z_{2,i}u^i, \\
        Z_3(u) = \sum_{i=0}^{r_3}z_{3,i}u^i, \\
        Z_4(u) = \sum_{i=0}^{r_4}z_{4,i}u^i. \end{array}$$
    Then $\psi_{\sigma_{id},1,\omega_1}$ corresponds to the seed $\Omega$.
    
    
    (iv) We consider the cluster symmetric map $\psi_{\sigma_{(1234)},1,\omega_1}$. Assume $\sigma_{(1234)}\mu_1 \in \mathcal{S}(\Omega)$, we have
    \begin{equation*}
        B = \begin{bmatrix}
        		0 & -1 & 2 & -1\\
			    1 &  0 & -3 & 2\\
			   -2 &  3 & 0 & -1\\
                1 &  -2 & 1 &  0
			\end{bmatrix}, ~ R = \begin{bmatrix}
        		1 &  &  & \\
			    & 1 & &  \\
			    &  & 1 & \\
                &  &  &  1  
			\end{bmatrix}, ~ \begin{array}{l}
        Z_1(u) = 1 + u, \\
        Z_2(u) = 1 + u, \\
        Z_3(u) = 1 + u, \\
        Z_4(u) = 1 + u. \end{array}
    \end{equation*}
    Then $\psi_{\sigma_{(1234)},1,\omega_1}$ corresponds to the seed $\Omega$. Let $\omega_4 = ((-1,2,-1,0),1,1+u)$ be a seedlet. It is easy to check that $\psi_{\sigma_{(13)},4,\omega_4}$ and $\psi_{\sigma_{(1234)},4,\omega_4}$ also corresponds to the seed $\Omega$.
\end{eg}

For the symmetric group and the symmetric polynomial, there is a well-known theorem, the fundamental theorem on symmetric polynomials.

\begin{thm}[(Fundamental theorem on symmetric polynomials\cite{Hungerford})]\label{FTSP} The set of all symmetric polynomials in $\mathbb{Q}[\mathbf{x}]$ is the polynomial ring in $S_{n,1}(\mathbf{x}), \cdots, S_{n,n}(\mathbf{x})$, that is,
    $$\mathbb{Q}[\mathbf{x}]^{\mathfrak{S}_n}[S_{n,1}(\mathbf{x}), \cdots, S_{n,n}(\mathbf{x})]$$
where $S_{n,i}$'s are the \textbf{elementary symmetric polynomials} of $n$ variables $\mathbf{x}$, that is, $S_{n,1}(\mathbf{x}) := x_1+ \cdots + x_n$, $S_{n,2}(\mathbf{x}) := x_1x_2 + \cdots x_{n-1}x_n, \cdots, S_{n,n}(\mathbf{x}) := x_1 \cdots x_n$.
\end{thm}

For the cluster symmetric group of a certain seed, its invariant Laurent polynomial ring has a similar structure.

\begin{prop}\label{prop len0 B=0}
    Given a seed $\Omega_0 := (B, \mathbf{x}, rI_n, \mathbf{Z})$, where $B = (0)_{n \times n}, Z_1(u)=\cdots=Z_n(u)$, $r \in \mathbb{Z}_{>0}$. Define a map $\varphi(\mathbf{x}) = (\varphi_1(\mathbf{x}), \dots, \varphi_n(\mathbf{x})),$ where $\varphi_k(\mathbf{x}) := \frac{x_k^2 + c}{x_k}$ for $k \in [1,n]$ and $c:= Z_1(1)$. Then the invariant Laurent polynomial ring $\mathbb{Q}[\mathbf{x}^{\pm}]^{\mathcal{G}(\Omega_0)}$ is the polynomial ring in $S_{n,1}(\varphi(\mathbf{x})), \cdots, S_{n,n}(\varphi(\mathbf{x}))$, that is, 
    \begin{equation} \label{eq ring omega0}
        \mathbb{Q}[\mathbf{x}^{\pm}]^{\mathcal{G}(\Omega_0)} = \mathbb{Q}[S_{n,1}(\varphi(\mathbf{x})), \cdots, S_{n,n}(\varphi(\mathbf{x}))].
    \end{equation}
\end{prop}
\begin{proof} It is clear that $\mathcal{S}(\Omega_0) = \{ \sigma\mu_i \mid \sigma \in \mathfrak{S}_n, i \in [1,n]\}$, $\mathbf{P}(\Omega_0)=(c,\dots,c)$ and $c \geq 2$. We claim that for $\sigma\mu_s \in \mathcal{S}(\Omega_0)$, the following relation holds,
    $$\varphi(\sigma\mu_s(\mathbf{x})) = \sigma(\varphi(\mathbf{x})).$$
It is true, since 
\begin{eqnarray*}
    \varphi_k(\sigma\mu_s(\mathbf{x})) = \begin{cases}
     \frac{({c}/{x_s})^2 + c}{{c}/{x_s}}, &\text{if } \sigma(k) = s,\\
     \frac{x_{\sigma(k)}^2 + c}{x_{\sigma(k)}}, &\text{if } \sigma(k) \ne s.
\end{cases}
    = \begin{cases}
    \varphi_s(\mathbf{x}), &\text{if } \sigma(k) = s,\\
    \varphi_{\sigma(k)}(\mathbf{x}), &\text{if } \sigma(k) \ne s.
\end{cases}
    = \varphi_{\sigma(k)}(\mathbf{x}),
\end{eqnarray*}
and $\varphi(\sigma\mu_s(\mathbf{x})) = (\varphi_1(\sigma\mu_s(\mathbf{x})), \cdots, \varphi_n(\sigma\mu_s(\mathbf{x})) = (\varphi_{\sigma(1)}(\mathbf{x}), \cdots, \varphi_{\sigma(n)}(\mathbf{x})) = \sigma(\varphi(\mathbf{x}))$.

($\supset$): Fix $k \in [1,n]$. For $\sigma\mu_s \in \mathcal{S}(\Omega_0)$, we have $S_{n,k}(\varphi(\sigma\mu_s(\mathbf{x}))) = S_{n,k}(\sigma(\varphi(\mathbf{x}))) = S_{n,k}(\varphi(\mathbf{x}))$. Hence $S_{n,k}(\varphi(\mathbf{x})) \in \mathbb{Q}[\mathbf{x}^{\pm}]^{\mathcal{G}(\Omega_0)}$.

($\subset$): Let $F_1(\mathbf{x}) := \frac{T_1(\mathbf{x})}{\mathbf{x}^{\mathbf{d}}}$ be a Laurent polynomial of type $\frac{\bm{\eta}}{\mathbf{d}}$ in $\mathbb{Q}[\mathbf{x}^{\pm}]^{\mathcal{G}(\Omega_0)}$. Since $F_1(\mathbf{x}) \in \mathbb{Q}[\mathbf{x}^{\pm}]^{\langle \mu_1 \rangle}$, by Equation \eqref{eq T/d in Q} in Theorem \ref{F in ring: id}, there exists a polynomial $T_2(\mathbf{x}) \in \mathbb{Q}[\mathbf{x}]$, such that
    $$\frac{T(\mathbf{x})}{\mathbf{x}^{d_1\mathbf{e}_1}} = T_2 \bigg(\frac{P_{\pi_1(\Omega_0)}(\mathbf{x}) + x_1^2}{x_1}, x_2, \cdots, x_n \bigg) = T_2(\varphi_1(\mathbf{x}), x_2, \cdots, x_n).$$
Let $\mathbf{x}_{(2)} := (\varphi_1(\mathbf{x}), x_2, \cdots, x_n)$ and $F_2(\mathbf{x}_{(2)}) := \frac{T_2(\mathbf{x}_{(2)})}{\mathbf{x}_{(2)}^{\mathbf{d} - d_1\mathbf{e}_1}}$. Since $F_1(\mathbf{x}) \in \mathbb{Q}[\mathbf{x}^{\pm}]^{\langle \mu_2 \rangle}$ and $\varphi_1(\mu_2(\mathbf{x}) = \varphi_1(\mathbf{x}) = \pi_1(\mu_2(\mathbf{x}_{(2)}))$, we have
\begin{eqnarray*}
  F_2(\mathbf{x}_{(2)}) &=& F_1(\mathbf{x})\\
                        &=& F_1(\mu_2(\mathbf{x}))\\
                        &=& \frac{T_2(\varphi_1(\mu_2(\mathbf{x})), \pi_2(\mu_2(\mathbf{x})), \cdots, \pi_n(\mu_2(\mathbf{x})))}{(\mu_2(\mathbf{x}))^{\mathbf{d} - d_1\mathbf{e}_1}}\\
                        &=& \frac{T_2(\pi_1(\mu_2(\mathbf{x}_{(2)})), \pi_2(\mu_2(\mathbf{x}_{(2)})), \cdots, \pi_n(\mu_2(\mathbf{x}_{(2)})))}{(\mu_2(\mathbf{x}_{(2)}))^{\mathbf{d} - d_1\mathbf{e}_1}}\\
                        &=& F_2(\mu_2(\mathbf{x}_{(2)})).
\end{eqnarray*}
Hence $F_2(\mathbf{x}_{(2)}) \in \mathbb{Q}[\mathbf{x}_{(2)}^{\pm}]^{\langle \mu_2 \rangle}$. By Equation \eqref{eq T/d in Q} in Theorem \ref{F in ring: id}, there exists a polynomial $T_3(\mathbf{x}) \in \mathbb{Q}[\mathbf{x}]$, such that
        $$\frac{T_2(\mathbf{x}_{(2)})}{\mathbf{x}_{(2)}^{d_2\mathbf{e}_2}} = T_3(x_{(2),1}, \varphi_2(\mathbf{x}_{(2)}), x_{{(2)},2}, \cdots, x_{{(2)},n}) = T_3(\varphi_1(\mathbf{x}), \varphi_2(\mathbf{x}), x_3, \cdots, x_n).$$
Repeating the above steps, we can find the polynomials  $T_4(\mathbf{x}), \cdots, T_n(\mathbf{x}) \in \mathbb{Q}[\mathbf{x}]$ in order, such that 
    $$F(\mathbf{x}) = \frac{T_2(\mathbf{x}_{(2)})}{\mathbf{x}_{(2)}^{\mathbf{d} - d_1\mathbf{e}_1}} = \frac{T_3(\mathbf{x}_{(3)})}{\mathbf{x}_{(3)}^{\mathbf{d} - d_1\mathbf{e}_1 - d_2\mathbf{e}_3}} = \cdots = \frac{T_{n-1}(\mathbf{x}_{({n-1})})}{\mathbf{x}_{({n-1})}^{\mathbf{d} - \sum_{i=1}^{n-1}d_i\mathbf{e}_i}} = T_n(\mathbf{x}_{(n)}),$$
where $\mathbf{x}_{(k)} := (\varphi_1(\mathbf{x}), \varphi_2(\mathbf{x}), \cdots,\varphi_{k}(\mathbf{x}), x_{k+1}, \cdots, x_n)$. Hence
    $$F(\mathbf{x}) = T_n(\mathbf{x}_{(n)}) = T_n(\varphi_1(\mathbf{x}), \cdots, \varphi_n(\mathbf{x})) = T_n(\varphi(\mathbf{x})).$$
For $\sigma \in \mathfrak{S}_n$, we have $T_n(\sigma(\varphi(\mathbf{x}))) = T_n(\varphi(\sigma\mu_1(\mathbf{x}))) = F(\sigma\mu_1(\mathbf{x})) = F(\mathbf{x}) = T_n(\varphi(\mathbf{x}))$. So, by Theorem \ref{FTSP}, there exists $H(\mathbf{x}) \in \mathbb{Q}[\mathbf{x}]$, such that 
    $$T_n(\varphi(\mathbf{x})) = H(S_{n,1}(\varphi(\mathbf{x})), \cdots, S_{n,n}(\varphi(\mathbf{x}))).$$
Hence $F(\mathbf{x}) \in \mathbb{Q}[S_{n,1}(\varphi(\mathbf{x})), \cdots, S_{n,n}(\varphi(\mathbf{x}))]$.
\end{proof}

Taking a more general case than Proposition \ref{prop len0 B=0}, we have the following example.

\begin{eg} \label{eg len0 B=0}
    Given a seed $\Omega := (B, \mathbf{x}, R, \mathbf{Z})$, where $B = (0)_{n \times n}$. Obviously,  $\mathcal{S}(\Omega) \supset \{ \mu_i \mid  i \in [1,n]\}$.
    By Definition \ref{def seed}, we know that $Z_i(1)$ is a positive integer greater than or equal to $2$ and the exchange polynomial tuple $\mathbf{P}(\Omega) = (Z_1(1), \dots,  Z_n(1))$. Take $H(\mathbf{x}) \in \mathbb{Q}[\mathbf{x}]$ and let $$F(\mathbf{x}) := H\bigg(\frac{x_1^2 + Z_1(1)}{x_1}, \cdots, \frac{x_n^2 + Z_n(1)}{x_n}\bigg).$$
    It is easy to check $F(\mu_i(\mathbf{x})) = F(\mathbf{x})$ for $i \in [1,n]$. So $F(\mathbf{x}) \in \mathbb{Q}[\mathbf{x}^{\pm}]^{\langle\mu_i \mid i \in [1,n]\rangle}$.
\end{eg}

\subsection{Existence of some cluster symmetric polynomials as examples} 

In this subsection, we consider the existence of nonconstant cluster symmetric polynomials related to some generalized cluster algebras and answer two questions posed by Gyoda and Matsushita in \cite{Gyoda}. We first recall their work. In \cite{Gyoda}, they show Table \ref{tab omega and f}. Observing the table, it is easy to check that for any $i \in [1,6]$, the Laurent polynomial $F_{3,i}(\mathbf{x}) \in \mathbb{Q}[\mathbf{x}^{\pm}]^{\langle\mu_1, \mu_2, \mu_3\rangle}$, where $\mu_1, \mu_2, \mu_3 \in \mathcal{S}(\Omega_{3,i})$ and $\Omega_{3,i} := (B_{3,i}, \mathbf{x}, R_{3,i}, \mathbf{Z}_{3,i})$.

\begin{table}[htbp]
\caption{Seed $\Omega_{3,i} := (B_{3,i}, \mathbf{x}, R_{3,i},  \mathbf{Z}_{3,i})$ and corresponding cluster symmetric polynomial $F_{3,i}$.} \label{tab omega and f}
\begin{center}
\renewcommand{\arraystretch}{1.2}
\resizebox{\textwidth}{!}{
$\begin{array}{|c|ccl|c|}
\hline
i & B_{3,i} & R_{3,i} & \qquad \quad  \mathbf{Z}_{3,i} & \begin{array}{c}
~\\[-3mm] 
F_{3,i}(x,y,z) \\[-3mm]
~ 
\end{array}\\
\hline
1 & \left[\begin{array}{rrr}
0&2&-2\\-2&0&2\\2&-2&0\end{array}\right]
& \begin{array}{c}
    r_1 = 1\\
    r_2 = 1\\
    r_3 = 1
\end{array}
& \begin{array}{l}
~\\[-3mm]
Z_1: 1+u \\
Z_2: 1+u \\
Z_3: 1+u \\[-3mm]
~
\end{array}
& \begin{array}{c}
     \underline{x^2+y^2+z^2} \\
     {xyz} 
\end{array} \\[1mm]
\hline
2&\left[\begin{array}{rrr}0&2&-1\\-2&0&1\\2&-2&0\end{array}\right]
& \begin{array}{c}
    r_1 = 1\\
    r_2 = 1\\
    r_3 = 2
\end{array}
& \begin{array}{l}
\\[-3mm]
Z_1: 1+u \\
Z_2: 1+u \\
Z_3: 1+k_3 u+u^2\\[-3mm]
~
\end{array} 
& \begin{array}{c}
     \underline{x^2+y^2+z^2+k_3xy} \\
     {xyz} 
\end{array} \\[1mm]
\hline
3&\left[\begin{array}{rrr}0&2&-1\\-1&0&1\\1&-2&0\end{array}\right]
& \begin{array}{c}
    r_1 = 2\\
    r_2 = 1\\
    r_3 = 2
\end{array}
& \begin{array}{l}
\\[-3mm]
Z_1: 1+k_1 u+u^2 \\
Z_2: 1+u \\
Z_3: 1+k_3 u+u^2 \\[-3mm]
~
\end{array}
& \begin{array}{c}
     \underline{x^2+y^2+z^2+k_1yz+k_3xy} \\
     {xyz} 
\end{array} \\[1mm]
\hline
4& \left[\begin{array}{rrr}0&1&-1\\-1&0&1\\1&-1&0\end{array}\right]
& \begin{array}{c}
    r_1 = 2\\
    r_2 = 2\\
    r_3 = 2
\end{array}
& \begin{array}{l}
\\[-3mm]
Z_1: 1+k_1 u+u^2 \\
Z_2: 1+k_2 u+u^2 \\
Z_3: 1+k_3 u+u^2 \\[-3mm]
~
\end{array} 
& \begin{array}{c}
     \underline{x^2+y^2+z^2+k_1yz+k_2zx+k_3xy} \\
     {xyz} 
\end{array} \\[1mm]
\hline
5 &\left[\begin{array}{rrr}0&1&-1\\-4&0&2\\4&-2&0\end{array}\right]
& \begin{array}{c}
    r_1 = 1\\
    r_2 = 1\\
    r_3 = 1
\end{array}
& \begin{array}{l}
\\[-3mm]
Z_1: 1+u \\
Z_2: 1+u \\
Z_3: 1+u \\[-3mm]
~
\end{array}
& \begin{array}{c}
     \underline{x^2+y^4+z^4+2xy^2+2xz^2}  \\
     {xy^2z^2} 
\end{array} \\[1mm]
\hline
6& \left[\begin{array}{rrr}0&1&-1\\-2&0&2\\2&-2&0\end{array}\right]
& \begin{array}{c}
    r_1 = 2\\
    r_2 = 1\\
    r_3 = 1
\end{array}
& \begin{array}{l}
\\[-3mm]
Z_1: 1+ku+u^2 \\
Z_2: 1+u \\
Z_3: 1+u \\[-3mm]
~
\end{array}
&  \begin{array}{c}
     \underline{x^2+y^4+z^4+2xy^2+ky^2z^2+2xz^2} \\
     {xy^2z^2} 
\end{array} \\
\hline
\end{array}
$}
\end{center}
\end{table}

Notice that in Table \ref{tab omega and f}, the matrix $B_{3,i}R_{3,i}$ is  
$$\text{ either } \left[\begin{array}{rrr}0 & 2 & -2 \\ -2 & 0 & 2 \\ 2 & -2 & 0\end{array}\right]  \text{ or } \left[\begin{array}{rrr}0 & 1 & -1 \\ -4 & 0 & 2 \\ 4 & -2 & 0\end{array}\right].$$
There is one more seed of rank 3 that would satisfy this condition, but Gyoda and Matsushita did not find the corresponding Diophantine equation, so they asked the following question.

\begin{que}[({\cite[Question 19]{Gyoda}})] \label{que: gyoda1}
    Given a seed $\Omega_{3,7} := (B,\mathbf{x}, R,\mathbf{Z})$ where
        $$ B =\left[\begin{array}{rrr}
        0 & 1 & -1 \\
        -1 & 0 & 2 \\
        1 & -2 & 0
        \end{array}\right],
        R =\left[\begin{array}{ccc}
        4 &  &  \\
         & 1 &  \\
         &  & 1
        \end{array}\right],\\
        \begin{array}{l}
        Z_1(u)=1+k_1u+k_2u^2+k_1 u^3+u^4, \\
        Z_2(u)=1+u, \\
        Z_3(u)=1+u.\end{array}$$
    Is there a Diophantine equation corresponding to the seed $\Omega_{3,7}$?
\end{que}

\begin{rmk} \label{rmk diop to seed}
    (i) In their paper\cite{Gyoda}, Gyoda and Matsushita did not give a strict definition of ``a Diophantine equation corresponding to the seed'' and also in their paper ``the Diophantine equations corresponding to the seed'' are all non-constant Laurent polynomials with initial vector $(1,1,1)$. 
    We understand that the above question is actually a search for non-constant Laurent polynomials $F(\mathbf{x}) \in \mathbb{Q}[\mathbf{x}^{\pm}]^{\langle \mu_1, \mu_2, \mu_3 \rangle}$, that is, whether $\mathbb{Q}[\mathbf{x}^{\pm}]^{\langle \mu_1, \mu_2, \mu_3 \rangle} \ne \mathbb{Q}$ holds.

    (ii) Table \ref{tab omega and f} is quoted from Table 1 in \cite{Gyoda}. There is a slight difference in that one of the columns in our table is about Laurent polynomials, while one of the columns there is about Diophantine equations. Both are the same when considering the positive integer of the equation $F_{3, i}(x,y,z) = F_{3, i}(1,1,1)$. For example, equation $F_{3,1}(x,y,z) = F_{3,1}(1,1,1)$ is the Markov equation $x^2 + y^2 + z^2 = 3xyz$ with $xyz \ne 0$. The Laurent polynomial $F_{3,5}$ was found by Lampe in \cite{Lampe} and the Laurent polynomials $F_{3,2}, F_{3,3}, F_{3,4}, F_{3,6}$ was found by Gyoda and Matsushita in \cite{Gyoda}.
\end{rmk}

We give an affirmative answer to this question.
\begin{prop}[(Answer to  Question \ref{que: gyoda1})] \label{poly gyoda}
    Given a seed $\Omega_{3,7}$ defined in Question \ref{que: gyoda1}. The Laurent polynomial
    $$F_{3,7}(\mathbf{x}) := a\frac{x_1^2 + x_2^4 + x_3^4 + 2x_1(x_2^2 + x_3^2)  + k_1x_2x_3(x_1+ x_2^2+ x_3^2) +  k_2x_2^2x_3^2}{x_1x_2^2x_3^2}+b,$$
     belongs to the invariant Laurent polynomial ring $\mathbb{Q}[\mathbf{x}^{\pm}]^{\langle\mu_1, \mu_2, \mu_3 \rangle}$, where $\mu_i \in \mathcal{S}(\Omega_{3,7})$ and $a,b\in \mathbb{Q}$ with $a \ne 0$.
\end{prop}
\begin{proof} We know that $\mu_1(\mathbf{x}) = (\frac{x_2^4+k_1x_2^3x_3+k_2x_2^2x_3^2 +k_1x_2x_3^3+x_3^4}{x_1}, x_2, x_3)$, $\mu_2(\mathbf{x}) = (x_1, \frac{x_1+x_3^2}{x_2}, x_3)$ and $\mu_3(\mathbf{x}) = (x_1, x_2, \frac{x_1+x_2^2}{x_3})$. It is easy to check that $F_{3,7}(\mu_i(\mathbf{x})) = F_{3,7}(\mathbf{x})$ for all $i =1,2,3$.

Although we have completed the proof, we show how we constructed $F_{3,7}(\mathbf{x})$. By Proposition \ref{is CSM}, we know that $\mu_i = \psi_{id, i, \omega_i}$ for all $i \in [1,3]$, where $\psi_{id, i, \omega_i}$'s are defined in Example \ref{eg: Gyoda} while assuming $k_0 = k_4 = 1$ and $k_1 = k_3$. Then by the result in Example \ref{eg: Gyoda}(iii), we know that the following Laurent polynomial $F(\mathbf{x})$ is invariant under $\mu_1$.
    $$F(\mathbf{x}) := a\frac{x_1^2 + x_2^4+k_{1}x_2^3x_3+k_{2}x_2^2x_3^2+k_{1}x_2x_3^3+x_3^4}{x_1} + H(\mathbf{x}),$$
where $a \in \mathbb{Q}_{\ne 0}$ and $H(\mathbf{x}) \in \mathbb{Q}[\mathbf{x}]$ with $\deg^1H(\mathbf{x}) = 0$, $\deg^2 H(\mathbf{x}) \leq 4$ and $\deg^3H(\mathbf{x}) \leq 4$. Let $T(\mathbf{x}) := a(x_1^2 + x_2^4+k_{1}x_2^3x_3+k_{2}x_2^2x_3^2+k_{1}x_2x_3^3+x_3^4) + x_1H(\mathbf{x})$. Since $\deg^2T(\mathbf{x}) = \deg^3T(\mathbf{x}) = 4$, we denote $\tilde{F}(\mathbf{x}) = \frac{F(\mathbf{x})}{x_2^{2}x_3^{2}}$. Clearly, $\tilde{F}(\mathbf{x}) \in \mathbb{Q}[\mathbf{x}^{\pm}]^{\langle\mu_1\rangle}$, since Proposition \ref{ddd}.

We consider when the Laurent polynomial $\tilde{F}(\mathbf{x})$ is invariant under the action $\mu_2$. Suppose $H(\mathbf{x}) = \sum_{i=0}^{4} h_{2,i}(\mathbf{x}) x_2^i$, where $h_{2,i}(\mathbf{x})$ is a polynomial with $\deg^2h_{2,i}(\mathbf{x}) = 0$.  Sorting the polynomial $T(\mathbf{x})$ by powers of $x_2$, we have
\begin{eqnarray*}
    T(\mathbf{x}) &= &(a + x_1h_{2,4}(\mathbf{x}))x_2^4 + (ak_1x_3 + x_1h_{2,3}(\mathbf{x}))x_2^3 + (ak_1x_3^2 + x_1h_{2,2}(\mathbf{x}))x_2^2\\
    & & +~(ak_1x_3^3 + x_1h_{2,1}(\mathbf{x}))x_2 + (ax_1^2+ax_3^4 + x_1h_{2,0}(\mathbf{x})).
\end{eqnarray*}
If $\tilde{F}(\mathbf{x}) \in \mathbb{Q}[\mathbf{x}^{\pm}]^{\langle \mu_2 \rangle}$, then by Theorem \ref{thm len1_T2f} and equation \eqref{len1_f=fp@s}, the following equations 
$$\begin{cases}
    ax_1^2+ax_3^4 + x_1h_{2,0}(\mathbf{x}) = (a + x_1h_{2,4}(\mathbf{x}))(x_1 + x_3^2)^2,\\
    ak_1x_3^3 + x_1h_{2,1}(\mathbf{x}) = (ak_1x_3 + x_1h_{2,3}(\mathbf{x}))(x_1 + x_3^2),
\end{cases}$$
must hold. By solving the above equations, we know that $h_{2,4}(\mathbf{x}) = h_{2,3}(\mathbf{x}) = 0$, $h_{2,1}(\mathbf{x}) = ak_1x_3$ and $h_{2,0}(\mathbf{x}) = 2ax_3^2$. So $H(\mathbf{x}) = (\sum_{i=0}^4 b_ix_3^i)x_2^2 + 2ak_1x_2x_3 + 2ax_3^2$ where $b_i \in \mathbb{Q}$. 

Then, similarly, we consider the case of the action $\mu_3$. We have $b_4 = b_3 = b_1 = 0$ and $b_0 = 2a$. That is, $H(\mathbf{x}) = (2a + b_2x_3^2)x_2^2 + 2ak_1x_2x_3 + 2ax_3^2$. Hence,
$$\tilde{F}(\mathbf{x}) = a\frac{x_1^2 + x_2^4 + x_3^4 + 2x_1(x_2^2 + x_3^2)  + k_1x_2x_3(x_1 + x_2^2 + x_3^2)+  k_2x_2^2x_3^2}{x_1x_2^2x_3^2} + b_2.$$

Therefore, the Laurent polynomial $\tilde{F}(\mathbf{x}) \in \mathbb{Q}[\mathbf{x}^{\pm}]^{\langle \mu_1,\mu_2,\mu_3\rangle}$.
\end{proof}

Note that the above cluster symmetric polynomial $F_{3,7}$ can be constructed using the MATLAB program in Appendix \ref{code ca2diop}, as shown in Code \ref{code: gyoda}.

Further, Gyoda and Matsushita ask the following question.

\begin{que}[({\cite[Question 20]{Gyoda}})]
\label{que gyoda2}
(1) Given a seed $\Omega := (B, \mathbf{x},  R, \mathbf{Z})$ of rank $n=3$ that satisfies the following two conditions
    \begin{eqnarray}
        \mu_i(B) = -B &\text{ for all } i \in [1,n], \label{cond muB = -B}\\
        Z_i(u) = u^{r_i}Z_i(u^{-1}) &\text{ for all } i \in [1,n]. \label{cond recip}
    \end{eqnarray}
Whether there exists a seed of $BR$ which is 
$$\text{ neither } \left[\begin{array}{rrr}0 & 2 & -2 \\ -2 & 0 & 2 \\ 2 & -2 & 0\end{array}\right]  \text{ nor } \left[\begin{array}{rrr}0 & 1 & -1 \\ -4 & 0 & 2 \\ 4 & -2 & 0\end{array}\right],$$
such that there exists a Diophantine equation corresponding to the seed $\Omega$?

(2) Is there a general way to construct a Diophantine equation from the information of the seed $\Omega := (B,\mathbf{x}, R, \mathbf{Z})$?
\end{que} 

We first consider question (1). As in Remark \ref{rmk diop to seed}(i), we consider whether the relation $\mathbb{Q}[\mathbf{x}^{\pm}]^{\langle \mu_1, \mu_2, \mu_3 \rangle} \ne \mathbb{Q}$ holds. The following examples give an affirmative answer.

\begin{eg}\label{eg existence of F}
(i) Given a seed $\Omega_{3,0} := (B,\mathbf{x}, R,\mathbf{Z})$ where
        $$ B =\left[\begin{array}{rrr}
        0 &  &  \\
         & 0 &  \\
         &  & 0
        \end{array}\right],
        R =\left[\begin{array}{ccc}
        r &  &  \\
         & r &  \\
         &  & r
        \end{array}\right],\\
        \begin{array}{l}
        Z_1(u)= 1 + z_1u + \cdots + z_1u^{r-1} + u^r, \\
        Z_2(u)=Z_1(u), \\
        Z_3(u)=Z_1(u).\end{array}$$
Let $c = Z_1(1)$. For any symmetric polynomial $\phi(\mathbf{x}) \in \mathbb{Q}[\mathbf{x}]^{\mathfrak{S}_3}$, by Proposition \ref{prop len0 B=0}, we know that the Laurent polynomial
    $$F_{3,0}(x,y,z) := \phi\bigg(\frac{x^2+c}{x}, \frac{y^2+c}{y},\frac{z^2+c}{z}\bigg)$$
is invariant under $\mu_1, \mu_2, \mu_3 \in \mathcal{S}(\Omega_{3,0})$. Hence $\mathbb{Q}[\mathbf{x}^{\pm}]^{\langle \mu_1, \mu_2, \mu_3 \rangle} \ne \mathbb{Q}$.

    (ii) Given a seed $\Omega_{3,8} := (B,\mathbf{x}, R,\mathbf{Z})$ where
        $$ B =\left[\begin{array}{rrr}
        0 & 4 & -4 \\
        -1 & 0 & 2 \\
        1 & -2 & 0
        \end{array}\right],
        R =\left[\begin{array}{ccc}
        1 &  &  \\
         & 1 &  \\
         &  & 1
        \end{array}\right],\\
        \begin{array}{l}
        Z_1(u)=1+u, \\
        Z_2(u)=1+u, \\
        Z_3(u)=1+u.\end{array}$$
It is easy to check that conditions \eqref{cond muB = -B} and \eqref{cond recip} are satisfied, and the Laurent polynomial
    $$F_{3,8}(x,y,z) := \frac{x^4 + y^2 + z^2 + 2yz}{x^2yz}$$
is invariant under $\mu_1, \mu_2, \mu_3 \in \mathcal{S}(\Omega_{3,8})$. Hence $\mathbb{Q}[\mathbf{x}^{\pm}]^{\langle \mu_1, \mu_2, \mu_3 \rangle} \ne \mathbb{Q}$. In fact, the Laurent polynomial $F_{3,8}$ was first constructed by Kaufman in \cite{Kaufman}.

    (iii) Given a seed $\Omega_{3,9} := (B,\mathbf{x}, R,\mathbf{Z})$ where
        $$ B =\left[\begin{array}{rrr}
        0 & 2 & -4 \\
        -1 & 0 & 2 \\
        1 & -1 & 0
        \end{array}\right],
        R =\left[\begin{array}{ccc}
        1 &  &  \\
         & 2 &  \\
         &  & 1
        \end{array}\right],\\
        \begin{array}{l}
        Z_1(u)=1+u, \\
        Z_2(u)=1+k_2+u^2, \\
        Z_3(u)=1+u.\end{array}$$
It is easy to check that conditions \eqref{cond muB = -B} and \eqref{cond recip} are satisfied, and the Laurent polynomial
    $$F_{3,9}(x,y,z) := \frac{x^4 + k_2x^2z + y^2 + z^2 + 2yz}{x^2yz}$$
is invariant under $\mu_1, \mu_2, \mu_3 \in \mathcal{S}(\Omega_{3,9})$. Hence $\mathbb{Q}[\mathbf{x}^{\pm}]^{\langle \mu_1, \mu_2, \mu_3 \rangle} \ne \mathbb{Q}$.

    (iv) Given a seed $\Omega_{3,10} := (B,\mathbf{x}, R,\mathbf{Z})$ where
        $$ B =\left[\begin{array}{rrr}
        0 & 2 & -2 \\
        -1 & 0 & 1 \\
        1 & -1 & 0
        \end{array}\right],
        R =\left[\begin{array}{ccc}
        1 &  &  \\
         & 2 &  \\
         &  & 2
        \end{array}\right],\\
        \begin{array}{l}
        Z_1(u)=1+u, \\
        Z_2(u)=1+k_2+u^2, \\
        Z_3(u)=1+k_3+u^2.\end{array}$$
It is easy to check that conditions \eqref{cond muB = -B} and \eqref{cond recip} are satisfied, and the Laurent polynomial
    $$F_{3,10}(x,y,z) = \frac{x^4 + k_3x^2y+k_2x^2z+y^2 + z^2 + 2yz}{x^2yz}$$
is invariant under $\mu_1, \mu_2, \mu_3 \in \mathcal{S}(\Omega_{3,10})$. Hence $\mathbb{Q}[\mathbf{x}^{\pm}]^{\langle \mu_1, \mu_2, \mu_3 \rangle} \ne \mathbb{Q}$.
\end{eg}

Are there any other seeds than the four mentioned above? To do so, we first prove the following theorem.

\begin{thm} \label{thm 4bounded}
    Given a seed $\Omega := (B, \mathbf{x}, R, \mathbf{Z})$.
    
    (i)  Suppose $\sigma\mu_i, \tau\mu_j \in \mathcal{S}(\Omega)$ with $i \ne j$ and  there exists a $\frac{\bm{
    \eta}}{\mathbf{d}}$ type Laurent polynomial $F(\mathbf{x}) \in \mathbb{Q}[\mathbf{x}^{\pm}]^{\langle \sigma\mu_i, \tau\mu_j\rangle}$. If $\eta_i \ne 0$, then 
    \begin{equation} \label{4bounded}
        4 \geq r_ir_j\max\{|b_{ji}|,|b_{\sigma(j)i}|\}\max\{|b_{ij}|,|b_{\tau(i)j}|,|b_{\sigma^{-1}(i)j}|,|b_{\tau(\sigma^{-1}(i))j}|\}.
    \end{equation}
    
    (ii) Suppose $\{\mu_i \mid i \in [1,n]\} \subset \mathcal{S}(\Omega)$. We denote a set
        \begin{equation*}
            \mathcal{I} := \{ i, j \in [1,n] \mid  i \ne j \text{ and } r_ir_j|b_{ij}b_{ji}| > 4 \}. 
        \end{equation*}
        
    If $\# \mathcal{I} = n$, then $\mathbb{Q}[\mathbf{x}^{\pm}]^{\langle\mu_i \mid i \in [1,n]\rangle} = \mathbb{Q}$.

    If $\# \mathcal{I} = n-1$. Let $s \in \{1, \cdots, n\} \setminus \mathcal{I}$. Then 
    $$\mathbb{Q}[\mathbf{x}^{\pm}]^{\langle\mu_i \mid i \in [1,n]\rangle} =  \begin{cases}
        \mathbb{Q}\bigg[\frac{Z_s(1) + x_s^2}{x_s}\bigg], &\text{if } b_{ks} = 0 \text{ for all } k \in [1,n],\\
        \mathbb{Q}, &\text{otherwise}.
    \end{cases}$$
\end{thm}
\begin{proof}
    (i) By Proposition \ref{is CSM}, we know $\sigma\mu_i = \psi_{\sigma, i, \pi_i(\Omega)}$ and $\tau\mu_j = \psi_{\tau, j, \pi_j(\Omega)}$, where $\pi_i(\Omega) = ((b_{1i}, \cdots, b_{ni}), r_i, Z_i)$ and $\pi_j(\Omega) = ((b_{1j}, \cdots, b_{nj}), r_j, Z_j)$. By Corollary \ref{cor 4 eta}, Relation \eqref{4bounded} holds.
    
    (ii) Let $F(\mathbf{x}) \in \mathbb{Q}[\mathbf{x}^{\pm}]^{\langle\mu_i \mid i \in [1,n]\rangle}$. Suppose $F(\mathbf{x})$ is of type $\frac{\bm{\eta}}{\mathbf{d}}$. By (i), we know $\eta_i = 0$ for all $i \in \mathcal{I}$. By Equation \eqref{len1_eta_s=2d_s}, we know $d_i = 0$ for all $i \in \mathcal{I}$. When $\# \mathcal{I} = n$, then $F(\mathbf{x})$ is of type $\frac{\mathbf{0}}{\mathbf{0}}$, that is, $F(\mathbf{x}) \in \mathbb{Q}$.

    When $\# \mathcal{I} = n-1$. If $b_{ks} = 0$ for all $k \in [1,n]$. Then, by Proposition \ref{is CSM} and Theorem \ref{F in ring: id}, we have $ \mathbb{Q}[\mathbf{x}^{\pm}]^{\langle\mu_s\rangle} = \mathbb{Q}[\frac{P_{\pi_s(\Omega)}(\mathbf{x}) + x_s^2}{x_s}, \mathbf{x}^{\pm}\setminus x_{s}^{\pm}] = \mathbb{Q}[\frac{Z_s(1) + x_s^2}{x_s}, \mathbf{x}^{\pm}\setminus x_{s}^{\pm}]$. 
    Since $\eta_i = 0$ for all $i \in \mathcal{I}$, we know that $\mathbb{Q}[\mathbf{x}^{\pm}]^{\mu_s} = \mathbb{Q}[\frac{Z_s(1) + x_s^2}{x_s}]$ and $\mathbb{Q}[\mathbf{x}^{\pm}]^{\mu_i} = \mathbb{Q}[x_s^{\pm}]$ for $i \ne s$.
        Hence, 
        \begin{eqnarray*}
            \mathbb{Q}[\mathbf{x}^{\pm}]^{\langle\mu_i \mid i \in [1,n]\rangle} = \mathbb{Q}[\mathbf{x}^{\pm}]^{\mu_s} \cap \bigcap_{i \in \mathcal{I}} \mathbb{Q}[\mathbf{x}^{\pm}]^{\langle\mu_i\rangle} = \mathbb{Q}\bigg[\frac{Z_s(1) + x_s^2}{x_s}\bigg].
        \end{eqnarray*}
    If there exists $k \in [1,n]$ such that $b_{ks} \ne 0$. Clearly, $k \ne s$. By Equation \eqref{eta_k>=drb}, we know that $0 = \eta_k \geq \frac{1}{2}\eta_sr_s|b_{ks}|$. Then $d_s = \eta_s = 0$. Hence $F(\mathbf{x})$ is of type $\frac{\mathbf{0}}{\mathbf{0}}$, that is, $F(\mathbf{x}) \in \mathbb{Q}$.
\end{proof}

\begin{table}[htbp]
\caption{Seed $\Omega_{2,i} := (B_{2,i}, \mathbf{x}, R_{2,i},  \mathbf{Z}_{2,i})$ and its cluster symmetric polynomial $F_{2,i}$.} \label{tab rank 2 1-cspoly}
\begin{center}
\renewcommand{\arraystretch}{1} \resizebox{\textwidth}{!}{
$\begin{array}{|c|ccl|c|}
\hline
i & B_{2,i} & R_{2,i} & \qquad \quad  \mathbf{Z}_{2,i} & \begin{array}{c}
~\\[-2mm] 
F_{2,i}(x,y) \\[-2mm]
~ 
\end{array}\\
\hline
1 & \left[\begin{array}{rrr}
0&2\\-2&0\end{array}\right]
& \begin{array}{c}
    r_1 = 1\\
    r_2 = 1
\end{array}
& \begin{array}{l}
~\\[-2mm]
Z_1: 1+u \\
Z_2: 1+u \\[-2mm]
~
\end{array}
& \begin{array}{c}
     \underline{x^2+y^2+1} \\
     {xy} 
\end{array} \\[3mm]
\hline
2&\left[\begin{array}{rrr}0&2\\-1&0\end{array}\right]
& \begin{array}{c}
    r_1 = 2\\
    r_2 = 1
\end{array}
& \begin{array}{l}
\\[-2mm]
Z_1: 1+k_1 u+u^2 \\
Z_2: 1+u \\ [-2mm]
~
\end{array}
& \begin{array}{c}
     \underline{x^2+y^2+k_1y + 1} \\
     {xy} 
\end{array} \\[1mm]
\hline
3& \left[\begin{array}{rrr}0&1\\-1&0\end{array}\right]
& \begin{array}{c}
    r_1 = 2\\
    r_2 = 2
\end{array}
& \begin{array}{l}
\\[-2mm]
Z_1: 1+k_1 u+u^2 \\
Z_2: 1+k_2 u+u^2 \\[-2mm]
~
\end{array} 
& \begin{array}{c}
     \underline{x^2+y^2+k_1y+k_2x + 1} \\
     {xy} 
\end{array} \\[1mm]
\hline
4 &\left[\begin{array}{rrr}0&1\\-4&0\end{array}\right]
& \begin{array}{c}
    r_1 = 1\\
    r_2 = 1\\
\end{array}
& \begin{array}{l}
\\[-2mm]
Z_1: 1+u \\
Z_2: 1+u \\[-2mm]
~
\end{array}
& \begin{array}{c}
     \underline{x^2+y^4+2x + 1}  \\
     {xy^2} 
\end{array} \\[1mm]
\hline
5& \left[\begin{array}{rrr}0&1\\-2&0\end{array}\right]
& \begin{array}{c}
    r_1 = 2\\
    r_2 = 1
\end{array}
& \begin{array}{l}
\\[-2mm]
Z_1: 1+ku+u^2 \\
Z_2: 1+u \\[-2mm]
~
\end{array}
&  \begin{array}{c}         \underline{x^2+y^4+ky^2+2x + 1} \\
     {xy^2} 
\end{array} \\
\hline
6& \left[\begin{array}{rrr}0&1\\-1&0\end{array}\right]
& \begin{array}{c}
    r_1 = 4\\
    r_2 = 1
\end{array}
& \begin{array}{l}
\\[-2mm]
Z_1: \sum_{i=0}^4k_iu^i\\
Z_2: 1+u \\[-2mm]
~
\end{array}
&  \begin{array}{c}         
\underline{x^2 + 2x  + k_1xy + Z_1(y) } \\
     {xy^2}
\end{array} \\
\hline
7& \left[\begin{array}{rrr}0&1\\-3&0\end{array}\right]
& \begin{array}{c}
    r_1 = 1\\
    r_2 = 1
\end{array}
& \begin{array}{l}
\\[-2mm]
Z_1: 1+u\\
Z_2: 1+u \\[-2mm]
~
\end{array}
&  \begin{array}{c}         
\underline{{(x^2+2x+Z_1(y))(y+1) + xy^3}}\\ 
     {xy^2}
\end{array} \\
\hline
8 & \left[\begin{array}{rrr}0&1\\-1&0\end{array}\right]
& \begin{array}{c}
    r_1 = 3\\
    r_2 = 1
\end{array}
& \begin{array}{l}
\\[-2mm]
Z_1: \sum_{i=0}^3k_iu^i\\
Z_2: 1+u \\[-2mm]
~
\end{array}
& \begin{array}{c}         
\underline{(x^2+2x+Z_1(y))(y+1) + xy(y^2+k_1)}\\
{xy^2}
\end{array} \\
\hline
9 & \left[\begin{array}{rrr}0&1\\-2&0\end{array}\right]
& \begin{array}{c}
    r_1 = 1\\
    r_2 = 1
\end{array}
& \begin{array}{l}
\\[-2mm]
Z_1: 1+u\\
Z_2: 1+u \\[-2mm]
~
\end{array}
&  \begin{array}{c}         
\underline{xy^2 + y^2 + x^2 + 2x + 1} \\
     {xy}
\end{array} \\
\hline
10 & \left[\begin{array}{rrr}0&1\\-1&0\end{array}\right]
& \begin{array}{c}
    r_1 = 2\\
    r_2 = 1
\end{array}
& \begin{array}{l}
\\[-2mm]
Z_1: 1+k_1u+u^2\\
Z_2: 1+u \\[-2mm]
~
\end{array}
&  \begin{array}{c}         
\underline{xy^2 + y^2 + k_1y + x^2 + 2x + 1} \\
     {xy}
\end{array} \\
\hline
11 & \left[\begin{array}{rr}0&1\\-1&0\end{array}\right]
& \begin{array}{c}
    r_1 = 1\\
    r_2 = 1
\end{array}
& \begin{array}{l}
\\[-2mm]
Z_1: 1+u\\
Z_2: 1+u \\[-2mm]
~
\end{array}
&  \begin{array}{c}         
\underline{x^2+y^2+2x+2y+x^2y+xy^2+1} \\
     {xy}
\end{array} \\
\hline
12 & \left[\begin{array}{rrr}0&0\\0&0\end{array}\right]
& \begin{array}{c}
    r_1 \geq 1\\
    r_2 \geq 1
\end{array}
& \begin{array}{l}
\\[-2mm]
Z_1: \sum_{i=0}^rk_{1,i}u^i\\
Z_2: \sum_{i=0}^rk_{2,i}u^i \\[-2mm]
~
\end{array}
&  \begin{array}{c}         
\underline{(x^2+Z_1(1))(y^2+Z_2(1)) } \\
     {xy}
\end{array} \\
\hline
\end{array}$}
\end{center}
\end{table}

We first use this theorem in the case of rank $n = 2$ to describe the equivalence condition that a generalized cluster algebra of rank $2$ has a non-constant cluster symmetric polynomial.

\begin{prop}\label{rank2 nonconstant}
    For any seed $\Omega := (B, \mathbf{x}, R, \mathbf{Z})$ of rank $n=2$. The relation $\mathbb{Q}[\mathbf{x}^{\pm}]^{\langle\mu_1, \mu_2\rangle} \ne \mathbb{Q}$ holds, if and only if, the seed $\Omega$ is permutation equivalent to one of the seeds $\Omega_{2,i}:=(\pm B_{2,i}, \mathbf{x}, R_{2,i}, \mathbf{Z}_{2,i})$ listed in Table \ref{tab rank 2 1-cspoly}, that is, $\Omega = \sigma(\Omega_{2,i})$ for some $i \in [1,12]$ and $\sigma \in \mathfrak{S}_2$.
\end{prop}
\begin{proof}
    It is easy to check $\mu_1, \mu_2 \in \mathcal{S}(\Omega)$. We denote a set
    $$\mathcal{I} := \{ i, j \in [1,2] \mid  i \ne j \text{ and } r_ir_j|b_{ij}b_{ji}| > 4 \}.$$
    Then $\# \mathcal{I}$ must be $2$ or $0$. When $\# \mathcal{I} = 2$, by Theorem \ref{thm 4bounded}(ii), we know $\mathbb{Q}[\mathbf{x}^{\pm}]^{\langle\mu_1, \mu_2\rangle} = \mathbb{Q}$. When $\# \mathcal{I} = 0$, we have $r_1r_2|b_{12}b_{21}| \leq 4$. Then the seed $\Omega$ is permutation equivalent to one of the seeds $\Omega_{2,i}:=(\pm B_{2,i}, \mathbf{x}, R_{2,i}, \mathbf{Z}_{2,i})$ for some $i \in [1,12]$.
    
    We claim that for any $i\in [1,12]$, the non-constant cluster symmetric polynomial $F_{2,i}$ belongs to the invariant ring $\mathbb{Q}[\mathbf{x}^{\pm}]^{\langle\mu_1, \mu_2\rangle}$, where $\mu_1, \mu_2 \in \mathcal{S}(\Omega_{2,i})$. 
    
    When $i = 1, \dots, 5$ or $11$, Gyoda and Matsushita show it in \cite[Table 3]{Gyoda}.
    
    When $i = 7$ or $9$, Chen and Li prove it in \cite[Example 2.20,2.21]{ChenLi}. 
    
    When $i = 6,8,10,12$, the non-constant cluster symmetric polynomial $F_{2,i}$ is constructed by our method (Theorem \ref{thm len1}). We prove that $F_{2,8} \in \mathbb{Q}[\mathbf{x}^{\pm}]^{\langle\mu_1, \mu_2\rangle}$, other is similar. Since
    \begin{eqnarray*}
        F_{2,8}(x,y) &=& \frac{x^2(1+y) + xy(2+2k+y^2) + Z_1(y)(1+y)}{xy^2}\\
         &=& \frac{y^4 + (1+k_1+x)y^3 + 2ky^2+ (1+k_1+x)Z_2(x)y+Z^2_2(x)}{xy^2},
    \end{eqnarray*}
    it is easy to check that $F_{2,8}(\mu_1(x,y)) = F_{2,8}(Z_1(y)/x,y)) = F_{2,8}(x,y)$ and $F_{2,8}(\mu_2(x,y)) = F_{2,8}(x,Z_2(x)/y)) = F_{2,8}(x,y)$.
\end{proof}

\begin{rmk}
    (i) Table \ref{tab rank 2 1-cspoly} lists all generalized cluster algebras of rank $2$ which have a non-constant cluster symmetric polynomial. Some of these generalized cluster algebras have other non-constant cluster symmetric polynomials. For example, the generalized cluster algebra $\mathcal{A}(\Omega_{2,8})$ has the cluster symmetric polynomial
    $$H_{2,8}(x,y) := \frac{x^4  + x(x^2 + Z_1(y))(y^3+k_1y+4) + x^2(k_1y^2+4k_1y+6) + Z^2_1(y)}{x^2y^3}$$
    and the generalized cluster algebra $\mathcal{A}(\Omega_{2,10})$ has the cluster symmetric polynomial
    $$H_{2,10}(x,y) := \frac{x^2y^2 + x^2 + k_1xy + 2x + y^4 + k_1y^3 + 2y^2 + k_1y + 1}{xy^2}.$$
    
    (ii) When the mutation degree matrix $R = I_2$, that is, when considering the cluster algebras of rank $2$, Proposition \ref{rank2 nonconstant} for this special case has been proved by Chen and Li in \cite[Theorem 2.36]{ChenLi}. They proved it using some general term formulas of $\mathbf{d}$-vectors. This $\mathbf{d}$-vector is about cluster variables, while our $\mathbf{d}$ is about Laurent polynomials.
\end{rmk}

Now, we consider the rank $n=3$. Case (i) of the following corollary answers Question \ref{que gyoda2}. Note that we do not need to discuss whether Condition \eqref{cond recip} holds, since it is clearly held by the reciprocal condition \eqref{recip}.

\begin{prop}\label{rank3 nonconstant}
For any seed $\Omega := (B, \mathbf{x}, R, \mathbf{Z})$ of rank $n=3$. Suppose $\mu_1, \mu_2, \mu_3 \in \mathcal{S}(\Omega)$.

(i) If $\mu_i(B) = -B$ for all $i \in [1,3]$. Then the relation
$\mathbb{Q}[\mathbf{x}^{\pm}]^{\langle\mu_1, \mu_2, \mu_3\rangle} \ne \mathbb{Q}$ holds, if and only if, the matrix $BR$ is permutation equivalent to one of the following four matrices 
    \begin{eqnarray*}
        A_1 = \left[\begin{array}{rrr}0 &  &  \\  & 0 &  \\  &  & 0\end{array}\right],  A_2 = \left[\begin{array}{rrr}0 & 2 & -2 \\ -2 & 0 & 2 \\ 2 & -2 & 0\end{array}\right],
        A_3 = \left[\begin{array}{rrr}0 & 1 & -1 \\ -4 & 0 & 2 \\ 4 & -2 & 0\end{array}\right], A_4 = \left[\begin{array}{rrr}0 & 4 & -4 \\ -1 & 0 & 2 \\ 1 & -2 & 0\end{array}\right].
    \end{eqnarray*}
That is, $BR = \sigma(A_k)$ for some $k \in [1,4]$ and $\sigma \in \mathfrak{S}_3$.

(ii) The relation $\mathbb{Q}[\mathbf{x}^{\pm}]^{\langle\mu_1, \mu_2, \mu_3\rangle} \ne \mathbb{Q}$ holds, if and only if, the matrix $BR$ is permutation equivalent to one of the following matrices
$$A_1,~ A_2,~ A_3,~ A_4,~ A_5 = \left[\begin{array}{rrr}0 & b & 0 \\ c & 0 & 0 \\ 0 & 0 & 0\end{array}\right],$$
where $b,c \in \mathbb{Z}_{\ne 0}$. That is, $BR = \sigma(A_k)$ for some $k \in [1,5]$ and $\sigma \in \mathfrak{S}_3$.
\end{prop}
\begin{proof} 
    Denote a set $\mathcal{I}_{B,R} := \{ i, j \in [1,3] \mid  i \ne j \text{ and } r_ir_j|b_{ij}b_{ji}| > 4 \}.$    
    
    (i)  When $\#\mathcal{I}_{B,R} = 3$. By Theorem \ref{thm 4bounded}(ii), we have $\mathbb{Q}[\mathbf{x}^{\pm}]^{\langle\mu_1, \mu_2, \mu_3 \rangle} = \mathbb{Q}$.

    When $\#\mathcal{I}_{B,R} = 2$. Let $s \in \{1,2,3\}\setminus \mathcal{I}$. Then by Theorem \ref{thm 4bounded}(ii), we have
    $$\mathbb{Q}[\mathbf{x}^{\pm}]^{\langle \mu_1, \mu_2, \mu_3 \rangle} =  \begin{cases}
        \mathbb{Q}[\frac{Z_s(1) + x_s^2}{x_s}], &\text{if } b_{ks} = 0 \text{ for all } k \in [1,n],\\
        \mathbb{Q}, &\text{otherwise}.
    \end{cases}$$
    Assume $b_{ks} = 0$ for all $k \in [1,n]$. Since $B$ is skew-symmetrizable, we know that $b_{sk} = 0$ for all $k \in [1,n]$. Then $\mu_s(B) = B$. Since $\mu_s(B) = -B$, we know that $B = A_1$. But it leads to a contradiction since $\#\mathcal{I}_{A_1, R} = 0$. Hence, the assumption is false; by Theorem \ref{thm 4bounded}(ii), we have $\mathbb{Q}[\mathbf{x}^{\pm}]^{\langle\mu_1, \mu_2, \mu_3 \rangle} = \mathbb{Q}$.

    When $\#\mathcal{I}_{B,R} = 1$, it is impossible.
    
    When $\#\mathcal{I}_{B,R} = 0$, without loss of generality, we assume that $r_1r_2|b_{12}b_{21}| \leq 4$ and $r_1r_3|b_{13}b_{31}| \leq 4$. Taking into account the entries $(2,3)$ and $(3,2)$ of the matrices on both sides of the equation $\mu_1(B) = -B$, we have
    $2b_{23} = -r_1 \operatorname{sgn}(b_{21})[b_{21}b_{13}]_+$ and $2b_{32} = -r_1 \operatorname{sgn}(b_{12})[b_{31}b_{12}]_+$. Then 
    $$4|b_{23}b_{32}| = r_1^2[b_{21}b_{13}]_+[b_{31}b_{12}]_+ \leq r_1^2|b_{12}b_{21}b_{13}b_{31}| \leq \frac{16}{r_2r_3}.$$
    So $r_2r_3|b_{23}b_{32}| \leq 4$. Hence, for all $\sigma \in \mathfrak{S}_3$, we have
    \begin{equation*}
        4|b_{ik}b_{ki}| = r_k^2[b_{ij}b_{jk}]_+[b_{kj}b_{ji}]_+,    
    \end{equation*}
    where $(i,j,k) = \sigma(1,2,3)$. We consider the symbol of $b_{ij}b_{jk}$. If without loss of generality we have $b_{12}b_{23} \leq 0$, then by the above equation, we have 
    $$b_{13}b_{31} = 0 \quad \text{ and }\quad 4|b_{12}b_{21}| = r_3^2[b_{13}b_{32}]_+[b_{23}b_{31}]_+  = 0.$$ 
    So $b_{12}b_{21} = 0$. Similarly, we have $b_{23}b_{32} = 0$. Since $B$ is skew-symmetrizable, we have $B = A_1$. For this case, Example \ref{eg existence of F}(i) shows the existence of a non-constant invariant Laurent polynomial, that is, $F_{3,0}(\mathbf{x})$.

    If for all $\sigma \in \mathfrak{S}_3$, we have $b_{ij}b_{jk} > 0$, where $(i,j,k) = \sigma(1,2,3)$. Then for any mutually unequal $i,j,k \in [1,3]$, we have $-4b_{ij}b_{ji} = r_k^2b_{ik}b_{kj}b_{jk}b_{ki}$. So 
    $$ -64 = r_1^2r_2^2r_3^2b_{12}b_{21}b_{13}b_{31}b_{23}b_{32}.$$ 
    Hence, by computation, we have $BR = \sigma_k(A_k)$ for some $k \in [2,4]$ and $\sigma_k \in \mathfrak{S}_3$. Without loss of generality, let $\sigma_k = id$.
    For the case $k = 2$, the seed $\Omega$ must be one of the seeds $\Omega_{3,1}, \dots, \Omega_{3,4}$ in Table \ref{tab omega and f} and $F_{3,1}, \dots, F_{3,4}$ are the corresponding invariant Laurent polynomials. 
    For the case $k = 3$, the seed $\Omega$ must be one of the seeds $\Omega_{3,5}, \Omega_{3,6}$ in Table \ref{tab omega and f} or the seed $\Omega_{3,7}$ in Proposition \ref{poly gyoda} and $F_{3,5}, \dots, F_{3,7}$ are the corresponding invariant Laurent polynomials.
    For the case $k = 4$, the seed $\Omega$ must be one of the seeds $\Omega_{3,8},\Omega_{3,9},\Omega_{3,10}$ in Example \ref{eg existence of F}(ii)-(iv) and $F_{3,8}, F_{3,9}, F_{3,10}$ are the corresponding invariant Laurent polynomials.

    (ii) Denote a set $\mathcal{J}_{B,R} := \{ i \in [1,3] \mid \mu_i(B) \ne -B\}$.

    When $\# \mathcal{J}_{B,R} = 0$, by (i) we know that 
    $\mathbb{Q}[\mathbf{x}^{\pm}]^{\langle\mu_1, \mu_2, \mu_3\rangle} \ne \mathbb{Q}$, if and only if, the matrix $BR$ is permutation equivalent to one of the matrices $A_1, A_2, A_3, A_4$.

    When $\# \mathcal{J}_{B,R} = 1$, without loss of generality, suppose $\mathcal{J}_{B,R} = \{ 3 \}$. Since $\mu_3 \in \mathcal{G}(\Omega)$, we know $\mu_3(B) = B$. Then by Definition \ref{def mu}, we know $BR = A_5$. Since $\mu_3(B) \ne -B$, we have $b_{12}b_{21} \ne 0$. If $r_1r_2|b_{12}b_{21}| > 4$, that is, $bc>4$, then $\#\mathcal{I}_{B,R} = 2$. Hence, by Theorem \ref{thm 4bounded}(ii), we have $\mathbb{Q}[\mathbf{x}^{\pm}]^{\langle\mu_1, \mu_2, \mu_3 \rangle} = \mathbb{Q}[\frac{Z_3(1) + x_3^2}{x_3}]$.
    If $ 0 < r_1r_2|b_{12}b_{21}| \leq 4$, that is, $0<bc<4$, then we know that
    $$BR = \sigma_i\bigg( \left[\begin{array}{rrr}\pm B_{2,i}R_{2,i} & \\  & 0\end{array}\right]\bigg),$$
    for some $i \in [1,11]$ and $\sigma_i \in \mathfrak{S}_3$ with $\sigma_i(3) = 3$, where $B_{2,i},R_{2,i}$ are listed in Table \ref{tab rank 2 1-cspoly}. Without loss of generality, we assume $\sigma_i = id$. For $i \in [1,11]$, we denote a seed $\Omega'_{3,i} := (B'_{3,i}, \mathbf{x}, R'_{3,i}, \mathbf{Z}'_{3,i})$, where $B'_{3,i} := \left[\begin{array}{rrr}\pm B_{2,i} & \\  & 0\end{array}\right]$, $R'_{3,i} := \left[\begin{array}{rrr} R_{2,i} & \\  & r\end{array}\right]$, $r \geq 1$ and $\mathbf{Z}_{3,i}' := (\pi_1(\mathbf{Z}_{2,i}), \pi_2(\mathbf{Z}_{2,i}), \sum_{i=0}^{r}k'_iu^i)$.
    Then the seed $\Omega$ must be one of the seeds $\Omega'_{3,1}, \dots, \Omega'_{3,11}$. For all $i \in [1,11]$, Laurent polynomial $F'_{3,i}(x,y,z) := F_{2,i}(x,y)$ is the invariant Laurent polynomial related to the seed $\Omega'_{3,i}$.
    
    When $\# \mathcal{J}_{B,R} \geq 2$, without loss of generality, suppose $2,3 \in \mathcal{J}_{B,R}$. Then $\mu_2(B) = \mu_3(B) = B$. Hence $b_{12} = \cdots = b_{n2} = b_{13} = \cdots = b_{n3} = 0$ and $B = A_1$. But it is impossible, since $\#\mathcal{J}_{A_1,R} = 0$.
\end{proof}

The above conclusions can also be immediately classified using the irreducibility of matrices. A matrix is called \textbf{irreducible} if it is not similar to a block upper triangular matrix with at least two blocks via a permutation.

\begin{cor}
\label{cor rank3 nonconstant}
For any seed $\Omega := (B, \mathbf{x}, R, \mathbf{Z})$ of rank $3$. Suppose $\mu_1, \mu_2, \mu_3 \in \mathcal{S}(\Omega)$.

(i) If $B$ is irreducible, then the relation
$\mathbb{Q}[\mathbf{x}^{\pm}]^{\langle\mu_1, \mu_2, \mu_3\rangle} \ne \mathbb{Q}$ holds, if and only if, the matrix $BR$ is permutation equivalent to one of the matrices $A_2, A_3, A_4$ defined in Proposition \ref{rank3 nonconstant}.

(ii) If $B$ is reducible, then the relation
$\mathbb{Q}[\mathbf{x}^{\pm}]^{\langle\mu_1, \mu_2, \mu_3\rangle} \ne \mathbb{Q}$ holds, if and only if, the matrix $BR$ is permutation equivalent to one of the matrices $A_1, A_5$ defined in Proposition \ref{rank3 nonconstant}.
\end{cor}

\begin{rmk}
\label{rmk irr B}
    We list the seeds that satisfy Corollary \ref{cor rank3 nonconstant}(i) in Table \ref{tab all}. In the next section, we discuss the solutions to the corresponding Diophantine equations. In fact, Table \ref{tab all} extends from Table \ref{tab omega and f} with $\Omega_{3,7}$ in Proposition \ref{poly gyoda} and $\Omega_{3,8},\dots,\Omega_{3,10}$ in Example \ref{eg existence of F}.
\end{rmk}

Finally, we discuss Question \ref{que gyoda2}(2). In fact, Proposition \ref{poly gyoda} is an example showing how a cluster symmetric polynomial $F_{3,7}(\mathbf{x})$ can be constructed from the seed $\Omega_{3,7}$, which in turn naturally has a corresponding cluster symmetric equation $F_{3,7}(\mathbf{x}) = c$. 

A general way to find a non-trivial cluster symmetric polynomial from a seed $\Omega$ is as follows:
\begin{enumerate}
    \item Computer the cluster symmetric set $\mathcal{S}(\Omega)$. If the set $\mathcal{S}(\Omega)$ is nonempty, then let $S:=\{\sigma_1\mu_{s_1}, \dots, \sigma_m\mu_{s_m} \}$ be a nonempty subset of $\mathcal{S}(\Omega)$.
    \item  If there exist $i,j \in [1,m]$ with $i \ne j$, such that the following relation 
    $$\qquad 4 \geq r_{s_i}r_{s_j}\max\{|b_{{s_j}{s_i}}|,|b_{\sigma_i({s_j}){s_i}}|\}\max\{|b_{{s_i}{s_j}}|,|b_{\sigma_j({s_i}){s_j}}|,|b_{\sigma_i^{-1}({s_i}){s_j}}|,|b_{\sigma_j(\sigma_i^{-1}({s_i})){s_j}}|\}$$
    does not hold. Then by Theorem \ref{thm 4bounded}, there is no non-trivial cluster symmetric polynomial in $\mathbb{Q}[\mathbf{x}^{\pm}]^{\langle \sigma_i\mu_{s_i},\sigma_j\mu_{s_j}\rangle}$.
    \item Otherwise, use the steps in Remark \ref{rmk step of constuct} to find a non-trivial cluster symmetric polynomial about $\psi_{\sigma_1,s_1,\pi_{s_1}(\Omega)}$. That is,
        \begin{enumerate}
            \item Choose a $n$-tuple $\mathbf{d} :=  d\sum_{j \in \langle \sigma_1 \rangle(s_1)}\mathbf{e}_j$, where $d > 0$.
            \item Choose a $n$-tuple $\bm{\eta}$ that satisfies the conditions $\eta_{s_1} = \eta_{\sigma_1^{-1}(s_1)} = 2d$ and $\min\{\eta_k, \eta_{\sigma_1^{-1}(k)}\} \geq dr_{s_1}|b_{k}| \geq |\eta_k - \eta_{\sigma_1^{-1}(k)}|$ for all $k \in [1,n]$.
            \item Solve the system $HLE(\sigma_1, s_1, \pi_{s_1}(\Omega), \bm{\eta}, \mathbf{d})$. We construct a Laurent polynomial by taking the general solution of the system as coefficients, and we denote it by $F(\mathbf{x})$.
        \end{enumerate}
    \item Suppose $F(\mathbf{x})$ is of type $\frac{\bm{\eta}'}{\mathbf{d}'}$, where $\eta'_i \leq \eta_i, d'_i \leq d_i$ for all $i \in [1,n]$. Then, by Proposition \ref{is CSM}, we know that $F(\mathbf{x}) \in \mathbb{Q}[\mathbf{x}^{\pm}]^{\langle \sigma_1\mu_{s_1}\rangle}$. If there is no $\widetilde{\mathbf{d}} \in \mathbb{Z}^n$, such that for all $i \in [1,m]$, the following relations $\mathbf{d}' + \widetilde{\mathbf{d}} = \sigma_i(\mathbf{d}' + \widetilde{\mathbf{d}})$ and
        $$\qquad\qquad \eta'_{s_i} = \eta'_{\sigma_i^{-1}(s_i)} = 2(d'_{s_i}+\widetilde{d}_{s_i}) = 2(d'_{\sigma_i^{-1}(s_i)}+\widetilde{d}_{\sigma_i^{-1}(s_i)})$$
    hold, then by Theorem \ref{thm len1_T2f}, we know that $F(\mathbf{x})/\mathbf{x}^{\widetilde{\mathbf{d}}} \notin \mathbb{Q}[\mathbf{x}^{\pm}]^{\langle \sigma_1\mu_{s_1},\dots,\sigma_m\mu_{s_m}\rangle}$.
    \item Otherwise, suppose that there exists such $\widetilde{\mathbf{d}}$. Let $\widetilde{F}(\mathbf{x}) := F(\mathbf{x})/\mathbf{x}^{\widetilde{\mathbf{d}}}$. Since some of the coefficients of the polynomial $\widetilde{F}(\mathbf{x})$ are free, we determine these coefficients using the following relations
        $$\qquad\qquad \widetilde{F}(\sigma_{i}\mu_{s_i}(\mathbf{x})) = \widetilde{F}(\mathbf{x}), \quad \text{for all } i \in[1,m].$$
    If the coefficients have a solution, then we find a non-trivial cluster symmetric polynomial in $\mathbb{Q}[\mathbf{x}^{\pm}]^{\langle \sigma_1\mu_{s_1}, \dots, \sigma_m\mu_{s_m}\rangle}$.
\end{enumerate}

\section{ Orbits of solution sets under action of cluster symmetric groups } \label{sec diop}

One of the fundamental goals for Diophantine equations is to study how to find all positive integer solutions when some initial solutions are known. Concretely, for the Diophantine equation $F(\mathbf{x}) = F(\mathbf{x}_0)$, how to describe the set of positive integer solutions $\mathcal{V}_{\mathbb{Z}_{>0}}(F(\mathbf{x}) - F(\mathbf{x}_0))$? In this section, we try to discuss this question for some concrete cluster symmetric equations. Surprisingly, the solution sets of these equations have similar structures.

\subsection{Orbits of solution sets of general cluster symmetric equations}
For cluster symmetric equations, a new solution can be obtained by applying a cluster symmetric map to a solution.

\begin{prop} \label{prop G subset V}
(i) Given a cluster symmetric map $\psi_{\sigma, s, \omega_s}$. Let $F(\mathbf{x})$ be a cluster symmetric polynomial about $\psi_{\sigma, s, \omega_s}$, that is, $F(\mathbf{x}) \in \mathbb{Q}[\mathbf{x}^{\pm}]^{\langle\psi_{\sigma, s, \omega_s}\rangle}$. For an $n$-tuple ${\mathbf{x}_0} \in \mathbb{Q}_{>0}^n$, then the orbit of the initial vector ${\mathbf{x}_0}$ under the group $\langle \psi_{\sigma, s, \omega_s} \rangle$ is a subset of the set of positive rational solutions of the equation $F(\mathbf{x}) = F({\mathbf{x}_0})$, that is, the following relation
    $$\langle \psi_{\sigma, s, \omega_s} \rangle(\mathbf{x}_0) \subset \mathcal{V}_{\mathbb{Q}_{>0}}(F(\mathbf{x}) - F(\mathbf{x}_0))$$    
holds, where $\mathcal{V}_{\mathbb{K}}(H(\mathbf{x})) := \{ \mathbf{x}' \in \mathbb{K}^n \mid H(\mathbf{x}') = 0 \}$.
    
    (ii) Given a seed $\Omega$. Let $G$ be a subgroup of the complete cluster symmetric group $\overline{\mathcal{G}}(\Omega)$. Let a Laurent polynomial $F(\mathbf{x}) \in \mathbb{Q}[\mathbf{x}^{\pm}]^{G}$. Denote an $n$-tuple $\mathbf{1} = (1,\dots, 1)$. Then the orbit of the initial vector $\mathbf{1}$ under the group $G$ is a subset of the set of positive integer solutions of  equation $F(\mathbf{x}) = F(\mathbf{1})$, that is, we have
    \begin{equation} \label{G subset V}
        G(\mathbf{1}) \subset \mathcal{V}_{\mathbb{Z}_{>0}}(F(\mathbf{x}) - F(\mathbf{1})).
    \end{equation}
\end{prop}

\begin{proof} (i) Obviously, since the exchange polynomial $P_{\omega_s}(\mathbf{x}) \in \mathbb{Z}_{\geq 0}[\mathbf{x}]$ and $P_{\omega_s}(\mathbf{x}) \ne 0$.

    (ii) Let $g \in G$. Clearly, $F(g(\mathbf{1})) = F(\mathbf{1})$. By Equation \eqref{musigma=sigmamu}, the action $g$ can be written as $g = \sigma\mu_{s_m} \cdots \mu_{s_1}$ for some $\sigma \in \mathfrak{S}_n$, $s_1, \cdots, s_m \in [1,n], m \in \mathbb{Z}_{\geq 0}$. Let $\mathbf{x}' := \mu_{s_m} \cdots \mu_{s_1}(\mathbf{1})$. Then by Theorem \ref{thm PLP}, we know that $\mathbf{x}' \in \mathbb{Z}_{>0}^n$. So $g(\mathbf{1}) = \sigma(\mathbf{x}') \in \mathbb{Z}_{>0}^n$.
\end{proof}

\subsection{Uniqueness of orbits of solution sets of some special Diophantine equations}

For some special cases, the sets on both sides of Relation \eqref{G subset V} are equal, that is, the set of positive integer solutions of the Diophantine equation $F(\mathbf{x}) = F(\mathbf{1})$ is exactly the orbit of the solution $\mathbf{1}$ under the group $G$. 
That is to say, the solution set of this equation is the \textbf{unique} orbit under the actions of the group.
For example, we have the following theorem about generalized cluster algebras of rank $3$, whose proof will be given in the second subsection.

\begin{thm}  \label{thm markov-cluster sol}
    Fix $i \in [1,10]$. Let $G_{3,i}$ be the group generated by the subset $\{ \mu_1, \mu_2, \mu_3 \}$ of the cluster symmetric set $\mathcal{S}(\Omega_{3,i})$. Then the set of positive integer solutions of the Markov-cluster equation $F_{3,i}(\mathbf{x}) = F_{3,i}(\mathbf{1})$ is exactly the orbit $G_{3,i}(\mathbf{1})$, that is, we have
    \begin{equation} \label{eqn rank3 G=V}
        G_{3,i}(1,1,1) = \mathcal{V}_{\mathbb{Z}_{>0}}(F_{3,i}(x, y, z) - F_{3,i}(1,1,1)),
    \end{equation}
    where $\Omega_{3,i}, F_{3,i}$ are listed in Table \ref{tab all}.
\end{thm}

\begin{table}[htbp]
\caption{Markov-cluster polynomial $F_{3,i}$ and its seed $\Omega_{3,i} := (B_{3,i}, \mathbf{x}, R_{3,i},  \mathbf{Z}_{3,i})$.}
\begin{center}
\renewcommand{\arraystretch}{1.1} \resizebox{\textwidth}{!}{
$\begin{array}{|c|ccl|c|}
\hline
i & B_{3,i} & R_{3,i} & \qquad \quad  \mathbf{Z}_{3,i} & \begin{array}{c}
~\\[-3mm] 
F_{3,i}(x,y,z) \\[-3mm]
~ 
\end{array}\\
\hline
1 & \left[\begin{array}{rrr}
0&2&-2\\-2&0&2\\2&-2&0\end{array}\right]
& \begin{array}{c}
    r_1 = 1\\
    r_2 = 1\\
    r_3 = 1
\end{array}
& \begin{array}{l}
~\\[-3mm]
Z_1: 1+u \\
Z_2: 1+u \\
Z_3: 1+u \\[-3mm]
~
\end{array}
& \begin{array}{c}
     \underline{x^2+y^2+z^2} \\
     {xyz} 
\end{array} \\[3mm]
\hline
2&\left[\begin{array}{rrr}0&2&-1\\-2&0&1\\2&-2&0\end{array}\right]
& \begin{array}{c}
    r_1 = 1\\
    r_2 = 1\\
    r_3 = 2
\end{array}
& \begin{array}{l}
\\[-3mm]
Z_1: 1+u \\
Z_2: 1+u \\
Z_3: 1+k_3 u+u^2\\[-3mm]
~
\end{array} 
& \begin{array}{c}
     \underline{x^2+y^2+z^2+k_3xy} \\
     {xyz} 
\end{array} \\[3mm]
\hline
3&\left[\begin{array}{rrr}0&2&-1\\-1&0&1\\1&-2&0\end{array}\right]
& \begin{array}{c}
    r_1 = 2\\
    r_2 = 1\\
    r_3 = 2
\end{array}
& \begin{array}{l}
\\[-3mm]
Z_1: 1+k_1 u+u^2 \\
Z_2: 1+u \\
Z_3: 1+k_3 u+u^2 \\[-3mm]
~
\end{array}
& \begin{array}{c}
     \underline{x^2+y^2+z^2+k_1yz+k_3xy} \\
     {xyz} 
\end{array} \\[3mm]
\hline
4& \left[\begin{array}{rrr}0&1&-1\\-1&0&1\\1&-1&0\end{array}\right]
& \begin{array}{c}
    r_1 = 2\\
    r_2 = 2\\
    r_3 = 2
\end{array}
& \begin{array}{l}
\\[-3mm]
Z_1: 1+k_1 u+u^2 \\
Z_2: 1+k_2 u+u^2 \\
Z_3: 1+k_3 u+u^2 \\[-3mm]
~
\end{array} 
& \begin{array}{c}
     \underline{x^2+y^2+z^2+k_1yz+k_2zx+k_3xy} \\
     {xyz} 
\end{array} \\[3mm]
\hline
5 &\left[\begin{array}{rrr}0&1&-1\\-4&0&2\\4&-2&0\end{array}\right]
& \begin{array}{c}
    r_1 = 1\\
    r_2 = 1\\
    r_3 = 1
\end{array}
& \begin{array}{l}
\\[-3mm]
Z_1: 1+u \\
Z_2: 1+u \\
Z_3: 1+u \\[-3mm]
~
\end{array}
& \begin{array}{c}
     \underline{x^2+y^4+z^4+2xy^2+2xz^2}  \\
     {xy^2z^2} 
\end{array} \\[3mm]
\hline
6& \left[\begin{array}{rrr}0&1&-1\\-2&0&2\\2&-2&0\end{array}\right]
& \begin{array}{c}
    r_1 = 2\\
    r_2 = 1\\
    r_3 = 1
\end{array}
& \begin{array}{l}
\\[-3mm]
Z_1: 1+ku+u^2 \\
Z_2: 1+u \\
Z_3: 1+u \\[-3mm]
~
\end{array}
&  \begin{array}{c}
     \underline{x^2+y^4+z^4+2xy^2+ky^2z^2+2xz^2} \\
     {xy^2z^2} 
\end{array}  \\[3mm]
\hline
7 &\left[\begin{array}{rrr}0&1&-1\\-1&0&2\\1&-2&0\end{array}\right]
& \begin{array}{c}
    r_1 = 4\\
    r_2 = 1\\
    r_3 = 1
\end{array}
& \begin{array}{l}
\\[-3mm]
Z_1: \sum_{i=0}^{4} k_iu^i \\
Z_2: 1+u \\
Z_3: 1+u \\[-3mm]
~
\end{array}
& \begin{array}{c}
     \underline{x^2+2x(y^2+z^2)+k_1xyz+z^4Z_1(y/z)}  \\
     {xy^2z^2} 
\end{array} \\[3mm]
\hline
8 &\left[\begin{array}{rrr}0&4&-4\\-1&0&2\\1&-2&0\end{array}\right]
& \begin{array}{c}
    r_1 = 1\\
    r_2 = 1\\
    r_3 = 1
\end{array}
& \begin{array}{l}
\\[-3mm]
Z_1: 1+u \\
Z_2: 1+u \\
Z_3: 1+u \\[-3mm]
~
\end{array}
& \begin{array}{c}
     \underline{x^4+y^2+z^2+2yz}  \\
     {x^2yz} 
\end{array}\\
\hline
9 &\left[\begin{array}{rrr}0&2&-4\\-1&0&2\\1&-1&0\end{array}\right]
& \begin{array}{c}
    r_1 = 1\\
    r_2 = 2\\
    r_3 = 1
\end{array}
& \begin{array}{l}
\\[-3mm]
Z_1: 1+u \\
Z_2: 1+k_2+u^2 \\
Z_3: 1+u \\[-3mm]
~
\end{array}
& \begin{array}{c}
     \underline{x^4+k_2x^2z+y^2+z^2+2yz}  \\
     {x^2yz} 
\end{array}\\
\hline
10 &\left[\begin{array}{rrr}0&2&-2\\-1&0&1\\1&-1&0\end{array}\right]
& \begin{array}{c}
    r_1 = 1\\
    r_2 = 2\\
    r_3 = 2
\end{array}
& \begin{array}{l}
\\[-3mm]
Z_1: 1+u \\
Z_2: 1+k_2+u^2 \\
Z_3: 1+k_3+u^2 \\[-3mm]
~
\end{array}
& \begin{array}{c}
     \underline{x^4+k_3x^2y+k_2x^2z+y^2+z^2+2yz}  \\
     {x^2yz} 
\end{array}\\
\hline
\end{array}
$}
\label{tab all}
\end{center}
\end{table}

Table \ref{tab all} is from Corollary \ref{cor rank3 nonconstant}; see Remark \ref{rmk irr B} for details. Note that the Laurent polynomial $F_{3,1}(\mathbf{x})$ in the table is related to the Markov equation \eqref{eq markov}, while the other Laurent polynomials share a certain similarity with $F_{3,1}(\mathbf x)$ and originate from the cluster theory. Therefore, we name these Laurent polynomials as follows.

\begin{DEF}
\label{def mc eqn}
    For any $i \in [1,10]$, we call each Laurent polynomial $F_{3,i}$ listed in Table \ref{tab all} a \textbf{Markov-cluster polynomial}\footnote{The reason we named it so, rather than ``Markov-cluster \textbf{Laurent} polynomial'', is similar to that in  footnote \ref{footnote name}.}, and the Diophantine equation $F_{3,i}(x, y, z) = F_{3,i}(1,1,1)$ a \textbf{Markov-cluster equation}.
\end{DEF}

According to Definition \ref{def lp}, there are three types of Markov-cluster polynomial,
    \begin{itemize}
        \item[$\bullet$] $\frac{(2,2,2)}{(1,1,1)}$ type: $F_{3,1}, F_{3,2}, F_{3,3}, F_{3,4}$;
        \item[$\bullet$] $\frac{(2,4,4)}{(1,2,2)}$ type: $F_{3,5}, F_{3,6}, F_{3,7}$;
        \item[$\bullet$] $\frac{(4,2,2)}{(2,1,1)}$ type: $F_{3,8}, F_{3,9}, F_{3,10}$.
    \end{itemize}

    The above classification of the Markov-cluster polynomials motivates us to define the height function, which plays an important role in the proof of Proposition \ref{prop Kaufman} and Proposition \ref{prop sol of gyoda} in the next two subsections.
    \begin{DEF} \label{def h ftn}
        Given a Laurent polynomial $F(\mathbf{x})$ of type $\frac{\bm{\eta}}{\mathbf{d}}$ in $\mathbb{Q}[\mathbf{x}^{\pm}]$. The \textbf{height function} $h_F$ of the Laurent polynomial $F(\mathbf{x})$ is $h_F(\mathbf{x}) := \max_{1 \leq i \leq n}\{ x_i^{d_i}\}$. We also call the function $h_F$ the \textbf{height function} of the Diophantine equation $F(\mathbf{x}) = F(\mathbf{x}_0)$, where $\mathbf{x}_0 \in \mathbb{Q}_{>0}^n$.
    \end{DEF}

It is interesting for us to note two conjectures about the Markov equation $F_{3,1}(\mathbf{x}) = F_{3,1}(\mathbf{1})$ as follows. \color{black}

\begin{conj}[(Markov Uniqueness conjecture, \cite{MUC, Markov100})]
\label{conj muc}
    Suppose $(a,b,c)$ and $(a',b',c')$ are positive integer solutions of the Markov equation $x^2+y^2+z^2 = 3xyz$. If $\max(a,b,c) = \max(a',b',c')$, then there exists a permutation $\sigma \in \mathfrak{S}_3$, such that $(a,b,c) = \sigma(a',b',c')$.
\end{conj}

This conjecture is a century old one and is still open so far.

The second is the strong approximation conjecture posed by Baragar \cite{Baragar}, which conjectures that the Markov graph over the finite field $\mathbb{F}_p$ is connected for any prime $p$. The first major progress on this conjecture is the work of Bourgain, Gamburd, and Sarnak \cite{Sarnak}, and then W. Y. Chen proved that the conjecture holds for all but finitely many primes in \cite{w.chen}.

From our perspective, since the Markov-cluster equations possess a solution structure similar to the Markov equation, we believe that it is worth studying the analogous versions of the above two conjectures for these Markov-cluster equations.

\subsection{Proof of Theorem \ref{thm markov-cluster sol}}

\subsubsection{Solutions of the equation $F_{3,10}(\mathbf{x}) = F_{3,10}(\mathbf{1})$}

We obtain the Diophantine equation \eqref{eq kaufman}, by substituting $(X,Y,Z)$ for $(x_1,x_2,x_3)$ in the equation $F_{3,10}(x_1,x_2,x_3) = F_{3,10}(1,1,1)$ defined in Table \ref{tab all}. Clearly, the set of positive integer solutions of this equation is exactly the set $\mathcal{V}_{\mathbb{Z}_{>0}}(F_{3, 10}(x, y, z) - F_{3, 10}(1,1,1))$. So we solve the Diophantine equation. Note that the following proposition states that Relation \eqref{eqn rank3 G=V} holds when $i = 8,9,10$.

\begin{prop} \label{prop Kaufman} 
    For any non-negative integers $k_2, k_3$. The set of positive integer solutions of the Diophantine equation
    \begin{equation} \label{eq kaufman}
        X^4 + k_3X^2Y+k_2X^2Z + Y^2 + Z^2 + 2YZ = (5+k_2+k_3){X^2YZ}
    \end{equation}
    is exactly the orbit of the  solution $(1,1,1)$ under the group $G := \langle \mu_1, \mu_2, \mu_3 \rangle$, where
    \begin{eqnarray*}
        \mu_1(X,Y,Z) &=& \bigg(\frac{Y + Z}{X}, Y, Z\bigg),\\
        \mu_2(X,Y,Z) &=& \bigg(X, \frac{X^4+k_2X^2Z+Z^2}{Y}, Z\bigg),\\
        \mu_3(X,Y,Z) &=& \bigg(X, Y, \frac{X^4+k_3X^2Y+Y^2}{Z}\bigg).
    \end{eqnarray*}
\end{prop}

\begin{proof}
    By Example \ref{eg existence of F}(iv) and Proposition \ref{prop G subset V}(ii), we know that $g(\mathbf{1})$ is a positive integer solution for any $g \in G$. So we only need to prove that for a positive integer solution $(x,y,z)$ of Equation \eqref{eq kaufman}, there exists $g \in G$, such that $g(1,1,1) = (x, y, z)$. To do this, we define a height function $h(X,Y,Z) := \max\{X^2, Y, Z\}$.
    
    We prove it in three steps.
    
    \textbf{Step 1.} Let $(x,y,z)$ be a positive integer solution. Suppose that at least two of $x^2,y,z$ are equal. We claim that $(x,y,z)$ must be one of the four solutions $(1,1,1)$, $\mu_1(1,1,1) =(2,1,1), \mu_2(1,1,1) = (1,2+k_2,1), \mu_3(1,1,1) = (1,1,2+k_3)$.

    (i) If $y = z$, then Equation \eqref{eq kaufman} becomes
                $$ x^4  - x^2vy +  4y^2 = 0.$$
    where $v = (5+k_2+k_3)y - (k_2+k_3)$. Since the discriminant $\Delta = y^2(v^2 - 16)$ must be a square, we let $v^2 - 16 = t^2$. Then, we have $(v - t)(v + t) = 16$. Since $v \geq 5$ and $v \pm t \mid 16$, it is easy to check $v = 5$. So $y = 1$ and $x = 1$ or $2$. Then $(x,y,z) = (1,1,1)$ or $(2,1,1)$.

    (ii) If $x^2 = y \ne z$, then Equation \eqref{eq kaufman} becomes
        $$(Az - B)y^2 - Cyz - z^2 = 0,$$
    where $A = 5+k_2+k_3$, $B = 2 + k_3$, $C = 2+k_2$. Clearly, $Az - B \ne 0$. By substituting $w = \frac{2(Az - B)y - Cz}{z}$ into the above equation, we have $w^2 = 4Az + C^2 - 4B$. Then $w$ is an integer. So $z \mid 2By$. Let $2By = zt$, where $t$ is a positive integer. Hence, the above equation becomes
    $$(2Ay -t - 2C)t = 4B.$$
    So $t \mid 4B$. Let $4B = ts$, then $2Ay - t - 2C = s$. Therefore, we have
    $$y = \frac{t+s+2C}{2A} \leq \frac{1 + 4B + 2C}{2A} = \frac{2(B+C+1)+2B-1}{2A} < 2.$$
Hence $y = 1$. It is easy to check that $z = 2+k_3$ or $z=1$(discard). So $(x,y,z) = (1,1,2+k_3)$.

    (iii) If $x^2 = z \ne y$. Similar to (ii), we know $(x,y,z) = (1,2+k_2,1)$.

    \textbf{Step 2.} Let $(x,y,z)$ be a positive integer solution. Suppose $x^2,y,z$ are not equal to each other. We claim that there exists $i \in [1,3]$ such that $h(\mu_i(x,y,z)) < h(x,y,z)$ and $\mu_i(x,y,z)$ is a positive integer solution.
    
    (i) If $h(x,y,z) = x^2$. Let $x' = (y+z)/x$ and $\omega = \max\{y,z\}$. Consider the function
    \begin{eqnarray*}
        f(\lambda) &:=& \lambda y^2z(F_{3,10}(\lambda,y,z)-F_{3,10}(1,1,1))\\
         &=& \lambda^2 - (Ayz+k_3y+k_2z)\lambda + (y+z)^2,
    \end{eqnarray*}
    where $A = 5+k_2+k_3$. Clearly, $f(x^2) = f(x'^2) = 0$ and $x'$ is a positive integer. Since $yz \geq \omega > (y+z)/2$, we have
        $$f(\omega) = \omega^2 + (k_3y+k_2z-Ayz)\omega + (y+z)^2 < \omega^2(1 + k_3 + k_2 - A + 4)= 0.$$
    Then $x'^2 < \omega < x^2$. Hence, $h(\mu_1(x,y,z)) = h(x',y,z) = \omega < x^2 = h(x,y,z)$.

    (ii) If $h(x,y,z) = y$. Let $y' = (x^4 + k_2x^2z + z^2)/y$ and $\omega = \max\{x^2, z\}$. Consider the function
    \begin{eqnarray*}
        f(\lambda) &:=& x\lambda^2z(F_{3,10}(x,\lambda,z)-F_{3,10}(1,1,1))\\
         &=& \lambda^2 + (2z + k_3x^2 - Ax^2z)\lambda + x^4 + k_2x^2z + z^2,
    \end{eqnarray*}
    where $A = 5+k_2+k_3$. Clearly, $f(y) = f(y') = 0$ and $y'$ is a positive integer. Since $x^2z \geq \omega$ and $ x^4 + z^2 < 2\omega^2$, we have
        $$f(\omega) = \omega^2 + (2z + k_3x^2 - Ax^2z)\omega + x^4 +  k_3x^2z + z^2 < \omega^2(1 + 2 + k_3 - A + 2 + k_3) = 0.$$
    Then $y' < \omega < y$. Hence, $h(\mu_2(x,y,z)) = h(x,y',z) = \omega < y = h(x,y,z)$.

    (iii) If $h(x,y,z) = z$. Similarly to (ii), we know $h(\mu_3(x,y,z)) < h(x,y,z)$.

    \textbf{Step 3.} Let $(x,y,z)$ be a positive solution. If $x^2,y,z$ are not equal to each other, then by step 2, we can find a finite sequence of $\mu_{s_1}, \cdots, \mu_{s_m}$ such that $$h(x,y,z) > h(\mu_{s_1}(x,y,z)) > h(\mu_{s_2}\mu_{s_1}(x,y,z)) \ldots > h(x_0, y_0, z_0),$$
    where $(x_0, y_0, z_0):= \mu_{s_m} \cdots \mu_{s_1}(x,y,z)$ and two of $x_0^2, y_0, z_0$ are equal. Then by step 1, we know that there exists $s_{m+1} \in [1,3]$ such that $\mu_{s_{m+1}}(x_0,y_0,z_0) = (1,1,1)$. Hence $(x,y,z) = \mu_{s_1} \cdots \mu_{s_{m+1}}(1,1,1)$. If at least two of $x^2,y,z$ are equal, in Step 1, we know that there exists $\mu_i \in G$, such that $(x,y,z) = \mu_i(1,1,1)$.
\end{proof}

There is another way to prove the above proposition. To do this, we need the following proposition which proves that the relation $G_{3,4}(1,1,1) = \mathcal{V}_{\mathbb{Z}_{>0}}(F_{3,4}(x, y, z) - F_{3,4}(1,1,1))$ holds in Table \ref{tab all}.
\begin{prop}[({\cite[Theorem 1]{Gyoda}})]
\label{prop geneMarkov}
        The set of positive integer solutions of the Diophantine equation
        \begin{equation} \label{eq geneMarkov}
        X^2 + Y^2 + Z^2 + k_1YZ + k_2ZX + k_3XY = (3+k_1+k_2+k_3)XYZ 
        \end{equation}
        is exactly the orbit of the  solution $(1,1,1)$ under the group $\widetilde{G} := \langle \widetilde{\mu}_1, \widetilde{\mu}_2, \widetilde{\mu}_3 \rangle$, where
        \begin{eqnarray*}
            \widetilde{\mu}_1(X,Y,Z) &=& \bigg(\frac{Y^2 + k_1YZ + Z^2}{X}, Y, Z\bigg),\\
            \widetilde{\mu}_2(X,Y,Z) &=& \bigg(X, \frac{X^2 + k_2XZ + Z^2}{Y}, Z\bigg),\\
            \widetilde{\mu}_3(X,Y,Z) &=& \bigg(X, Y, \frac{X^2 + k_3XY + Y^2}{Z}\bigg).
        \end{eqnarray*}
    \end{prop}

\begin{proof}[Another proof of Proposition \ref{prop Kaufman}] 
    We only need to prove that for a positive integer solution $(x,y,z)$ of Equation \eqref{eq kaufman}, there exists $g \in G$, such that $g(1,1,1) = (x, y, z)$. Denote a map $\varphi(x,y,z) := (\sqrt{x},y,z)$. Consider the case $k_1 = 2$ in Proposition \ref{prop geneMarkov}, it is easy to check $\mu_i\varphi(x,y,z) = \varphi\widetilde{\mu}_i(x,y,z)$ for $i = 1,2,3$.
    
    Let $(x,y,z)$ be a positive integer solution of Equation \eqref{eq kaufman}.  Clearly, $(x^2,y,z)$ is a positive integer solution of equation $X^2 + Y^2 + Z^2 + 2YZ = 5XYZ$. Then, by Proposition \ref{prop geneMarkov}, there exists $\widetilde{\mu}_{s_1}\cdots\widetilde{\mu}_{s_q}\in\widetilde{G}$, such that $(x^2,y,z) = \widetilde{\mu}_{s_1}\cdots\widetilde{\mu}_{s_q}(1,1,1)$. Hence,
    \begin{eqnarray*}
        (x,y,z) = \varphi(x^2,y,z) = \varphi\widetilde{\mu}_{s_1}\cdots\widetilde{\mu}_{s_q}(1,1,1) =  \mu_{s_1}\cdots\mu_{s_q}\varphi(1,1,1) = \mu_{s_1}\cdots\mu_{s_q}(1,1,1).
    \end{eqnarray*}
\end{proof}

\subsubsection{Solutions of the equation $F_{3,7}(\mathbf{x}) = F_{3,7}(\mathbf{1})$}

By substituting $(X,Y,Z)$ for $(x_1,x_2,x_3)$ in the equation $F_{3,7}(x_1,x_2,x_3) = F_{3,7}(1,1,1)$ which is defined in Table \ref{tab all}, we obtain the Diophantine equation \eqref{eq gyoda}. Clearly, the set of positive integer solutions of this equation is exactly the set $\mathcal{V}_{\mathbb{Z}_{>0}}(F_{3, 7}(x, y, z) - F_{3, 7}(1,1,1))$. So we solve the Diophantine equation.

\begin{prop} \label{prop sol of gyoda}
    Let $k_1, k_2 \in \mathbb{Z}_{\ge 0}$. The set of positive integer solutions of the following Diophantine equation
    \begin{equation} \label{eq gyoda}
        X^2 + Y^4 + Z^4 + 2X(Y^2 + Z^2)  + k_1YZ(X+Y^2+Z^2) +  k_2Y^2Z^2 = (7+3k_1+k_2){XY^2Z^2}
    \end{equation}
    is exactly the orbit of the  solution $(1,1,1)$ under the group $G := \langle \mu_1, \mu_2, \mu_3 \rangle$, where
    \begin{eqnarray*}
        \mu_1(X,Y,Z) &=& \bigg(\frac{Y^4 + k_1Y^3Z+k_2Y^2Z^2 + k_1YZ^3+Z^4}{X}, Y, Z\bigg),\\
        \mu_2(X,Y,Z) &=& \bigg(X, \frac{X+Z^2}{Y}, Z\bigg),\\
        \mu_3(X,Y,Z) &=& \bigg(X, Y, \frac{X+Y^2}{Z}\bigg).
    \end{eqnarray*}
\end{prop}

\begin{proof}
   By Proposition \ref{poly gyoda} and Proposition \ref{prop G subset V}(ii), we know that $g(\mathbf{1})$ is a positive integer solution for any $g \in G$. So we only need to prove that for a positive integer solution $(x,y,z)$ of Equation \eqref{eq gyoda}, there exists $g \in G$, such that $g(1,1,1) = (x, y, z)$. To do this, we define a height function $h(X,Y,Z) := \max\{X, Y^2, Z^2\}$.
   
   We prove it in three steps.
    
\textbf{Step 1.} Let $(x,y,z)$ be a positive integer solution. Suppose that at least two of $x,y^2,z^2$ are equal. We claim that $(x,y,z)$ must be one of the four solutions $(1,1,1)$, $\mu_1(1,1,1) =(2+2k_1+k_2,1,1), \mu_2(1,1,1) = (1,2,1), \mu_3(1,1,1) = (1,1,2)$.

(i) If $y=z$. Equation \eqref{eq gyoda} becomes 
    $$(Ax - B)y^4 - Cxy^2 - x^2 = 0,$$
where $A = 7+3k_1+k_2$, $B = 2 + 2k_1 + k_2$, $C = 4+k_1$. Clearly, $Ax - B \ne 0$. By substituting $w = \frac{2(Ax - B)y^2 -Cx}{x}$ into the above equation, we have $w^2 = 4Ax + C^2 - 4B$. Then $w$ is an integer. So $x \mid 2By^2$. Let $2By^2 = xt$, where $t$ is a positive integer. Hence, the above equation becomes
    $$(2Ay^2 -t - 2C)t = 4B.$$
So $t \mid 4B$. Let $4B = ts$, then $2Ay^2 - t - 2C = s$. Therefore, we have
    $$y^2 = \frac{t+s+2C}{2A} \leq \frac{1 + 4B + 2C}{2A} = \frac{2(B+C+1)+2B-1}{2A} < 2.$$
Hence $y = 1$. It is easy to check $x = 1$ or $x = A-C-1 = 2+2k_1+k_2$. So
$(x,y,z) = (1,1,1)$ or $(2 + 2k_1 + k_2,1,1)$.

(ii) If $x = y^2 \ne z^2$.  Equation \eqref{eq gyoda} becomes 
\begin{eqnarray}\label{eq x=y2}
    z^4 + k_1yz^3 + (k_2+2 - Ay^2)y^2z^2 + 2k_1y^3z + 4y^4 = 0.
\end{eqnarray}
where $A = 7+3k_1+k_2$. By substituting $w = k_1y + 2z + 4y^2/z$ into the above equation, we have $w^2 = y^2(4Ay^2 + k_1^2 - 4k_2 +8)$. Then $w$ is an integer. So $z \mid 4y^2$. Let $4y^2 = zt$, where $t$ is a positive integer. Hence, the above equation becomes
    $$\bigg(\frac{4y}{t}\bigg)^4 + k_1\bigg(\frac{4y}{t}\bigg)^3 + (k_2 + 2 - Ay^2)\bigg(\frac{4y}{t}\bigg)^2 + 2k_1\bigg(\frac{4y}{t}\bigg) + 4 = 0.$$
Since it is a monic polynomial with integer coefficients, we know that $\frac{4y}{t}$ is a positive integer and $\frac{z}{y} = \frac{4y}{t} \mid 4$. So $z=2y$ or $z=4y$. If $z = 2y$, then $y=1$, $z=2$. If $z = 4y$, then $4Ay^2 = 4^3 + 4^2k_1 + 4(k_2+2) + 4k_1 + 1$ is odd, it is impossible. Hence $(x,y,z) = (1,1,2)$.

(iii) If $x = z^2 \ne y^2$. Similarly to (ii), we know $(x,y,z) = (1,2,1)$.

\textbf{Step 2.} For a positive integer solution $(x,y,z)$.
Suppose $x,y^2,z^2$ are not equal to each other. We claim that there exists $i \in [1,3]$ such that $h(\mu_i(x,y,z)) < h(x,y,z)$ and $\mu_i(x,y,z)$ is a positive integer solution.

(i) If $h(x,y,z) = x$. Denote $x' = ({y^4+k_1y^3z+k_2y^2z^2 +k_1yz^3+z^4})/{x}$. Let $w := \max \{y^2, z^2\}$. Then $y^2z^2 = w\min\{y^2,z^2\} \geq w$. We denote a function 
\begin{eqnarray*}
    f(\lambda) &:=& \lambda y^2z^2(F_{3,7}(\lambda,y,z) - F_{3,7}(1,1,1))\\ 
        &=& \lambda^2 + (2y^2 + 2z^2 + k_1yz - (7+3k_1+k_2)y^2z^2)\lambda + y^4 + z^4 + k_1yz^3 + k_2y^2z^2 + k_1y^3z.
\end{eqnarray*}
Clearly, $f(x') = f(x) = 0$, $x'$ is a positive integer and
$$f(w) \leq w^2 + (4w + k_1w - (7+3k_1+k_2)w)w +  (2+2k_1+k_2)w^2 = 0.$$ 
We know that $x' \leq w$. Hence $h(\mu_1(x,y,z)) = h(x', y, z) = w < x = h(x, y, z)$.

(ii) If $h(x,y,z) = y^2$. Denote $y' = \frac{x+z^2}{y}$. Let $w := \sqrt{\max \{x, z^2\}}$. Then we have 
$$x^2 \leq w, z \leq w, xz^2 \geq w^2 \text{ and } y > w > 0.$$
We denote a function
\begin{eqnarray*}
    f(\lambda) &:=& x\lambda^2z^2(F_{3,7}(x,\lambda,z) - F_{3,7}(1,1,1))\\
    &=& \lambda^4 + k_1z\lambda^3 + (2x + k_2z^2 - (7+3k_1+k_2)xz^2)\lambda^2 + k_1z(x+z^2)\lambda + (x+z^2)^2.
\end{eqnarray*}
Clearly, $f(y') = f(y) = 0$, $y'$ is a positive integer, $f(0) > 0$, $f(-y) < 0$, $f(-\infty) > 0$ and
\begin{eqnarray*}
    f(w) = w^4 + k_1zw(x + z^2 + w^2) + (2x + k_2z^2 - (7+3k_1+k_2)xz^2)w^2 + (x+z^2)^2\\
        \leq w^4 + k_1w^4 + (2w^2 + k_2w^2 - (7+3k_1+k_2)w^2)w^2 + 2k_1w^4 + 4w^4 = 0.~~~\quad
\end{eqnarray*}
and $f(-w) = f(w) - 2k_1zw^3 - 2k_1z(x+z^2)w < 0.$
Then there exists $y_1 \in (0, w]$, $y_2 \in (-w, 0)$, $y_3 \in (-\infty, -y)$ such that $f(y_1) = f(y_2) = f(y_3) = 0$. Since $y' >0 $, we know that $y' = y_1 \leq w$. Hence $h(\mu_2(x,y,z)) = h(x, y', z)  = w < y = h(x,y,z)$.

(iii) If $h(x,y,z) = z^2$. Denote $z' = \frac{x+y^2}{z}$. Similarly to (ii), we know that $h(\mu_3(x,y,z)) = h(x, y, z') < h(x, y, z)$.

\textbf{Step 3.} Let $(x,y,z)$ be a positive solution. If $x,y^2,z^2$ are not equal to each other, then by Step 2, we can find a finite sequence of $\mu_{s_1}, \cdots, \mu_{s_m}$ such that $h(x,y,z) > h(\mu_{s_1}(x,y,z)) > h(\mu_{s_2}\mu_{s_1}(x,y,z)) \ldots > h(x_0, y_0, z_0)$, where $(x_0, y_0, z_0):= \mu_{s_m} \cdots \mu_{s_1}(x,y,z)$ and two of $x_0, y_0^2, z_0^2$ are equal. Then by Step 1, we know that there exists $s_{m+1} \in [1,3]$ such that $\mu_{s_{m+1}}(x_0,y_0,z_0) = (1,1,1)$. Hence $(x,y,z) = \mu_{s_1} \cdots \mu_{s_{m+1}}(1,1,1)$. 
If at least two of $x,y^2,z^2$ are equal, by Step 1, we know that there exists $g \in G$, such that $(x,y,z) = g(1,1,1)$.

\end{proof}

Lastly, we can finish the proof of the main theorem of this section as follows.

\begin{proof}[Proof of Theorem \ref{thm markov-cluster sol}]




        
    For the case $i = 1$, it was proved by Markov in \cite{Markov}.
    For the case $i = 5$, it was proved by Lampe in \cite{Lampe}.
    For the case $i = 2,3,4,6$, it was proved by Gyoda and Matsushita in \cite{Gyoda}. 
    For the case $i = 7$, it is true, since Proposition \ref{prop sol of gyoda}.
    For the case $i = 8, 9, 10$, it is true, since Proposition \ref{prop Kaufman}.
    
\end{proof}

\section{Cluster symmetric sets of Diophantine equations}
\label{sec CSP}

In the previous section, we showed that cluster symmetry maps play an important role in solving Diophantine equations. If a Diophantine equation is invariant under a cluster symmetry map, we can obtain new solutions from the old solutions of the equation through the map. Therefore, the key questions are how to determine whether a given Diophantine equation has a cluster symmetric map and how to find all of its cluster symmetric maps. Furthermore, can we relate a given Diophantine equation to a generalized cluster algebra? In this section, we answer these questions.

Note that any Diophantine equation can be expressed as $F(\mathbf{x}) = c$, where $F(\mathbf{x})$ is a Laurent polynomial. Therefore, in this section, we study the cluster symmetry of Laurent polynomials.

\subsection{Cluster symmetric maps of a Laurent polynomial}

We collect all cluster symmetric maps of a given Laurent polynomial into a set.
\begin{DEF}\label{def S1 of F}
    Given a Laurent polynomial $F(\mathbf{x}) \in \mathbb{Q}[\mathbf{x}^{\pm}]$. Suppose $F(\mathbf{x})$ is of type $\frac{\bm{\eta}}{\mathbf{d}}$. The \textbf{cluster symmetric set $\mathcal{S}(F)$} of $F(\mathbf{x})$ is defined as 
    $$\mathcal{S}(F) := \{ \psi_{\sigma,s,\omega_s} \mid F(\psi_{\sigma,s,\omega_s}(\mathbf{x})) = F(\mathbf{x}), \eta_s  \ne 0\}.$$
    The \textbf{cluster symmetric group } of $F(\mathbf{x})$ be the group $\mathcal{G}(F)$ generated by the set $\mathcal{S}(F)$. 
    When the cluster symmetric set $\mathcal{S}(F)$ is nonempty, we call $F(\mathbf{x})$ has \textbf{cluster symmetry}.
\end{DEF}

\begin{rmk} \label{why eta ne 0}
    We require condition $\eta_s \ne 0$ because if $\eta_s = 0$, then the cluster symmetric map $\psi_{\sigma,s,\omega_s}$ actually only serves as the permutation $\sigma$ which can be obtained directly from the symmetries of the Laurent polynomial. For example, we consider the Laurent polynomial $F(x_1, x_2, x_3) := x_2^2 + x_3^2$. It is easy to check that $F(\mathbf{x})$ is invariant under the cluster symmetric map $\psi_{\sigma_{(23)}, 1, \omega_1}$, where $\omega_1$ is an arbitrary seedlet.
\end{rmk}

When $\mathcal{S}(F) \ne \emptyset$, by Definition \ref{def 1-cspoly}, $F(\mathbf{x})$ is a cluster symmetric polynomial, the Diophantine equation $F(\mathbf{x}) = c$ is a cluster symmetric equation.  Our goal is to obtain the set $\mathcal{S}(F)$, but we can do more than that. Some Laurent polynomials do not have cluster symmetry, but when they are adjusted, the resulting new Laurent polynomials may have cluster symmetry. For example, the polynomials $x_1^2 + x_2^2 + x_3^2$ and $x_1^2 + x_2^2 + x_3^2 + cx_1x_2x_3$ do not have any cluster symmetry. However, the Laurent polynomials $\frac{x_1^2 + x_2^2 + x_3^2}{x_1x_2x_3}$ and $\frac{x_1^2 + x_2^2 + x_3^2}{x_1x_2x_3} + c$ do, since they are, respectively, the Laurent polynomials $F_1(\mathbf{x})$ and $F_1(\mathbf{x}) + c$ in Equation \eqref{eq Laurent Markov} associated with the Markov equation \eqref{eq markov}. 
Based on this observation, we add the tuple $\widetilde{\mathbf{d}}$ in the following definition.

\begin{DEF}\label{def csp}
        Let $F(\mathbf{x})$ be a Laurent polynomial in $\mathbb{Q}[\mathbf{x}^{\pm}]$. For any 
        cluster symmetric map $\psi_{\sigma, s, \omega_s}$ and $n$-tuple $\widetilde{\mathbf{d}} \in \mathbb{Z}^n$, the pair $(\psi_{\sigma, s, \omega_s}, \widetilde{\mathbf{d}})$ is the \textbf{cluster symmetric pair of} $F(\mathbf{x})$, if 
        $\widetilde{F}(\psi_{\sigma, s, \omega_s}(\mathbf{x})) = \widetilde{F}(\mathbf{x})$,
        where $\widetilde{F}(\mathbf{x}) := \mathbf{x}^{-\widetilde{\mathbf{d}}} F(\mathbf{x})$. 
\end{DEF}

There is a class of cluster symmetric pairs that can be constructed directly.
\begin{prop} \label{RMK non-trivial}
    Given a $\frac{\bm{\eta}}{\mathbf{d}}$ type Laurent polynomial $F(\mathbf{x}) := \frac{T(\mathbf{x})}{\mathbf{x}^{\mathbf{d}}} \in \mathbb{Q}[\mathbf{x}^{\pm}]$.
    
    (i) Let $I := \{ i \in [1,n] \mid \eta_i = 0\}$. Fix $s \in I$. For an arbitrary seedlet $\omega_s$, a permutation $\sigma \in \{ \sigma \in \mathfrak{S}_n \mid \sigma^{-1}(s) \in I, T(\sigma(\mathbf{x})) = T(\mathbf{x}) \}$ and an $n$-tuple $\widetilde{\mathbf{d}} \in \{\widetilde{\mathbf{d}} \in \mathbb{Z}^n \mid \sigma(\mathbf{d} + \widetilde{\mathbf{d}}) = \mathbf{d} + \widetilde{\mathbf{d}}, d_s + \widetilde{d}_s = 0 \}$. Then, the pair $(\psi_{\sigma, s, \omega_s}, \widetilde{\mathbf{d}})$ is a cluster symmetric pair of $F(\mathbf{x})$. In this case, we call the pair $(\psi_{\sigma, s, \omega_s}, \widetilde{\mathbf{d}})$ is a \textbf{trivial cluster symmetric pair of} $F(\mathbf{x})$.

    (ii) If the pair $(\psi_{\sigma, s, \omega_s}, \widetilde{\mathbf{d}})$ is an \textbf{non-trivial} cluster symmetric pair of $F(\mathbf{x})$. Then $\eta_s = \eta_{\sigma^{-1}(s)} \ne 0$. We denote $\mathcal{M}(F)$ to be the set of non-trivial cluster symmetric pairs of $F(\mathbf{x})$.
\end{prop}

\begin{proof}
    Let $\widetilde{F}(\mathbf{x}) := \mathbf{x}^{-\widetilde{\mathbf{d}}}F(\mathbf{x}) = \frac{T(\mathbf{x})}{\mathbf{x}^{\mathbf{d}+ \widetilde{\mathbf{d}}}}$ and $t = \sigma^{-1}(s)$. If $\eta_s = \eta_{\sigma^{-1}(s)} = 0$, it is clear that
    \begin{equation}\label{eta=0 FF}
        T(\psi_{\sigma, s, \omega_s}(\mathbf{x})) = T\bigg(\bigg(\sigma(\mathbf{x})\bigg)\bigg|_{x_s\gets \frac{P_{\omega_s}(\mathbf{x})}{x_s}}\bigg) = T(\sigma(\mathbf{x})).
    \end{equation}

    (i) By the above equation and Equation \eqref{sigma_mu_d}, we have 
        $$\widetilde{F}(\psi_{\sigma, s, \omega_s}(\mathbf{x})) =  \frac{T(\psi_{\sigma, s, \omega_s}(\mathbf{x}))}{(\psi_{\sigma, s, \omega_s}(\mathbf{x}))^{\mathbf{d}+ \widetilde{\mathbf{d}}}} = \frac{T(\sigma(\mathbf{x}))}{\mathbf{x}^{\sigma^{-1}(\mathbf{d}+ \widetilde{\mathbf{d}})}\left(\frac{P_{\omega_s}(\mathbf{x})}{x_s^2}\right)^{d_{t}+\tilde{d}_{t}}} = \frac{T(\mathbf{x})}{\mathbf{x}^{\mathbf{d}+ \widetilde{\mathbf{d}}}} = \widetilde{F}(\mathbf{x}).$$
    Hence, the pair $(\psi_{\sigma, s, \omega_s}, \widetilde{\mathbf{d}})$ is a cluster symmetric pair of $F(\mathbf{x})$.
    
    (ii) Since $(\psi_{\sigma, s, \omega_s}, \widetilde{\mathbf{d}})$ is a cluster symmetric pair, by Proposition \ref{thm len1_T2f}, we know that $\eta_s = \eta_{t} = 2(d_s + \tilde{d}_s) = 2(d_{t} + \tilde{d}_t)$ and $\sigma(\mathbf{d} + \widetilde{\mathbf{d}}) = \mathbf{d} + \widetilde{\mathbf{d}}$. Assume $\eta_s = 0$. Then by equations \eqref{eta=0 FF} and \eqref{sigma_mu_d}, we have 
        $$\frac{T(\mathbf{x})}{\mathbf{x}^{\mathbf{d}+ \widetilde{\mathbf{d}}}} =  \frac{T(\psi_{\sigma, s, \omega_s}(\mathbf{x}))}{(\psi_{\sigma, s, \omega_s}(\mathbf{x}))^{\mathbf{d}+ \widetilde{\mathbf{d}}}} = \frac{T(\sigma(\mathbf{x}))}{\mathbf{x}^{\sigma^{-1}(\mathbf{d}+ \widetilde{\mathbf{d}})}\left(\frac{P_{\omega_s}(\mathbf{x})}{x_s^2}\right)^{d_{t}+\tilde{d}_{t}}} = \frac{T(\sigma(\mathbf{x}))}{\mathbf{x}^{\mathbf{d}+ \widetilde{\mathbf{d}}}}.$$
    So $T(\sigma(\mathbf{x})) = T(\mathbf{x})$. Hence, the cluster symmetric pair $(\psi_{\sigma, s, \omega_s}, \widetilde{\mathbf{d}})$ is trivial, it is a contradiction. Hence, $\eta_{s} = \eta_t \ne 0$.
\end{proof}

The non-trivial ones can be found by following algorithm.

\begin{alg} \label{algo diop2cs}   
    Given a Laurent polynomial $F(\mathbf{x}) \in \mathbb{Q}[\mathbf{x}^{\pm}]$, we can obtain the set $\mathcal{M}(F)$ of the non-trivial cluster symmetric pairs of $F(\mathbf{x})$ by the following routine:
    \begin{enumerate}
        \item[(1)] Set $\mathcal{M}(F) \gets \emptyset$. Suppose $F(\mathbf{x}) := \frac{T(\mathbf{x})}{\mathbf{x}^{\mathbf{d}}}$ is of type $\frac{\bm{\eta}}{\mathbf{d}}$. Let $S := \{ i \in [1, n] \mid \eta_i \text{ is even and nonzero}\}$ be a set of directions. 
        \item[(2)] If the set $S$ is empty, then stop.
        \item[(3)] Pick a direction $s \in S$, and compute the set of permutations
        $\Sigma_s := \{ \sigma \in \mathfrak{S}_n \mid \sigma^{-1}(s) \in S, \eta_s = \eta_{\sigma^{-1}(s)} \}.$
        \begin{enumerate}
            \item If the set $\Sigma_s$ is empty, then set $S \gets S\setminus\{s\}$ and go to subroutine (2).
            \item Pick a permutation $\sigma \in \Sigma_s$. Let $t = \sigma^{-1}(s)$ and $\widetilde{\mathbf{d}} := (\widetilde{d}_1, \cdots, \widetilde{d}_n)$ be the $n$-tuple of free variables satisfying $\sigma(\mathbf{d} + \widetilde{\mathbf{d}}) = \mathbf{d} + \widetilde{\mathbf{d}}$, $2(d_s+\widetilde{d}_s) = \eta_s$. 
            \begin{enumerate}
                \item Denote $\eta_s$ polynomials $f_{s,0}(\mathbf{x}), \cdots, f_{s,\eta_s}(\mathbf{x})$, such that $T(\mathbf{x}) = \sum_{i=0}^{\eta_s} f_{s, i}(\mathbf{x}) \mathbf{x}^{i\mathbf{e}_s}$, and $\deg^s f_{s,i}(\mathbf{x}) = 0$ for any $i \in [0, \eta_s]$.
                \item Denote $\eta_t$ polynomials $f_{t,0}(\mathbf{x}), \cdots, f_{t,\eta_t}(\mathbf{x})$, such that $T(\mathbf{x}) = \sum_{i=0}^{\eta_t} f_{t, i}(\mathbf{x}) \mathbf{x}^{i\mathbf{e}_t}$, and $\deg^t f_{t,i}(\mathbf{x}) = 0$ for any $i \in [0, \eta_t]$.
                \item Pick an integer $k_0 \in [0,\eta_s/2]$, such that $f_{t,k_0}(\mathbf{x}), f_{s,\eta_s - k_0}(\mathbf{x}) \ne 0$. Let $P(\mathbf{x}) := ({f_{t,k_0}(\sigma(\mathbf{x}))}/{f_{s,\eta_s - k_0}(\mathbf{x})})^{-(d_s  + \tilde{d}_s - k_0)}$. For any $i \in [0,\eta_s]$, check the relation
                $f_{t,i}(\sigma(\mathbf{x})) = f_{s,\eta_s-i}(\mathbf{x})P^{d_s+\widetilde{d}_s - i}(\mathbf{x})$
                hold. If this subroutine fails, then set $\Sigma_s \gets \Sigma_s\setminus\{\sigma\}$ and go to subroutine (b).
                \item Set $\mathcal{W} \gets \{ (\mathbf{b}, r) \in  \mathbb{Z}^n \times \mathbb{Z}_{>0} \mid b_s = 0, \min\{\eta_k, \eta_{\sigma^{-1}(k)}\} \geq \frac{1}{2}\eta_sr|b_k| \geq |\eta_k - \eta_{\sigma^{-1}(k)}|$, for any $k \in [1,n]\}$.
                \begin{enumerate}
                    \item[(A)] If the set $\mathcal{W}$ is empty, then set $\Sigma_s \gets \Sigma_s\setminus\{\sigma\}$ and go to subroutine (b).
                    \item[(B)] Pick a pair $(\mathbf{b},r) \in \mathcal{W}$, if there exists a polynomial $Z(u) = \sum_{i=0}^{r}z_iu^i$ such that $z_0, z_r \ne 0$ and $P(\mathbf{x}) = \mathbf{x}^{r[-\mathbf{b}]_+}Z(\mathbf{x}^{\mathbf{b}})$, then let
                    $\omega_s:=(\mathbf{b},r,Z)$ be a seedlet and
                    set $\mathcal{M}(F) \gets \mathcal{M}(F) \cup \{(\psi_{\sigma, s,\omega_s}, \widetilde{\mathbf{d}})\}$. Otherwise, set $\mathcal{W} \gets \mathcal{W} \setminus\{(\psi_{\sigma, s,\omega_s}, \widetilde{\mathbf{d}})\}$ and go to subroutine (A).
                \end{enumerate}
            \end{enumerate}
        \end{enumerate}
    \end{enumerate}
    This algorithm is complete in finite steps since the set $S$ in subroutine (1), the set $\Sigma_s$ in subroutine (3), and the set $\mathcal{W}$ in subroutine (iv) are finite. This algorithm can obtain the set $\mathcal{M}(F)$ of non-trivial cluster pairs of $F(\mathbf{x})$ by following reasons,
    \begin{itemize}
        \item By Definition \ref{def csm}(ii), the subroutine (B) determines whether the polynomial $P(\mathbf{x})$ found in subroutine (iii) is an exchange polynomial of the seedlet $\omega_s := (\mathbf{b}, r, Z)$. If so, the map $\psi_{\sigma, s,\omega_s}$ is a cluster symmetric map.
        \item In subroutine (iv), if $(\mathbf{b}, r) \notin \mathcal{W}$, then by Proposition \ref{d and n}, we know that $\widetilde{F}(\psi_{\sigma, s, \omega_s}(\mathbf{x})) \ne \widetilde{F}(\mathbf{x})$. Hence $(\mathbf{b}, r)$ should belongs to the set $\mathcal{W}$.
        \item The subroutines (3), (b), (i), (ii) and (iii) determine whether the relations \eqref{len1_f=fp@s}, \eqref{len1_d=sigmad} and \eqref{len1_eta_s=2d_s} hold. If so, then by Theorem \ref{thm len1}, we have $\widetilde{F}(\psi_{\sigma, s, \omega_s}(\mathbf{x})) = \widetilde{F}(\mathbf{x})$, where $\widetilde{F}(\mathbf{x}) := \mathbf{x}^{-\widetilde{\mathbf{d}}} F(\mathbf{x})$. Hence, by Definition \ref{def csp}, we know that the pair $(\psi_{\sigma, s,\omega_s}, \widetilde{\mathbf{d}})$ is a cluster symmetric pair of $F(\mathbf{x})$.
        \item The subroutines (1) and (3) ensure the relation $\eta_s = \eta_{\sigma^{-1}(s)} \ne 0$, then by Proposition \ref{RMK non-trivial}, we know the pair $(\psi_{\sigma, s,\omega_s}, \widetilde{\mathbf{d}})$ is non-trivial. Hence $(\psi_{\sigma, s,\omega_s}, \widetilde{\mathbf{d}}) \in \mathcal{M}(F)$.
    \end{itemize}
\end{alg}

Based on the above algorithm, we write a MATLAB program attached to Appendix \ref{code diop2ca}. Thus, readers can efficiently and conveniently obtain the non-trivial cluster symmetric set $\mathcal{M}(F)$. When the set $\mathcal{M}(F)$ is nonempty, the Laurent polynomial $\widetilde{F}(\mathbf{x})$ in the set $\{ F(\mathbf{x})/\mathbf{x}^{\widetilde{\mathbf{d}}} \mid (\psi_{\sigma, s, \omega_s},\widetilde{\mathbf{d}}) \in \mathcal{M}(F) \}$ has cluster symmetry.

\begin{prop}\label{prop S(F)}
Given a Laurent polynomial $F(\mathbf{x}) \in \mathbb{Q}[\mathbf{x}^{\pm}]$. The cluster symmetric set $\mathcal{S}(F)$ of $F(\mathbf{x})$ can be obtained by Algorithm \ref{algo diop2cs}, that is, $\mathcal{S}(F) = \{ \psi_{\sigma,s,\omega_s} \mid (\psi_{\sigma,s,\omega_s}, \mathbf{0}) \in \mathcal{M}(F) \}$.
\end{prop}
\begin{proof} $(\supset):$ Obviously.
    $(\subset):$ Clearly the pair $(\psi_{\sigma,s,\omega_s}, \mathbf{0})$ is a cluster symmetric pair of $F(\mathbf{x})$. By Theorem \ref{thm len1_T2f}, we know $\eta_{\sigma^{-1}(s)} = \eta_s \ne 0$. By Proposition \ref{RMK non-trivial}, we know that $(\psi_{\sigma,s,\omega_s}, \mathbf{0})$ is non-trivial, and hence it belongs to the set $\mathcal{M}(F)$.
\end{proof}

\begin{eg}\label{eg F to 1-CSpair}
(i) Consider the polynomial $$T_1(\mathbf{x}) := ax_2x_3^2 + x_1^2x_4 + bx_2^2x_4,$$
where $a,b\in \mathbb{Z}_{>0}$. Using Algorithm \ref{algo diop2cs}, or running the corresponding MATLAB program in Appendix \ref{code diop2ca}, we know that the non-trivial pairs of the polynomial $T_1(\mathbf{x})$ are $$(\psi^{\pm}_{\sigma_{(24)}, 1, \omega_1}, \mathbf{d}),\quad (\psi^{\pm}_{\sigma_{(24)}, 1, \omega'_1}, \mathbf{d}),$$
where $\omega_1 := ((0,1,-2,1), 1, a+bu)$, $\omega'_1 := ((0,-1,2,-1), 1, b+au)$, $\mathbf{d} := (1, d_2, d_3, d_2)$ and $d_2,d_3 \in \mathbb{Z}$. 
By Property \ref{pty ex-poly}(i), we know $(\psi_{\sigma_{(24)}, 1, \omega_1}, \mathbf{d}) = (\psi_{\sigma_{(24)}, 1, \omega'_1}, \mathbf{d})$. Hence $\mathcal{M}(T_1) = \{(\psi^{\pm}_{\sigma_{(24)}, 1, \omega_1}, \mathbf{d})\}$. Let $F_1(\mathbf{x}) := \frac{T_1(\mathbf{x})}{\mathbf{x}^{\mathbf{d}}}$. Then we have $\mathcal{G}(F_1) = \langle \psi_{\sigma_{(24)}, 1, \omega_1} \rangle$ and, by Proposition \ref{prop G subset V}(i), we have $\mathcal{G}(F_1)(\mathbf{x}_0) \subset \mathcal{V}_{\mathbb{Q}_{>0}}(F_1(\mathbf{x}) - F_1(\mathbf{x}_0))$ for any tuple ${\mathbf{x}_0} \in \mathbb{Q}_{>0}^4$.

(ii) Consider the polynomial
    $$T_2(\mathbf{x}) := (x_1x_2+ax_3^2+b^2x_4^2)(x_1+x_2)  + bx_4(x_1^2 + x_2^2) + abx_3^2x_4,$$
    where $a,b\in \mathbb{Z}_{>0}$. We know that the non-trivial pairs of the polynomial  $T_2(\mathbf{x})$ are
    \begin{eqnarray*}
        (\psi^{\pm}_{\sigma_{(12)}, 1, \omega_1}, \mathbf{d}_1),~ (\psi^{\pm}_{id, 1, \omega_1}, \mathbf{d}_2),~
        (\psi^{\pm}_{\sigma_{(12)}, 2, \omega_2}, \mathbf{d}_3),~ (\psi^{\pm}_{id, 2, \omega_2}, \mathbf{d}_4),
    \end{eqnarray*}
    where $\omega_1 := ((0,1,-2,1), 1, a+bu)$, $\omega_2 := ((1,0,-2,1), 1, a+bu)$, $\mathbf{d}_1 := (1, 1, d_3, d_4) $, $\mathbf{d}_2 := (1, d_2, d_3, d_4)$, $\mathbf{d}_3 := (1, 1, d_3, d_4) $, $\mathbf{d}_4 := (d_1, 1, d_3, d_4)$ for every
    $d_i \in \mathbb{Z}$. Let $F_2(\mathbf{x}) := \frac{T_2(\mathbf{x})}{\mathbf{x}^{\mathbf{d}_1}}$. 
    Then we have 
    $$\mathcal{G}(F_2) = \langle \psi_{\sigma_{(12)}, 1, \omega_1}, \psi_{id, 1, \omega_1}, \psi_{\sigma_{(12)}, 2, \omega_2}, \psi_{id, 2, \omega_2} \rangle$$
    and $\mathcal{G}(F_2)(\mathbf{x}_0) \subset \mathcal{V}_{\mathbb{Q}_{>0}}(F_2(\mathbf{x}) - F_2(\mathbf{x}_0))$ for any tuple ${\mathbf{x}_0} \in \mathbb{Q}_{>0}^4$.

    (iii) Consider a polynomial $$T_3(\mathbf{x}) := x_1^2x_4^2 + x_2^2x_3^2 + x_1x_3^3 + x_2^3x_4.$$ 
    We know the non-trivial cluster symmetric pairs of the polynomial $T_3(\mathbf{x})$ are
    \begin{eqnarray*}
        (\psi^{\pm}_{\sigma_{(1234)}, 1, \omega_1}, \mathbf{d}),~
        (\psi^{\pm}_{\sigma_{(24)}, 1, \omega_1}, \mathbf{d}_1),~ (\psi^{\pm}_{\sigma_{(13)}, 4, \omega_4}, \mathbf{d}_2),
    \end{eqnarray*}
    where $\omega_1 := ((0,1,-2,1), 1, 1+u)$, $\omega_4 := ((-1,2,-1,0), 1, 1+u)$, $\mathbf{d} := (1,1,1,1)$, $\mathbf{d}_1 := (1, d_2, d_3, d_2) $, $\mathbf{d}_2 := (d_1, d_2,d_1,1)$ and
    $d_i\in \mathbb{Z}$. 
    Let $F_3(\mathbf{x}) := \frac{T_3(\mathbf{x})}{\mathbf{x}^{\mathbf{d}}}$.
    Then we have 
    $$\mathcal{G}(F_3) = \langle \psi_{\sigma_{(1234)}, 1, \omega_1}, \psi_{\sigma_{(24)}, 1, \omega_1}, \psi_{\sigma_{(13)}, 4, \omega_4} \rangle$$ and $\mathcal{G}(F_3)(\mathbf{x}_0) \subset \mathcal{V}_{\mathbb{Q}_{>0}}(F_3(\mathbf{x}) - F_3(\mathbf{x}_0))$ for any tuple ${\mathbf{x}_0} \in \mathbb{Q}_{>0}^4$.
\end{eg}

\subsection{Generalized cluster algebra associated to a Laurent polynomial}
Based on the results of the previous subsection, we can further determine whether a given Laurent polynomial can be realized within a generalized cluster algebra. This approach enables us to leverage the positive Laurent phenomenon of the generalized cluster algebra (Theorem \ref{thm PLP}), as discussed in Section \ref{sec diop}, to solve Diophantine equations $F(\mathbf{x}) = F(\mathbf{1})$. To this end, inspired by Proposition \ref{prop G subset V}(ii), we give the following definition.

\begin{DEF} \label{def 1-cs seed}
    Given a Laurent polynomial $F(\mathbf{x}) \in \mathbb{Q}[\mathbf{x}^{\pm}]$ with nonempty cluster symmetric set $\mathcal{S}(F)$.
    If there exists a seed $\Omega$, such that the cluster symmetric set $\mathcal{S}(F)$ of $F(\mathbf{x})$ is a subset of the cluster symmetric set $\mathcal{S}(\Omega)$ of $\Omega$, then we call the seed $\Omega$ a \textbf{cluster symmetric seed} of $F(\mathbf{x})$ and the generalized cluster algebra $\mathcal{A}(\Omega)$ a \textbf{generalized cluster algebra associated to} $F(\mathbf{x})$. 
\end{DEF}
In this situation, the Laurent polynomial $F(\mathbf x)$ is a cluster symmetric polynomial about any cluster symmetric map in $\mathcal{S}(F)$. And, by Proposition \ref{prop G subset V}(ii), we have $$\mathcal{G}(F)(\mathbf{1}) \subset \mathcal{V}_{\mathbb{Z}_{>0}}(F(\mathbf{x}) - F(\mathbf{1})).$$

By Definition \ref{def seed and 1-csm}, we only need to check whether all cluster symmetric maps in the set $\mathcal{S}(F)$ correspond to the same seed, which can determine the cluster symmetric seed of $F(\mathbf{x})$. 

\begin{prop} \label{prop find 1-csSeed}
    Given a Laurent polynomial $F(\mathbf{x}) \in \mathbb{Q}[\mathbf{x}^{\pm}]$ with nonempty cluster symmetric set $\mathcal{S}(F)$. If there exists a seed $\Omega$, such that it corresponds to any cluster symmetric map in the set $\mathcal{S}(F)$, then the seed $\Omega$ is a cluster symmetric seed of $F(\mathbf{x})$, that is, $\mathcal{S}(F) \subseteq \mathcal{S}(\Omega)$. Followingly, $\mathcal{G}(F) $ is a subgroup of $\mathcal{G}(\Omega)$.
\end{prop}
\begin{proof}
    Let $\psi_{\sigma,s,\omega_s} \in \mathcal{S}(F)$. Since $\Omega$ corresponds to $\psi_{\sigma,s,\omega_s}$, by Definition \ref{def seed and 1-csm}, we know that $\psi_{\sigma,s,\omega_s} \in \mathcal{S}(\Omega)$. Hence $\mathcal{S}(F) \subseteq \mathcal{S}(\Omega)$.
\end{proof}

\begin{eg} (i) Consider the Laurent polynomial 
    $$F_1(\mathbf{x}) := \frac{ax_2x_3^2 + x_1^2x_4 + bx_2^2x_4}{x_1x_2x_3x_4}.$$
By Example \ref{eg F to 1-CSpair}(i), we know the cluster symmetric group $\mathcal{G}(F_1) = \langle \psi_{\sigma_{(24)}, 1, \omega_1} \rangle$. When $(a,b) \ne (1,1)$, by Example \ref{eg 1-csm not to seed}, the cluster symmetric map $\psi_{\sigma,1,\omega_1}$ does not correspond to any seeds. When $(a,b) = (1,1)$, by Example \ref{eg 1-csm to seed}(i),
we know that $\psi_{\sigma_{(24)}, 1, \omega_1}$ corresponds to the seed $\Omega := (B, \mathbf{x}, R, \mathbf{Z})$ where
    \begin{equation*}
        B = \begin{bmatrix}
        		0 & -1 & 2 & -1\\
			    1 &  0 & -c & -d\\
			   -2 &  c & 0 & 2-c\\
                1 &  d & c-2 &  0
			\end{bmatrix},  R = \begin{bmatrix}
        		1 &  &  & \\
			    & r_2 & &  \\
			    &  & r_3 & \\
                &  &  &  r_2  
			\end{bmatrix},  \begin{array}{l}
        Z_1(u) = 1 + u, \\
        Z_2(u) = \sum_{i=0}^{r_2}z_{2,i}u^i, \\
        Z_3(u) = \sum_{i=0}^{r_3}z_{3,i}u^i, \\
        Z_4(u) = Z_2(u), \end{array}
    \end{equation*}
    where $c,d \in \mathbb{Z}$. Hence, the seed $\Omega$ is a cluster symmetric seed of $F(\mathbf{x})$.

(ii) Consider the Laurent polynomial
    $$F_2(\mathbf{x}) := \frac{(x_1x_2 + x_3^2 + x_4^2)(x_1+x_2)  + x_4(x_1^2 + x_2^2) + x_3^2x_4}{x_1x_2x_3x_4}.$$
    By Example \ref{eg F to 1-CSpair}(ii), we know the cluster symmetric group of $F(\mathbf{x})$ is
    $$\mathcal{G}(F_2) = \langle \psi_{\sigma_{(12)}, 1, \omega_1}, \psi_{id, 1, \omega_1}, \psi_{\sigma_{(12)}, 2, \omega_2}, \psi_{id, 2, \omega_2} \rangle.$$
 By Example \ref{eg 1-csm to seed}(iii), the cluster symmetric map $\psi_{\sigma_{id},1,\omega_1}$ corresponds to the seed $\Omega := (B, \mathbf{x}, R, \mathbf{Z})$, where
    $$B = \begin{bmatrix}
        		0 & -1 & 2 & -1\\
			    1 &  0 & -1 & 0\\
			    -2 & 1 & 0 & 1\\
                1 & 0 & -1 &  0
			\end{bmatrix},  R = \begin{bmatrix}
        		1 &  &  & \\
			    & r_2 & &  \\
			    &  & r_3 & \\
                 &  &  &  r_4  
			\end{bmatrix}, \begin{array}{l}
        Z_1(u) = 1 + u, \\
        Z_2(u) = \sum_{i=0}^{r_2}z_{2,i}u^i, \\
        Z_3(u) = \sum_{i=0}^{r_3}z_{3,i}u^i, \\
        Z_4(u) = \sum_{i=0}^{r_4}z_{4,i}u^i. \end{array}$$
But, by Example \ref{eg 1-csm to seed}(ii), the cluster symmetric map $\psi_{\sigma,1,\omega_1}$ does not correspond to any seeds.
Hence, there is no cluster symmetric seed of $F(\mathbf{x})$.

(iii) Consider the Laurent polynomial 
$$F_3(\mathbf{x}) := \frac{x_1^2x_4^2 + x_2^2x_3^2 + x_1x_3^3 + x_2^3x_4}{x_1x_2x_3x_4}.$$ 
By Example \ref{eg F to 1-CSpair}(iii), we know the cluster symmetric group of $F(\mathbf{x})$ is
$$\mathcal{G}(F_3) = \langle \psi_{\sigma_{(1234)}, 1, \omega_1}, \psi_{\sigma_{(24)}, 1, \omega_1}, \psi_{\sigma_{(13)}, 4, \omega_4} \rangle.$$
By Example \ref{eg 1-csm to seed}(iv), we know that $\psi_{\sigma_{(1234)}, 1, \omega_1}, \psi_{\sigma_{(24)}, 1, \omega_1}, \psi_{\sigma_{(13)}, 4, \omega_4}$ corresponds to the seed $\Omega := (B, \mathbf{x}, R, \mathbf{Z})$, where
\begin{equation*}
    B = \begin{bmatrix}
            0 & -1 & 2 & -1\\
            1 &  0 & -3 & 2\\
           -2 &  3 & 0 & -1\\
            1 &  -2 & 1 &  0
        \end{bmatrix}, ~ R = \begin{bmatrix}
            1 &  &  & \\
            & 1 & &  \\
            &  & 1 & \\
            &  &  &  1  
        \end{bmatrix}, ~ \begin{array}{l}
    Z_1(u) = 1 + u, \\
    Z_2(u) = 1 + u, \\
    Z_3(u) = 1 + u, \\
    Z_4(u) = 1 + u. \end{array}
\end{equation*}
Hence, the seed $\Omega$ is a cluster symmetric seed of $F_3(\mathbf{x})$. In fact, the cluster symmetric map $\psi_{\sigma_{(1234)}, 1, \omega_1}$ is related to the Somos $4$ sequence in \cite{HoneSwart}.
In their paper, Hone and Swart constructed the Laurent polynomial $F_3(\mathbf{x})$ which remains invariant under the action of the cluster symmetric map $\psi_1$, and also proved that the Somos $4$ sequence is related to an elliptic curve. 
\end{eg}

By Proposition \ref{prop find 1-csSeed}, we have $\mathcal{S}(F) \subseteq \mathcal{S}(\Omega)$. 
Finally, we investigate the reverse inclusion relationship between $\mathcal{S}(\Omega)$ and $\mathcal{S}(F)$ from the perspective of the generalized cluster algebra.

\begin{prop}
    Given a seed $\Omega$ with nonempty cluster symmetric set $\mathcal{S}(\Omega)$. 
    Let $F(\mathbf{x})$ be a Laurent polynomial of type $\frac{\bm{\eta}}{\mathbf{d}}$ that is invariant under the cluster symmetric group $\mathcal{G}(\Omega)$. 
    If $\eta_s \ne 0$ for all $\sigma\mu_s \in \mathcal{S}(\Omega)$, then $S(\Omega) \subseteq S(F)$. Followingly, $\mathcal{G}(\Omega)$ is a subgroup of $\mathcal{G}(F)$.
\end{prop}
\begin{proof}
    Let $\sigma\mu_s \in \mathcal{S}(\Omega)$. Since $F(\sigma\mu_s(\mathbf{x})) = F(\mathbf{x})$ and $\eta_s \ne 0$, by Proposition \ref{is CSM}, we know that $\sigma\mu_s \in \mathcal{S}(F)$.    
\end{proof}

To summarize this subsection, find a cluster symmetric seed of a given Laurent polynomial $F(\mathbf{x})$ in the following steps:
\begin{enumerate} 
    \item Using Algorithm \ref{algo diop2cs} or MATLAB program in Appendix \ref{code diop2ca}, we obtain the set $\mathcal{M}(F)$ of non-trivial cluster symmetric pairs of $F(\mathbf{x})$;
    \item By Definition \ref{def S1 of F}, construct the cluster symmetric set $\mathcal{S}(F)$ from the set $\mathcal{M}(F)$;
    \item When the set $\mathcal{S}(F)$ is nonempty, find a seed $\Omega$, such that $\psi_{\sigma,s,\omega_s}$ is corresponds to the seed $\Omega$ for any $\psi_{\sigma,s,\omega_s} \in \mathcal{S}(F)$. If one can find, by Proposition \ref{prop find 1-csSeed}, the seed $\Omega$ is a cluster symmetric seed of $F(\mathbf{x})$.
\end{enumerate}

\subsection{Summary} 
\begin{figure}
\[\begin{tikzcd} 
	\begin{array}{c} \text{Definition } \ref{def geneCA}:\\ \textbf{generalized cluster algebra }\mathcal{A}(\Omega) \end{array}\\
	\begin{array}{c} \text{Definition } \ref{def seed}:\\ \textbf{seed } \Omega \end{array}\\
	\\
    \\
    \\
    \\
	\boxed{\begin{array}{c} \text{Definition } \ref{seedlet}: \\ \textbf{cluster symmetric map } \psi_{\sigma,s,\omega_s} \end{array}} & \begin{array}{c} \text{Definition }\ref{def S1 of F}:\\ \textbf{cluster symmetric set } \mathcal{S}(F)
    \end{array} 
    \\
	& \begin{array}{c} \text{Definition }\ref{def csp}: \\ \textbf{cluster symmetric pair }\\ (\psi_{\sigma,s,\omega_s},\widetilde{\mathbf{d}}) \text{ of } F(\mathbf{x}). \end{array} \\
	\\
	\\
    \\
	\\
	\\
	{\begin{array}{c} \text{Definition } \ref{def 1-cspoly}: \\ 
    \textbf{cluster symmetric polynomial }\\ F(\mathbf{x}) \in \mathbb{Q}[\mathbf{x}^{\pm}] \end{array}} & \begin{array}{c}
	     \text{Describing invariant rings } \mathbb{Q}[\mathbf{x}^{\pm}]^{G}\\
         \begin{cases}
	   \text{ Proposition }\ref{F in ring: id} \\
       \text{ Proposition } \ref{prop len0 B=0}\\
       \text{ Theorem } \ref{thm 4bounded}(ii)
	   \end{cases}     
\end{array}  \\
	\begin{array}{c}\textbf{cluster symmetric equation } \\ F(\mathbf{x}) = c \end{array} &  \begin{array}{c}
	     \text{Solving equations}\\
	   \left\{\begin{array}{c} 
       \text{Theorem }  \ref{thm markov-cluster sol}:\\
       \textbf{Markov-cluster equations }
	\end{array}\right.
    \end{array}
	\arrow[tail reversed, from=1-1, to=2-1] 
	\arrow["\begin{array}{c} \text{Definition }\ref{def csm}:\\ \textbf{cluster symmetric set } \mathcal{S}(\Omega)\\ 
    ~\\ \text{Proposition }\ref{is CSM}:\\
    \text{Any elements in } \mathcal{S}(\Omega) \text{ are }\\
    \text{cluster symmetric maps}.
    \end{array}"', curve={height=12pt}, from=2-1, to=7-1]
	\arrow["\begin{array}{c} \text{Definition } \ref{def seed and 1-csm}:\\
    \psi_{\sigma, s, \omega_s} \textbf{corresponds to } \Omega\text{,}\\
    \text{if } \sigma\mu_s \in \mathcal{G}(\Omega)\text{, } \omega_s = \pi_s(\Omega^{\pm}).
    \end{array}"', curve={height=12pt}, from=7-1, to=2-1]
	\arrow[from=7-2, to=7-1]
	\arrow["\begin{array}{c} 
    \text{Theorem } \ref{thm len1}:\\ 
    \text{Construct }F(\mathbf{x})\in \mathbb{Q}\lbrack\mathbf{x}^{\pm}\rbrack^{\langle\psi_{\sigma,s,\omega_s}\rangle}.\\
    \text{Appendix }\ref{code ca2diop}:\\ \text{Corresponding MATLAB program.}\\
    \text{Remark } \ref{rmk step of constuct}:\\
    \text{Concrete construction steps.}
    \end{array}"', from=7-1, to=14-1]
    \arrow[curve={height=12pt}, from=8-2, to=7-2]
	\arrow["\begin{array}{c} \text{Algorithm }\ref{algo diop2cs}:\\ \text{Find the set $\mathcal{M}(F)$ of }\\
    \text{ non-trivial cluster symmetric pairs of } F(\mathbf{x}).\\
    \text{Appendix }\ref{code diop2ca}: \\ \text{Corresponding MATLAB program.}
    \end{array}"'{pos=0.6}, curve={height=-40pt}, from=14-1, to=8-2]
	\arrow[from=14-1, to=15-1]
    \arrow[dashed, from=14-1, to=14-2]
    \arrow[dashed, from=15-1, to=15-2]
    \arrow["\begin{array}{c}
	     \text{Definition } \ref{def 1-cs seed}:  \\
         \textbf{cluster symmetric seed of } F(\mathbf{x})\\
         ~~\\
         \text{Proposition } \ref{prop find 1-csSeed}:\\
         \text{Find a seed } \Omega  \text{, such that }\\
         \mathcal{S}(F) \subset \mathcal{S}(\Omega)
	   \end{array}"'{pos=0.5}, curve={height=40pt}, from=7-2, to=2-1]
\end{tikzcd}\]
\caption{Relationship between cluster theory and Diophantine equations}
\label{fig: CAN}
\end{figure}
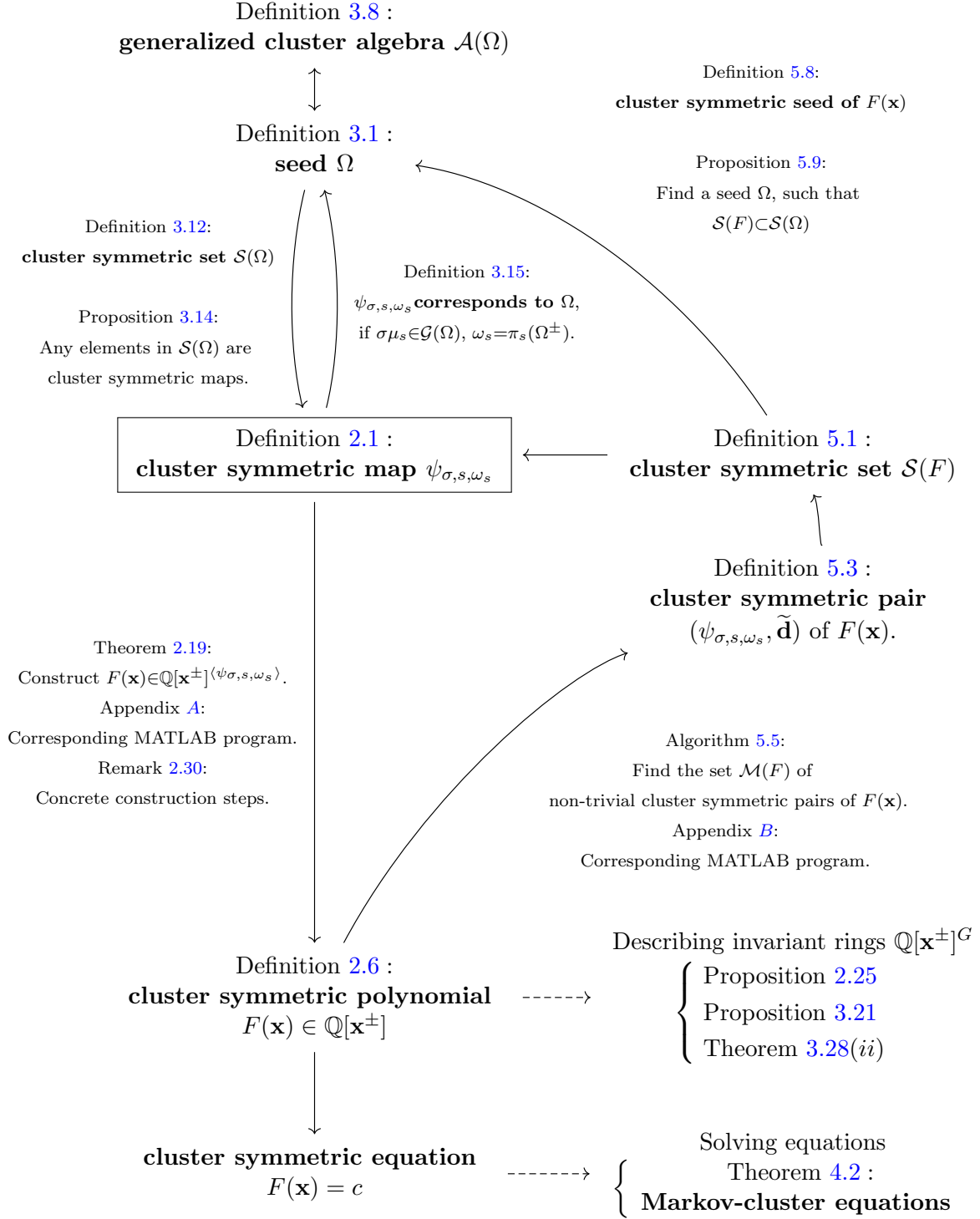

As a summary, we describe the relationship between the main notations and the main results of this paper in Figure \ref{fig: CAN}. Through cluster symmetry, we establish a connection between cluster theory and Diophantine equations. Additionally, as an application, we describe three classes of invariant rings and solve two Diophantine equations.

Number theory, first systematically investigated by Diophantus, has spanned nearly two millennia; invariant theory, established by Hilbert, has evolved for over a century; whereas cluster theory, pioneered by Fomin and Zelevinsky, has been developed for more than two decades.
Each of these three disciplines has yielded abundant achievements and profound insights within their respective domains. 
An in-depth study of cluster-theoretic approaches to Diophantine equations contributes to a deeper understanding of the intrinsic connections between number theory, invariant theory, and cluster theory. Our future work will continue to focus on this intersection.

\newpage
\appendix
\section{MATLAB program of Theorem \ref{thm len1}}\label{code ca2diop}
The associated MATLAB program can be downloaded at this \href{https://mathbao.github.io}{link}. All programs can be run on \href{https://matlab.mathworks.com}{MATLAB Online}. Here we show some examples.

We first consider Example \ref{somos5} for $\tilde{\alpha} = 5, \tilde{\beta} = 3$. In the command line window, enter the following code.

\begin{figure}[htbp]
    \centering
    \includegraphics[width=1\linewidth]{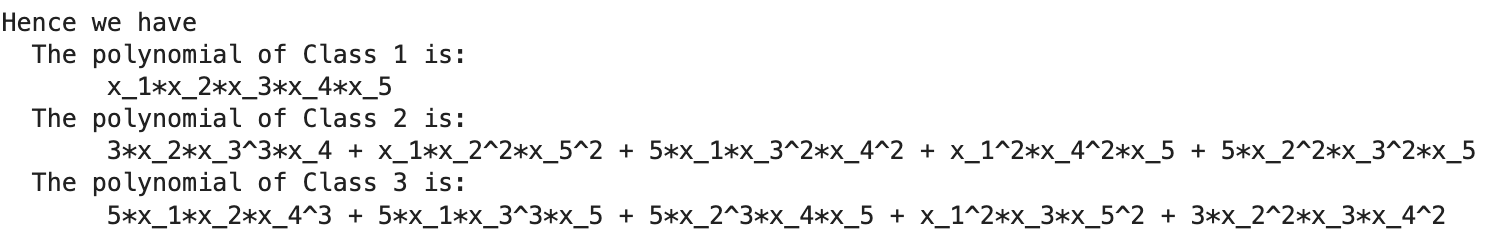}
    \caption{Result of Code \ref{code: somos5}}
    \label{fig: somos5}
\end{figure}

\begin{matlab}{Find cluster symmetric polynomials of Example \ref{somos5}}{code: somos5}
%% Input Data
%% the seedlet $\omega_s = (b, r, Z)$
b = [0,1,-1,-1,1];  % tuple $\mathbf{b}$
r = 1;               
Z = [3, 5]; % coefficients of the polynomial Z(u)=3+5u

s = 1;  % direction $s$
sigma = [2,3,4,5,1]; % permuation $\sigma_{(12345)}$, where sigma(i) is $\sigma_{(12345)}(i)$
eta = [2,3,4,3,2]; % $\bm\eta$
d = [1,1,1,1,1];   % $\mathbf{d}$

%% Solve the system of homogeneous linear equation $HLE(\sigma, s, \omega_s, \bm{\eta}, \mathbf{d})$
FindTheLaurentPolyOf(b,r,Z,s,sigma,eta,d);
\end{matlab}
After 34.32 seconds of computation, we get the result shown in Figure \ref{fig: somos5}. 

From this result, we obtain a monomial $x_1x_2x_3x_4x_5$ and the following two polynomials
\begin{eqnarray*}
    T_1(\mathbf{x}) := x_1x_2^2x_5^2 + x_1^2x_4^2x_5 + 5(x_1x_3^2x_4^2  + x_2^2x_3^2x_5) + 3x_2x_3^3x_4,\\
    T_2(\mathbf{x}) := x_1^2x_3x_5^2 + 5(x_1x_2x_4^3 + x_1x_3^3x_5 + x_2^3x_4x_5) + 3x_2^2x_3x_4^2.
\end{eqnarray*}
We let $F_1(\mathbf{x}) := \frac{T_1(\mathbf{x})}{x_1x_2x_3x_4x_5}$ and $F_2(\mathbf{x}) := \frac{T_2(\mathbf{x})}{x_1x_2x_3x_4x_5}$. Then the Laurent polynomial $a_1F_1(\mathbf{x}) + a_2F_2(\mathbf{x}) + a_3$ is invariant under the cluster symmetric map $\psi_{\sigma_{(12345)},1,\omega_1}$.

We then consider Question \ref{que: gyoda1}. In the command line window, enter the following code.

\begin{matlab}{cluster symmetric polynomials of Question \ref{que: gyoda1}}{code: gyoda}
%% Input the seed $\Omega = (B, R, \mathbf{Z})$
B = [0, 1,-1;
    -1, 0, 2;
     1,-2, 0];
R = [4,1,1];
syms k1 k2
Z = [1, k1,k2,k1,1;
     1, 1, 0, 0, 0;
     1, 1, 0, 0, 0];

S = [1;2;3]; % direction list S
Sigma = [1:3; 1:3; 1:3]; % permuation list Sigma
eta = [2,4,4]; % $\bm\eta$
d = [1,2,2];   % $\mathbf{d}$

%% Solve the systems of homogeneous linear equation $HLE(\sigma, s, \omega_s, \bm{\eta}, \mathbf{d})$ for all $s \in S$
FindTheLaurentPolyOf(B,R,Z,S,Sigma,eta,d);
\end{matlab}

\begin{figure}[htbp]
    \centering
    \includegraphics[width=1\linewidth]{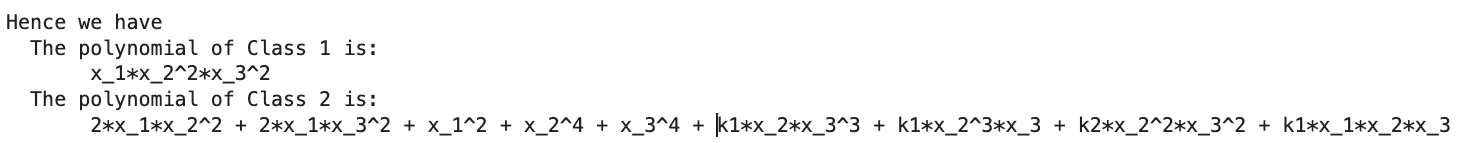}
    \caption{Result of Code \ref{code: gyoda}}
    \label{fig: gyoda}
\end{figure}

After 256.59 seconds of computation, we get the result of the MATLAB program shown in Figure \ref{fig: gyoda}.

From this result, we obtain a monomial $x_1x_2^2x_3^2$ and the following polynomial
\begin{eqnarray*}
        T_{3,7}(\mathbf{x}) := x_1^2 + x_2^4 + x_3^4 + 2x_1x_2^2 + 2x_1x_3^2  + k_1x_2x_3^3 +  k_2x_2^2x_3^2 + k_1x_2^3x_3 + k_1x_1x_2x_3.
\end{eqnarray*}
We let $F_{3,7}(\mathbf{x}) := a\frac{T_{3,7}(\mathbf{x})}{x_1x_2^2x_3^2} + b$, where $a,b \in \mathbb{Q}$. Then the Laurent polynomial $F_{3,7}(\mathbf{x})$ is invariant under mutations $\mu_1, \mu_2, \mu_3$, where $\mu_i \in \mathcal{S}(\Omega_{3,7})$. This provides an affirmative answer to Question \ref{que: gyoda1}, that is, the seed $\Omega_{3,7}$ has a corresponding Diophantine equation $F_{3,7}(\mathbf{x}) = F_{3,7}(\mathbf{x}_0)$ for some $\mathbf{x}_0 \in \mathbb{Q}^3_{>0}$.

\section{MATLAB program of Algorithm \ref{algo diop2cs}}
\label{code diop2ca}

The associated MATLAB program can be downloaded at this \href{https://mathbao.github.io}{link}. All programs can be run on \href{https://matlab.mathworks.com}{MATLAB Online}. Here we show some examples.

We first consider Example \ref{eg F to 1-CSpair}(i). In the command line window, enter the following code.

\begin{matlab}{Find all non-trivial cluster symmetric pairs of Example \ref{eg F to 1-CSpair}(i)}{code: pair1}
n = 4;                 % rank n
x = sym('x_', [1,n]);  % varibles \mathbf{x}
syms alpha beta        % constants
Tpower = [0,1,2,0;
          2,0,0,1;
          0,2,0,1];
Tcoeff = [alpha,1,beta];
T = Tcoeff*prod(x.^Tpower,2); % poly $T(\mahtbf{x})$
d = zeros(1,n);   % $\mathbf{d}$

%%  Find the set of non-trivial cluster symmetric pair of $x^{-\mathbf{d}}T(\mathbf{x})$.
M = FindTheClusterSymPairOf(T,d,x);  % M is the set
\end{matlab}

\begin{figure}[htbp]
\begin{center}
\includegraphics[width=1\linewidth]{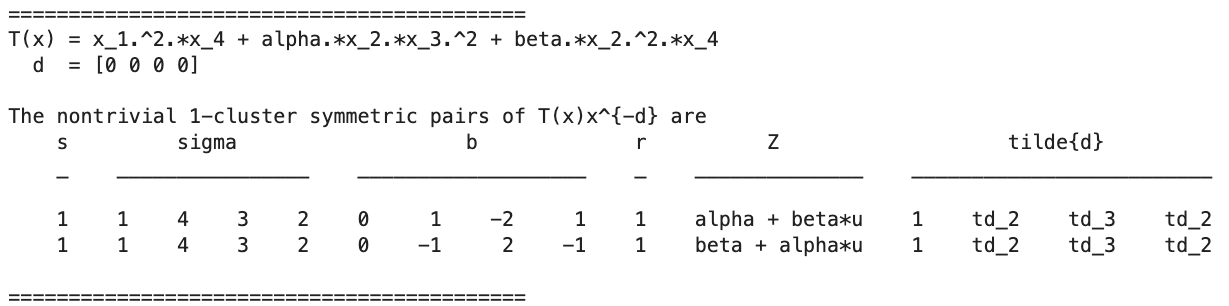}
\caption{Result of Code \ref{code: pair1}.\label{fig: pair1}}
\end{center}
\end{figure}

After 1.31 seconds of computation, we get the result shown in Figure \ref{fig: pair1}.

Each row of the table in Figure \ref{fig: pair1} is a non-trivial cluster symmetric pair. For example, the first row corresponds to the non-trivial cluster symmetric pair $(\psi_{\sigma_{(1234)}, 1, \omega_1},\mathbf{d})$
where $\omega_1 := ((0,1,-2,1), 1, \alpha+\beta u)$, $\mathbf{d} := (1, d_2, d_3, d_2)$ and $d_2,d_3 \in \mathbb{Z}$.

We then consider Example \ref{eg F to 1-CSpair}(ii). In the command line window, enter the following code.

\begin{matlab}{Find all non-trivial cluster symmetric pairs of Example \ref{eg F to 1-CSpair}(ii)}{code: pair2}
n = 4;                 % rank n
x = sym('x_', [1,n]);  % varibles \mathbf{x}
syms a b        % constants
Tpower = [1,2,0,0;
          2,1,0,0;
          1,0,2,0;
          1,0,0,2;
          0,1,2,0;
          2,0,0,1;
          0,1,0,2;
          0,2,0,1;
          0,0,2,1];
Tcoeff = [1,1,a,b^2,a,b,b^2,b,a*b];
T = Tcoeff*prod(x.^Tpower,2); % poly $T(\mahtbf{x})$
d = zeros(1,n);   % $\mathbf{d}$

%%  Find the set of non-trivial cluster symmetric pair of $x^{-\mathbf{d}}T(\mathbf{x})$.
M = FindTheClusterSymPairOf(T,d,x);  % M is the set
\end{matlab}

\begin{figure}[htbp]
\begin{center}
\includegraphics[width=1\linewidth]{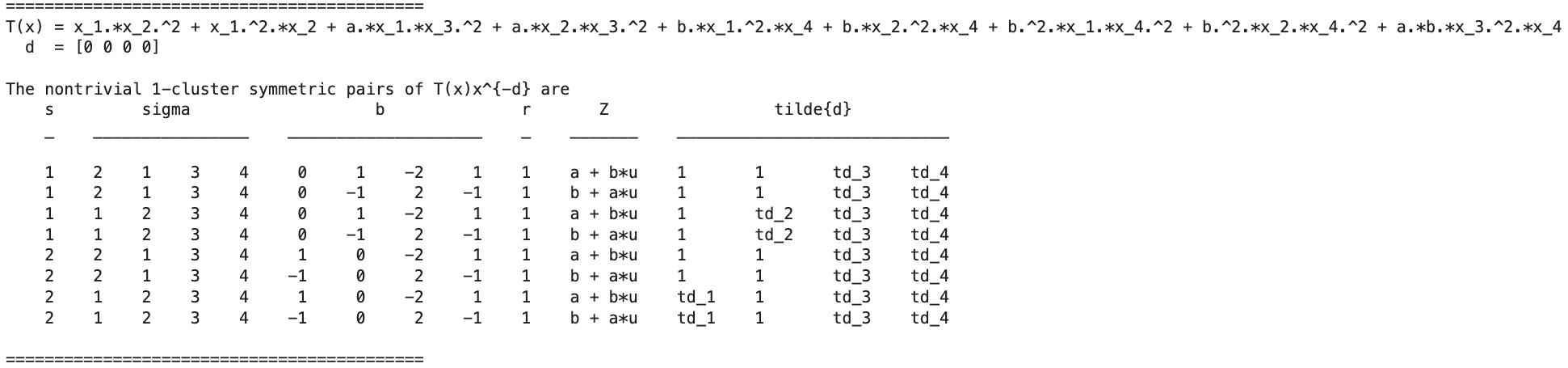}
\caption{Result of Code \ref{code: pair2}.\label{fig: pair2}}
\end{center}
\end{figure}

After 1.47 seconds of computation, we get the result shown in Figure \ref{fig: pair2}.

Each row of the table in Figure \ref{fig: pair1} is a non-trivial cluster symmetric pair which is already shown in Example \ref{eg F to 1-CSpair}(ii).
\vspace{4mm}

\begin{acknowledgements}
    This project is supported by the National Natural Science Foundation of China (No.12131015).
\end{acknowledgements}

\bibliographystyle{amsalpha}
\providecommand{\bysame}{\leavevmode\hbox to3em{\hrulefill}\thinspace}
\providecommand{\MR}{\relax\ifhmode\unskip\space\fi MR }
\providecommand{\MRhref}[2]{%
  \href{http://www.ams.org/mathscinet-getitem?mr=#1}{#2}
}
\providecommand{\href}[2]{#2}

\end{document}